\documentclass[10pt]{article}

\usepackage{changepage}

\usepackage[utf8x]{inputenc}

\usepackage{textcomp,marvosym}

\usepackage{fixltx2e}

\usepackage{amsmath,amssymb}

\usepackage{cite}

\usepackage{nameref,hyperref}

\usepackage[right]{lineno}

\usepackage{microtype}
\DisableLigatures[f]{encoding = *, family = * }

\usepackage{graphicx}

\usepackage{setspace} 

\usepackage{setspace}

\usepackage[labelfont=bf,labelsep=period,justification=raggedright]{caption}
\date{}
\usepackage[aboveskip=1pt,labelfont=bf,labelsep=period,justification=raggedright,singlelinecheck=off]{caption}

\makeatletter
\renewcommand{\@biblabel}[1]{\quad#1.}
\makeatother

\date{}

\newtheorem{defn}{Definition}

\newtheorem{rem}{Remark}

\usepackage{siunitx}

 \newcommand{\hl}[1]{#1}

\usepackage{enumitem}
\setlist[itemize]{leftmargin=*}

\usepackage{epigraph}

\begin{document}
\vspace*{0.2in}

\begin{flushleft}
{\Large
\textbf\newline{Extending Topological Surgery to Natural Processes and Dynamical Systems} 
}
\newline
\\
Stathis Antoniou\textsuperscript{1*},
Sofia Lambropoulou\textsuperscript{2}
\\
\bigskip
\textbf{1,2} Department of Mathematics, National Technical University of Athens, Athens, Greece.
\bigskip

%
%





* santoniou@math.ntua.gr

\epigraph{\hl{\textit{People who wish to analyze nature without using mathematics must settle for a reduced understanding.}}}{\hl{Richard P. Feynman}}

\end{flushleft}
\section*{Abstract}
Topological surgery is a mathematical technique used for creating new manifolds out of known ones. We observe that it occurs in natural phenomena where a sphere of dimension 0 or 1 is selected, forces are applied and the manifold in which they occur changes type. For example, 1-dimensional surgery happens during chromosomal crossover, DNA recombination and when cosmic magnetic lines reconnect, while 2-dimensional surgery happens in the formation of tornadoes, in the phenomenon of Falaco solitons, in drop coalescence and in the cell mitosis. Inspired by such phenomena, we introduce new theoretical concepts which enhance topological surgery with the observed forces and dynamics. To do this, we first extend the formal definition to a continuous process caused by local forces. Next, for modeling phenomena which do not happen on arcs or surfaces but are 2-dimensional or 3-dimensional, we fill in the interior space by defining the notion of solid topological surgery. We further introduce the notion of embedded surgery in $S^3$ for modeling phenomena which involve more intrinsically the ambient space, such as the appearance of knotting in DNA and phenomena where the causes and effect of the process lies beyond the initial manifold, such as the formation of black holes. Finally, we connect these new theoretical concepts with a dynamical system and we present it as a model for both 2-dimensional 0-surgery and natural phenomena exhibiting a `hole drilling' behavior. We hope that through this study, topology and dynamics of many natural phenomena, as well as topological surgery itself, will be better understood.



\let\thefootnote\relax\footnotetext{

{\noindent}\textit{2010 Mathematics Subject Classification}: 57R65, 57N12, 57M25, 57M99, 37B99, 92B99.

{\noindent}\textit{Keywords}: three-space, ambient space, three-sphere, embedded topological surgery, natural phenomena, natural processes, dynamics, continuous process, attracting, forces, joining thread, topological drilling, recoupling, reconnection, mathematical model, Falaco solitons, tornadoes, whirls, DNA recombination, magnetic reconnection, mitosis, meiosis, chromosomal crossover, solid topological surgery, black holes, necking,fracture, gene transfer, dynamical system, bifurcation, Lotka-Volterra, predator-prey, numerical solutions, magnetic field, trajectories, topology.}

\newpage
{\small\tableofcontents}

\newpage
\doublespacing
\section{Introduction} \label{Intro}
\hl{Topological surgery is a mathematical technique used for changing the homeomorphism type, or simply the shape, of a manifold, which is a} `\hl{nice}' \hl{topological space}. This technique creates new manifolds out of known ones. For example, all orientable surfaces may arise from the 2-dimensional sphere using surgery. \hl{Topological surgery can happen in any dimensions, but they all share the same features.} An \textit{ $n$-dimensional topological surgery} on an $n$-manifold $M$ is, roughly, the topological procedure whereby an appropriate $n$-manifold with boundary is removed from  $M$ and  is replaced  by another $n$-manifold with the same boundary, using a `gluing' homeomorphism along the common boundary, thus creating a new  $n$-manifold $\chi(M)$. \hl{The mathematical notions needed for understanding the definition of surgery can be found in} Section~\ref{Useful}. References to illustrations of examples of $1$- and $2$-dimensional surgery can be found in Sections~\ref{1D_FormalE}, ~\ref{2D_FormalE0} and ~\ref{2D_FormalE1}.

\smallbreak 
In this paper we observe that topological surgery is exhibited in nature in numerous, diverse processes of various scales for ensuring new results. Surgery in nature is usually performed on basic manifolds with or without boundary, that undergo merging and recoupling. Such processes are initiated by attracting forces acting on a sphere of dimension 0 (that is, two points) or 1 (that is, a circle). A large part of this work is dedicated to setting the topological ground for modeling such phenomena in dimensions 1,2 and 3. Namely, we introduce new theoretical concepts which are better adapted to the phenomena and which enhance the formal definition of surgery. This work extends significantly the preliminary results and early ideas presented in \cite{An},\cite{SS} and \cite{SSNI}. \hl{With our enhanced definitions of topological surgery in hand, we pin down several physical phenomena undergoing surgery. Furthermore, we present a dynamical system that performs a specific type of surgery.} More precisely, the new concepts are:

\smallbreak
{\noindent} $\bullet$  {\bf The introduction of forces:} A sphere of dimension 0 or 1 is selected in space and attracting forces act on it. These dynamics explain the intermediate steps of the formal definition of surgery and extend it to a continuous process caused by local forces. Note that these intermediate steps can also be explained by Morse theory \hl{but this approach does not involve the forces. On the other hand,} the theoretical forces that we introduce are also observed in the phenomena exhibiting surgery. For example, in dimension 1, during chromosomal crossover the pairing is caused by mutual attraction of the parts of the chromosomes that are similar or homologous, as detailed and  illustrated in Section~\ref{1DDynamics}. In dimension 2, the creation of tornadoes is caused by attracting forces between the cloud and the earth (as detailed and illustrated in Section~\ref{2D0}), while soap bubble splitting is caused by the surface tension of each bubble which acts as an attracting force (this is discussed and illustrated in Section~\ref{2D}).

{\noindent} $\bullet$  {\bf Solid surgery:} The interior of the initial manifold is now filled in. For example, in dimension 1 this allows to model phenomena happening on surfaces such as the merging of oil slicks. An oil slick is seen as a disc, which is a continuum of concentric circles together with the center. An example in dimension 2 is the process of mitosis, whereby a cell splits into two new cells (this is discussed and illustrated in Section~\ref{2D1}). The cell is seen as a 3-ball, that is, a continuum of concentric spheres together with the central point. Other examples comprise the formation of waterspouts where we see the formation of the tornado's cylindrical `cork' (as described and illustrated in Section~\ref{2D0}) and the creation of Falaco solitons where the creation of two discs joined with an 'invisible' thread is taking place in a water pool (as detailed and illustrated in Section~\ref{2D0}).

{\noindent} $\bullet$   {\bf Embedded surgery:} All phenomena exhibiting surgery take place in the ambient 3-space. For this reason we introduce the notion of embedded 1- or 2-dimensional surgery, which is taking place on an embedding of the initial manifold in $3$-space, \hl{instead of happening abstractly}. The ambient 3-space leaves room for the initial manifold to assume a more complicated configuration and allows the complementary space of the initial manifold to participate actively in the process. \hl{For example, in dimension 1 during DNA recombination, the initial DNA molecule which is recombined can also be knotted (see description and illustration in Section~\ref{1DDynamics}). In other words,  the initial $1$-manifold can be a knot (an embedding of the circle) instead of an abstract circle. Examples in dimension 2 comprise the processes of tornado and black hole formation (see Section~\ref{2D0} and illustration therein), which are not confined to the initial manifold, and topological surgery is causing (or is caused by) a change in the whole space.} 

{\noindent} $\bullet$   {\bf Connection between $1$- and $2$-dimensional surgeries:} As we explain then, the appearance of forces, enhanced with the notions of solid 1- and 2-dimensional surgery, can be all connected via appropriate (planar, spherical or toroidal) cross-sections. In fact all the above culminate to the notion of embedded solid 2-dimensional surgery and can be derived from there.

{\noindent} $\bullet$   {\bf Connection with a dynamical system:} Finally, we establish a connection between these new notions applied on 2-dimensional topological surgery and the dynamical system presented in \cite{SaGr1}. We analyze how, with a slight perturbation of parameters, trajectories pass from spherical to toroidal shape through a `hole drilling' process. We show that our new topological notions are verified by both the local behavior of the steady state points of the system and the numerical simulations of its trajectories. This result gives us on the one hand a mathematical model for 2-dimensional surgery and on the other hand a system that can model natural phenomena exhibiting these types of surgeries.

The paper is organized as follows: In Section~\ref{Useful} we recall the topological notions  that will be used and provide specific examples that will be of great help to readers that are not familiar with these mathematical notions. In Section~\ref{definitions}, we present and discuss the formal definition of topological surgery. In Section~\ref{1D}, we introduce dynamics to 1-dimensional surgery, we define solid 1-dimensional surgery and we discuss 1-dimensional natural processes exhibiting these types of surgeries. In Section~\ref{2D} we extend these definitions to 2-dimensional surgery and discuss related 2-dimensional natural processes. We then use these new theoretical concepts in Section~\ref{Connecting} to pin down the relations among topological surgeries of different dimensions. As all natural phenomena exhibiting surgery (1 or 2-dimensional, solid or usual) take place in the ambient 3-space, in Section~\ref{S3} we present the 3-sphere $S^3$ and the duality of its descriptions. This allows us to define in Section~\ref{es3sp} the notion of embedded surgery. Finally, our connection of solid 2-dimensional surgery with a dynamical system is established in Section~\ref{ds}. 

\section{Useful mathematical notions}\label{Useful}
\hl{In this section we introduce basic notions related to topological surgery. Reader that are familiar with the formalism of the topic can directly move to the formal definition in  Section~{\ref{definitions}}.}

\subsection{Manifolds}\label{Manifolds}
\begin{itemize} 
\item{\hl{An \textit{$n$-manifold without boundary} is a `nice' topological space with the property that each point in it has a neighborhood topologically equivalent to the usual $n$-dimensional Euclidean space ${\mathbb R}^n$. In other words an $n$-manifold resembles locally ${\mathbb R}^n$.}}
\item {\hl{Similarly, an \textit{$n$-manifold with boundary} is `nice' topological space with the property that each point in it has a neighborhood topologically equivalent either to ${\mathbb R}^n$ (if the point lies in the interior) or ${\mathbb R}_+^n$ (if the point lies on the boundary).}}
\end{itemize}

\subsection{Homeomorphisms} \label{Homeomorphisms}
\hl{In Section}~\ref{Manifolds} \hl{by} `\hl{topologically equivalent}' \hl{we mean the following: two $n$-manifolds $X$ and $Y$ are \textit{homeomorphic} or \textit{topologically equivalent} if there exists a homeomorphism between them, namely a function $f: X \rightarrow Y$ with the properties that}:
\begin{itemize} 
\item{\hl{$f$ is continuous}}
\item {\hl{There exists the inverse function $f^{-1}: Y \rightarrow X$ (equivalently $f$ is 1-1 and onto)}}
\item {\hl{$f^{-1}$ is also continuous}}
\end{itemize}
\hl{Intuitively the homeomorphism $f$ is an elastic deformation of the space $X$ to the space $Y$, not involving any self-intersections or any `cutting' and `regluing' (see also Appendix {\ref{Appendix}})}.

\subsection{Properties of manifolds}\label{ProperMan}
\hl{An $n$-manifold, $M$, is said to be:}
\begin{itemize} 
\item{\hl{\textit{connected} if it consists of only one piece}},
\item {\hl{\textit{compact} if it can be enclosed in some $k$-dimensional ball}},
\item {\hl{\textit{orientable} if any oriented frame that moves along any closed path in $M$ returns to a position that can be transformed to the initial one by a rotation}}.
\end{itemize}
\hl{The above notions are more rigorously defined in Appendix {\ref{Appendix}}}.

\subsection{$n$-spheres and $n$-balls}\label{NshperesAndNballs}
\begin{figure}[!h]
\begin{center}
\centering
\captionsetup{justification=centering}
\includegraphics[width=5cm]{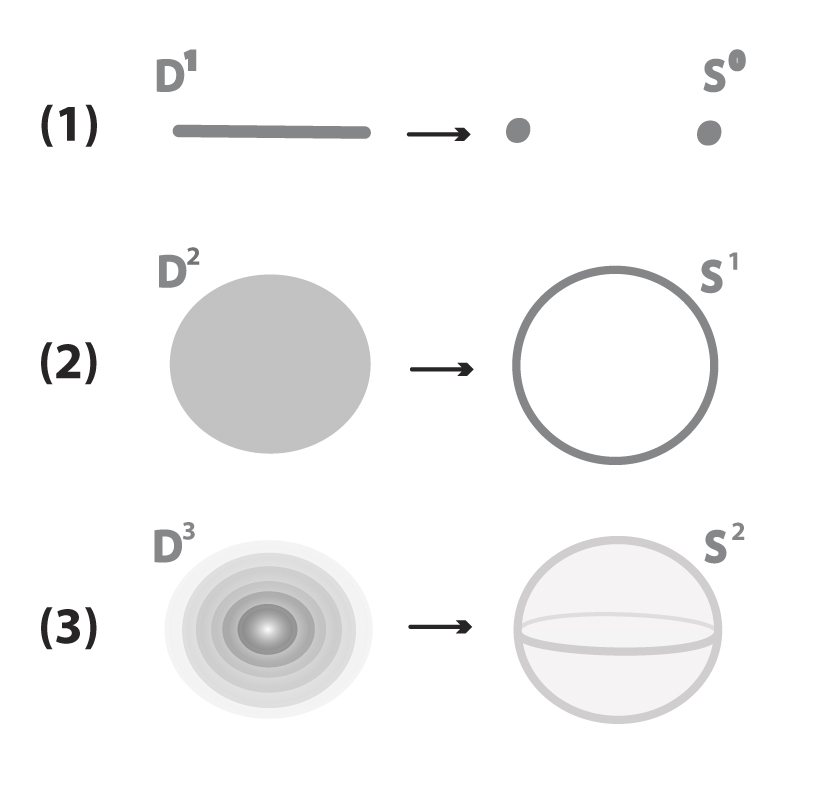}
\caption{\hl{\textbf{ (1) A segment $D^1$ is bounded by two points $S^0$ (2) A disc $D^2$ is bounded by a circle $S^1$ (3) A 3-ball $D^3$ is bounded by a sphere $S^2$}.}}
\label{Fig1}
\end{center}
\end{figure}

In each dimension the basic connected, oriented, compact $n$-manifold without boundary is the $n$-sphere, $S^n$. Also, the basic connected, oriented $n$-manifold with boundary is the $n$-ball, $D^n$. \hl{The boundary of a $n$-dimensional ball is a $n-1-$dimensional sphere, $\partial D^n=S^{n-1}, n>=1$. In Fig~{\ref{Fig1}}, this relation is shown for $n=1,2$ and $3$. As shown in Fig~{\ref{Fig1}} (1), the space ${S^0}$ is the disjoint union of two points. By convention, we consider these two one-point spaces to be $\{+1\}$ and  $\{-1\}$: ${S^0}=\{+1\} \amalg \{-1\}$.}  
  
\smallbreak

Besides the relation of $S^n$ with $D^{n+1}$ described above, the $n$-sphere $S^n$ is also intrinsically related to the Euclidian space ${\mathbb R}^n$ via the notion of compactification.

\subsection{The compactification of  ${\mathbb R}^n$}\label{CompactificationRn}

\begin{figure}[!h]
\begin{center}
\centering
\captionsetup{justification=centering}
\includegraphics[width=6cm]{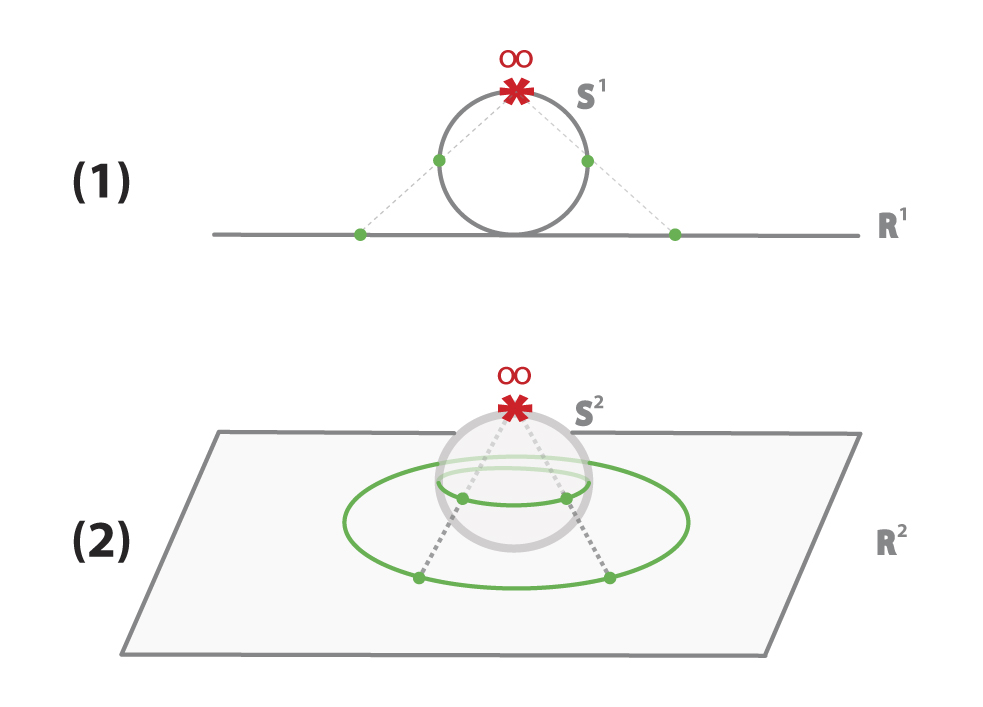}
\caption{\hl{\textbf{(1) $S^1$ onto ${\mathbb R}^1$ (2) $S^2$ onto ${\mathbb R}^2$.}}}
\label{Fig2}
\end{center}
\end{figure}

\hl{\textit{Compactification} is the process of making a topological space into a compact space.  
For each dimension $n$, the space ${\mathbb R}^n$ with all points at infinity compactified to one single point is homeomorphic to $S^n$. So, $S^n$ is also called the one-point compactification of ${\mathbb R}^n$. Conversely, a sphere $S^n$ can be decompactified to the space ${\mathbb R}^n$ by the so-called \textit{stereographic projection.} For example, for $n=1$ we have that the circle $S^1$ is the one-point compactification of the real line ${\mathbb R}^1$, see Fig~{\ref{Fig2}}~(1), while for $n=2$ the sphere $S^2$ is the one-point compactification of the plane ${\mathbb R}^2$, see  Fig~{\ref{Fig2}}~(2). The compactification of ${\mathbb R}^3$ is discussed and illustrated  in Section~{\ref{ViaR3}} (see Appendix {\ref{Appendix}} for details on the one-point compactification of ${\mathbb R}^n$).}

\subsection{Product spaces}\label{Productspaces}
\hl{The \textit{product space} of two manifolds $X$ and $Y$ is the manifold made from their Cartesian product $X \times Y$ (see also Appendix {\ref{Appendix}}). If $X, Y$ are manifolds with boundary, the boundary of product space $X \times Y$ is  $\partial (X \times Y)=(\partial X \times Y) \cup (X \times \partial Y)$.}

\hl{For example the next common connected, oriented, compact $2$-manifold without boundary after $S^2$ is the \textit{torus}, which can be perceived as the boundary of a doughnut, and it is the product space $S^1 \times S^1$. Analogously, a \textit{solid torus}, which can be perceived as a whole doughnut, is the product space $S^1 \times D^2$. A solid torus is a 3-manifold with boundary a torus}: 

\begin{samepage} 
 \begin{center}
\centering
\captionsetup{justification=centering}
\hl{$\partial (S^1 \times D^2)=S^1 \times \partial D^2=S^1 \times S^1$}
\end{center}  
\end{samepage}

\begin{figure}[!h]
\begin{center}
\centering
\captionsetup{justification=centering}
\includegraphics[width=6cm]{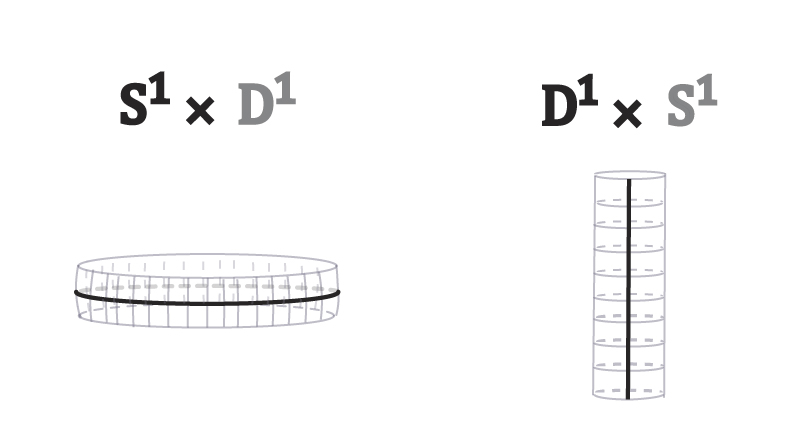}
\caption{\hl{\textbf{Two ways of viewing a cylinder}}}
\label{Fig3}
\end{center}
\end{figure}

\hl{Other product spaces that we will be using here are: the \textit{cylinder} $S^1 \times D^1$ or $D^1 \times S^1$ (see Fig{~\ref{Fig3}}), the \textit{solid cylinder} $D^2 \times D^1$ which is homeomorphic to the 3-ball and the spaces of the type ${S^0} \times {D^n}$, which are the disjoint unions of two $n$-balls ${D^n} \amalg {D^n}$.}

\smallbreak

\hl{All the above examples of product spaces that are of the form ${S^p} \times {D^q}$ can be viewed as $q$-thickenings of the $p$-sphere. For example the $2$-thickening of $S^0$ comprises two discs, while the $3$-thickening of $S^0$ comprises two 3-balls. It is also worth noting that the product spaces $S^p \times D^q$ and $D^{p+1} \times S^{q-1}$ have the same boundary: $\partial (S^p \times D^q)=\partial (D^{p+1} \times S^{q-1})=S^p \times S^{q-1} (\star)$.   }

\subsection{Embeddings}\label{EmbeddingsMan}
\begin{figure}[!h]
\begin{center}
\centering
\captionsetup{justification=centering}
\includegraphics[width=6cm]{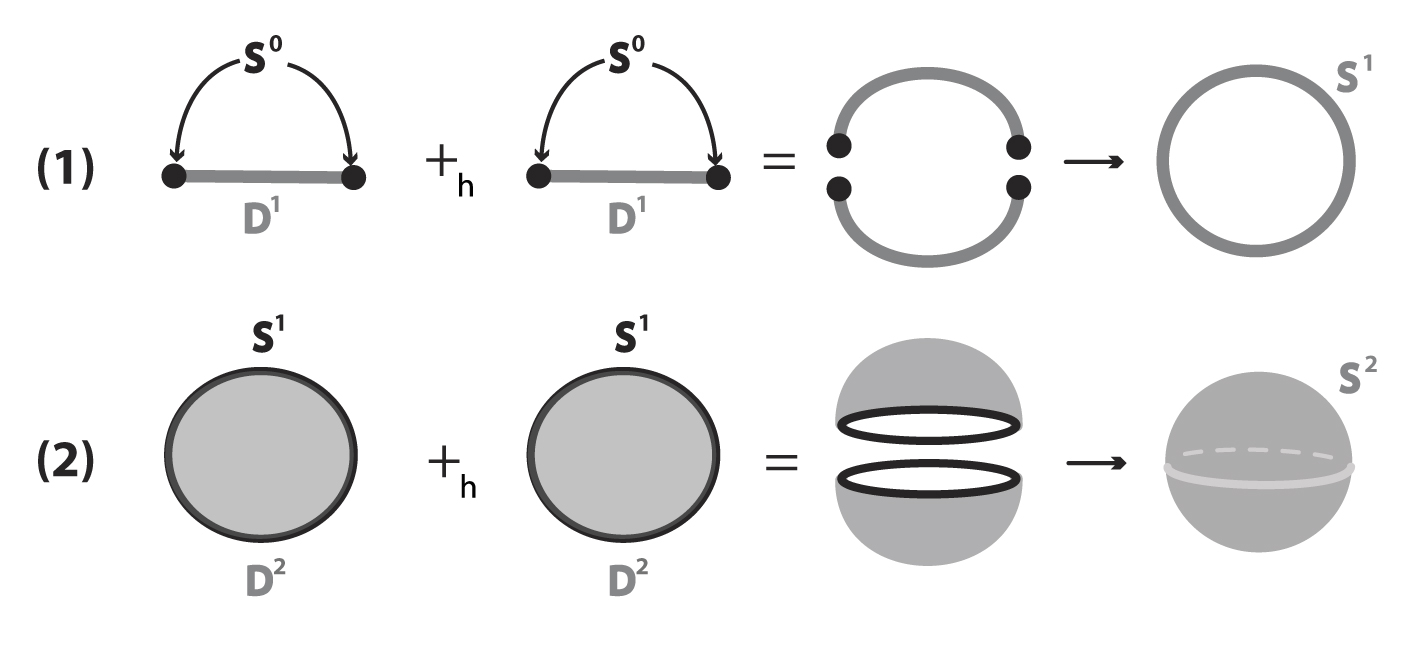}
\caption{\hl{\textbf{(1) $D^1 \cup_h D^1 = S^1$ (2) $D^2 \cup_h D^2 = S^2$ }}}
\label{Fig4}
\end{center}
\end{figure}

\begin{itemize} 
\item{\hl{An \textit{embedding} of an $n$-manifold $N^n$ in an $m$-manifold $M^m$ is a map $f: N \hookrightarrow M$ such that its restriction on the image $f(N)$ is a homeomorphism between $N$ and $f(N)$. The notion of embedding allows to view spaces inside specific manifolds instead of abstractly.
Embeddings even of simple manifolds can be very complex. For example, the embeddings of the circle $S^1$ in the $3$-space  ${\mathbb R}^3$ are the well-known knots whose topological classification is still an open problem of low-dimensional topology.}}
\item{\hl{An embedding of a submanifold $N^n \hookrightarrow  M^m$ is \textit{framed} if it extends to an embedding $N^{n} \times D^{m-n} \hookrightarrow  M$.}}
\item{\hl{A \textit{framed $n$-embedding in $M$} is an embedding of the $(m-n)$-thickening of the $n$-sphere, $h:S^n\times D^{m-n}\hookrightarrow  M$, with core $n$-embedding $e=h_{|}: {S^n}={S^n} \times \{0\} \hookrightarrow  M$. For example, the framed $1$-embeddings in ${\mathbb R}^3$ comprise embedded solid tori in the $3$-space with core $1$-embeddings being knots.}}
\item{\hl{Let $X$, $Y$ be two $n$-manifolds with homeomorphic boundaries $\partial X$ and $\partial Y$ (which are $(n-1)-$manifolds). Let also $h$ denote a homeomorphism $h:\partial X \rightarrow \partial Y$. Then, from $X \cup Y$ one can create a new $n$-manifold without boundary by `gluing' $X$ and $Y$ along their boundaries. The gluing is realized by identifying each point $x \in \partial X$ to the point $h(x) \in  \partial Y$. The map $h$ is called \textit{gluing homeomorphsim}, see Appendix {\ref{Appendix}}. One important example is the gluing of two $n$-discs along their common boundary which gives rise to the $n$-sphere, see Fig{~\ref{Fig4}} for n=1,2. For n=3, the gluing of two $3$-balls yielding the $3$-sphere $S^3$ is illustrated and discussed in Section~{\ref{ViaTwoBalls}}. Another interesting example is the gluing of solid tori which also yield the $3$-sphere. This is illustrated and discussed in Section~{\ref{ViaTwoTori}} }}
\end{itemize}

\smallbreak

As we will see in next section, the notions of embedding and gluing homeomorphism together with property $(\star)$ described in {\ref{Productspaces}} are the key ingredients needed to define topological surgery. It is roughly the procedure of removing an embedding of $S^p \times D^q$ and gluing back $D^{p+1} \times S^{q-1}$ along their common boundary.

\section{The formal definition of surgery} \label{definitions}

We recall the following well-known definition of surgery:

\begin{defn} \label{surgery} \rm An \textit{m-dimensional n-surgery} is the topological procedure of creating a new $m$-manifold $M'$ out of a given $m$-manifold $M$ by removing a framed $n$-embedding $h:S^n\times D^{m-n}\hookrightarrow  M$, and replacing it with $D^{n+1}\times S^{m-n-1}$, using the `gluing' homeomorphism $h$ along the common boundary $S^n\times S^{m-n-1}$. Namely, and denoting surgery by  $\chi$ :
\[M' = \chi(M) = \overline{M\setminus h(S^n\times D^{m-n})} \cup_{h|_{S^n\times S^{m-n-1}}} (D^{n+1}\times S^{m-n-1}). \]
The symbol `$\chi$' of surgery comes from the Greek word `$\chi\epsilon\iota\rho o\upsilon\rho\gamma\iota\kappa\acute{\eta}$' (cheirourgiki) whose term `cheir' means hand. Note that from the definition, we must have $n+1 \leq m$. \hl{Also, the horizontal bar in the above formula indicates the topological closure of the set underneath.} \end{defn}

{\noindent}Further, \rm the \textit{dual m-dimensional $(m-n-1)$-surgery} on $M'$ removes a dual framed $(m-n-1)$-embedding  $g:D^{n+1}\times S^{m-n-1}\hookrightarrow  M'$ such that $g|_{S^n\times S^{m-n-1}}=h^{-1}|_{S^n\times S^{m-n-1}}$, and replaces it with $S^n\times D^{m-n}$, using the `gluing' homeomorphism $g$ (or $h^{-1}$) along the common boundary $S^n\times S^{m-n-1}$. That is:
\[M = \chi^{-1}(M') = \overline{M'\setminus g(D^{n+1}\times S^{m-n-1})} \cup_{h^{-1}|_{S^n\times S^{m-n-1}}} (S^n\times D^{m-n}). \]

{\noindent}Note that resulting manifold $\chi(M)$ may or may not be homeomorphic to $M$. From the above definition, it follows that $M = \chi^{-1}(\chi(M))$. Preliminary definitions behind the definitions of surgery such as topological spaces, homeomorphisms, embeddings and other related notions are provided in Section~\ref{Useful} and Appendix \ref{Appendix}. For further reading, excellent references on the subject are \cite{Ra, PS, Ro}. \hl{We shall now apply the above definition to dimensions 1 and 2.}

\subsection{1-dimensional 0-surgery}\label{1D_FormalE}
\begin{figure}[!h]
\begin{center}
\centering
\captionsetup{justification=centering}
\includegraphics[width=11cm]{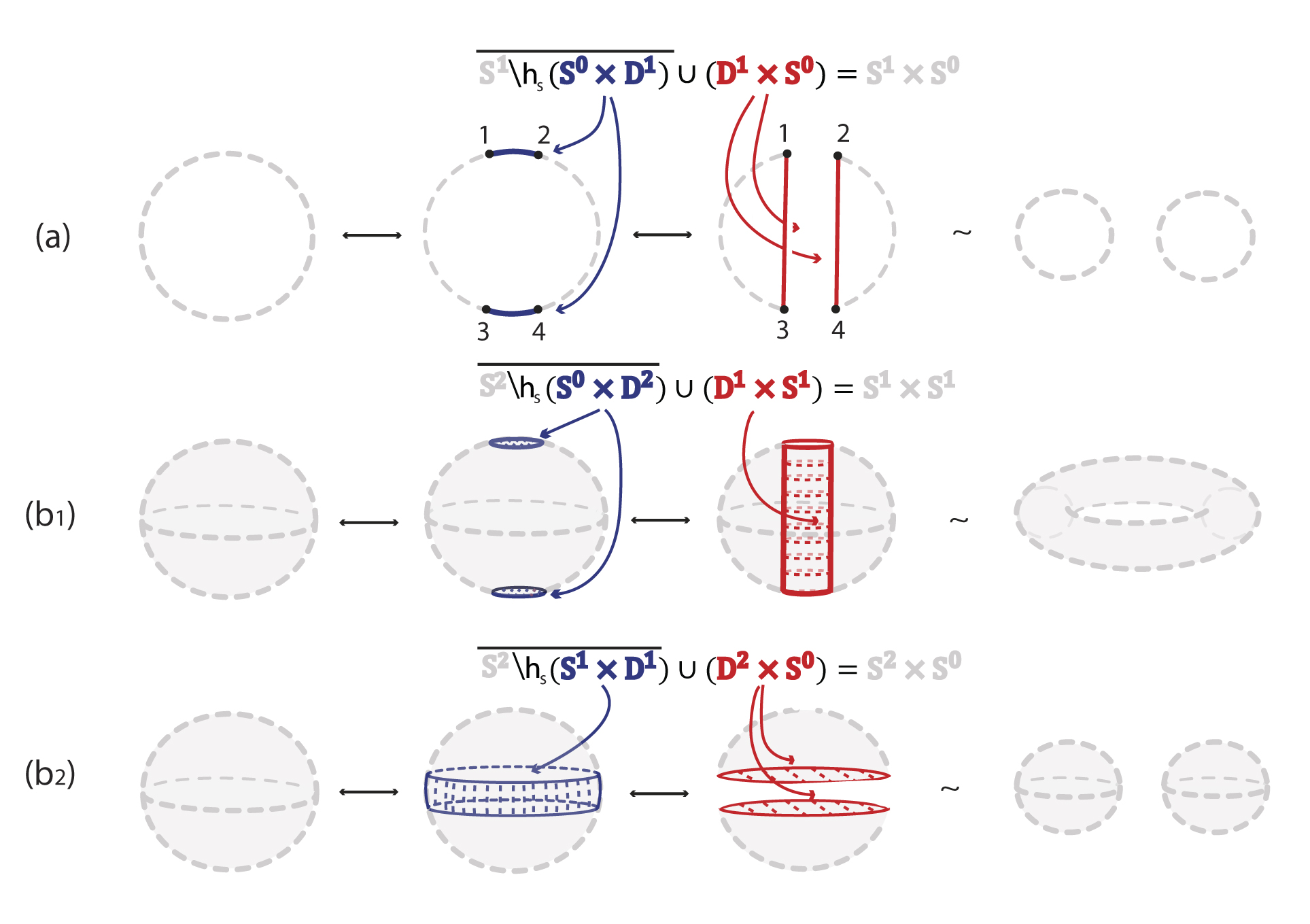}
\caption{\hl{\textbf{Formal (a) 1-dimensional 0-surgery (b{\textsubscript{1}}) 2-dimensional 0-surgery and (b{\textsubscript{2})} 2-dimensional 1-surgery using the standard embedding $h_s$}.}}
\label{Fig5}
\end{center}
\end{figure}

\begin{figure}[!h]
\begin{center}
\centering
\captionsetup{justification=centering}
\includegraphics[width=11cm]{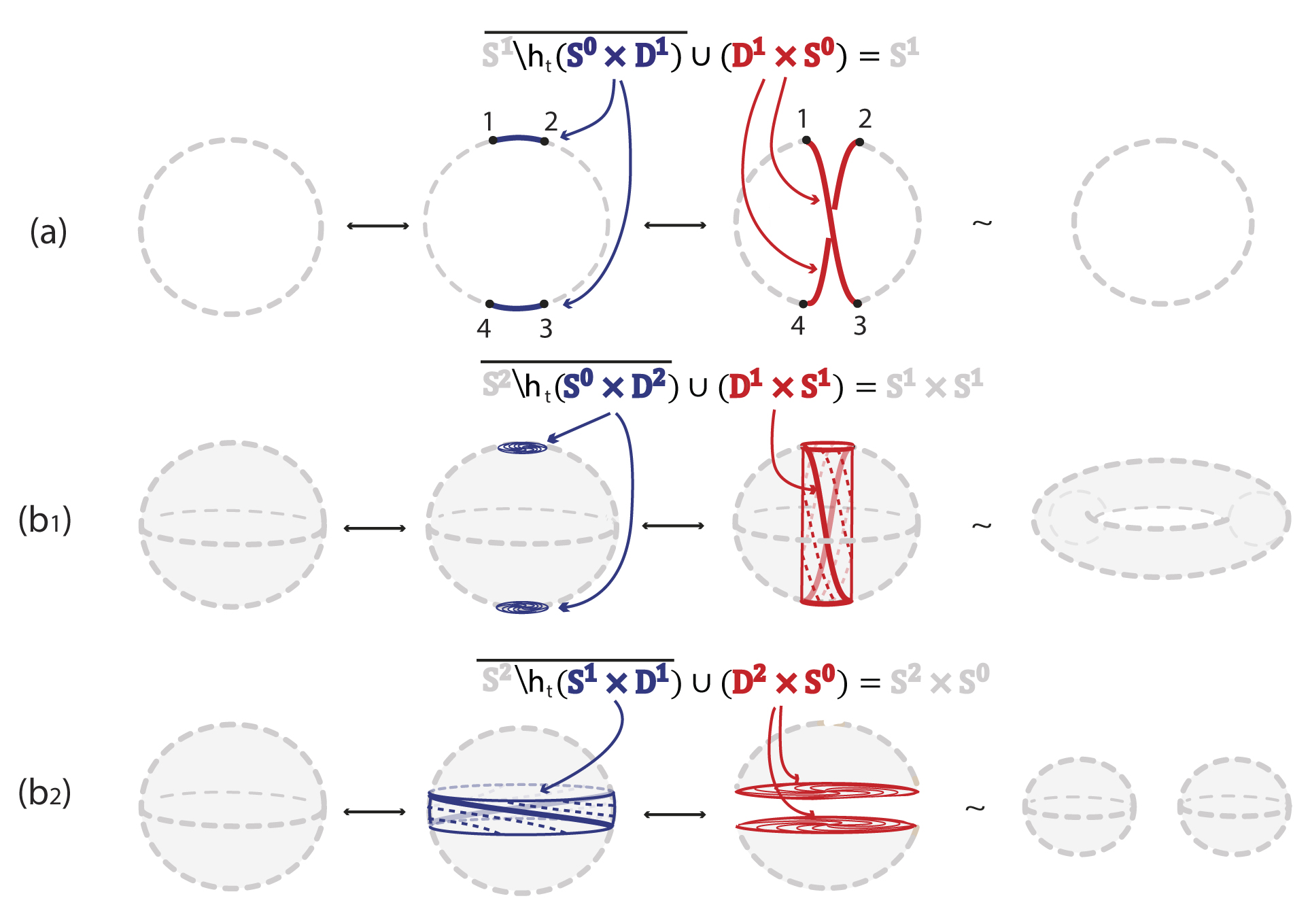}
\caption{\hl{\textbf{Formal (a) 1-dimensional 0-surgery (b{\textsubscript{1}}) 2-dimensional 0-surgery and (b{\textsubscript{2})} 2-dimensional 1-surgery using a twisting embedding $h_t$}.}}
\label{Fig6}
\end{center}
\end{figure}

We only have one kind of surgery on a 1-manifold $M$, the   \textit{1-dimensional 0-surgery} where $m=1$ and $n=0$: \[M' = \chi(M) = \overline{M\setminus h(S^0\times D^{1})} \cup_{h|_{S^0\times S^{0}}} D^{1}\times S^{0}. \] The above definition means that two segments $S^0\times D^1$ are removed from $M$ and they are replaced by two different segments $D^1 \times S^0$ by reconnecting the four boundary points $S^0\times S^0$ in a different way. \hl{In Fig~{\ref{Fig5}} (a) and ~{\ref{Fig6}} (a),  $S^0\times S^0 =\{1, 2, 3, 4\}$. As one possibility, if we start with $M=S^1$ and use as $h$ the standard (identity) embedding denoted with $h_s$,  we obtain two circles $S^1 \times S^0$. Namely,  denoting by $1$ the identity homeomorphism,  we have $h_s:S^0\times D^{1}=D^{1} \amalg D^{1} \xrightarrow{1 \amalg  1}  S^0\times D^{1} \hookrightarrow M$,  see  Fig~{\ref{Fig5}} (a). However,  we can also obtain one circle $S^1$ if $h$ is an embedding $h_t$ that reverses the orientation of one of the two arcs of $S^0\times D^1$. Then in the substitution,  joining endpoints 1 to 3 and 2 to 4,  the two new arcs undergo a half-twist, see Fig~{\ref{Fig6}} (a). More specifically,  if we take ${D^1}=[-1, +1]$ and define the homeomorphism  $\omega:D^{1}\to D^{1} ; t \to -t$,  the embedding used in Fig~{\ref{Fig6}} (a) is $h_t:S^0\times D^{1}=D^{1} \amalg D^{1} \xrightarrow{1 \amalg  \omega}  S^0\times D^{1} \hookrightarrow M $ which rotates one $D^{1}$ by 180{\si{\degree}}. The difference between the embeddings $h_s$ and $h_t$ of $S^0\times D^1$ can be clearly seen by comparing the four boundary points $1, 2, 3$ and $4$ in Fig~{\ref{Fig5}} (a) and  Fig~{\ref{Fig6}} (a).}

\smallbreak

\smallbreak

Note that in dimension one,  the dual case is also an 1-dimensional 0-surgery. For example,  looking at the reverse process of Fig~\ref{Fig5} (a),  we start with two circles $M'=S^1\amalg S^1$ and,  if each  segment of $D^1 \times S^0$ is embedded in a different circle,  the result of the (dual) 1-dimensional 0-surgery is one circle: $ \chi^{-1}(M')=M=S^1$.

\subsection{2-dimensional 0-surgery}\label{2D_FormalE0}
Starting with a 2-manifold $M$,  there are two types of surgery. One type is the \textit{ 2-dimensional 0-surgery},  whereby two discs $S^0\times D^2$ are removed from  $M$ and are replaced in the closure of the remaining manifold by a cylinder $D^1\times S^1$,  which gets attached via a homeomorphism along the common boundary $S^0\times S^1$ comprising two copies of $S^1$. The gluing homeomorphism of the common boundary may twist one or both copies of $S^1$. For $M=S^2$ the above operation changes its homeomorphism type from the 2-sphere to that of the torus. \hl{View Fig~{\ref{Fig5}~({b\textsubscript{1}}}) for the standard embedding  $h_s$ and Fig~{\ref{Fig6}~({b\textsubscript{1}}}) for a twisting embedding $h_t$. For example,  the  homeomorphism  $\mu:D^{2}\to D^{2} ; (t_1, t_2) \to (-t_1, -t_2)$ induces the 2-dimensional analogue $h_t$ of the embedding defined in the previous example,  namely: $h_t:S^0\times D^{2}=D^{2} \amalg D^{2} \xrightarrow{1 \amalg  \mu}  S^0\times D^{2} \hookrightarrow M $ which rotates one $D^{2}$ by 180{\si{\degree}}. When,  now,  the cylinder $D^1\times S^1$ is glued along the common boundary  $S^0\times S^1$,  the twisting of this boundary induces the twisting of the cylinder,  see Fig~{\ref{Fig6}~({b\textsubscript{1}}}).}

\subsection{2-dimensional 1-surgery}\label{2D_FormalE1}
The other possibility of 2-dimensional surgery on $M$ is the \textit{2-dimensional 1-surgery}: here a cylinder (or annulus) $S^1 \times D^1$ is removed from $M$ and is replaced in the closure of the remaining manifold by two discs  $D^2 \times S^0$ attached along the common boundary $S^1 \times S^0$. For $M=S^2$ the result is two copies of  $S^2$,  \hl{see Fig~{\ref{Fig5}~({b\textsubscript{2}}}) for the standard embedding  $h_s$. Fig~{\ref{Fig6}~({b\textsubscript{2}}}) illustrates a twisting embedding $h_t$,  where a twisted cylinder is being removed. In that case,  taking $D^{1}=\{h:h \in [-1, 1]\}$ and homeomorphism $\zeta$:} 

\begin{samepage} 
 \begin{center}
\hl{$\zeta:S^1\times D^{1}\to S^1\times D^{1};$ }
 \nopagebreak 
 \\[5pt] 
\hl{$\zeta: (t_1, t_2, h) \to (t_1\cos{\frac{(1-h)\pi}{2}}-t_2\sin{\frac{(1-h)\pi}{2}}, t_1\sin{\frac{(1-h)\pi}{2}}+t_2\cos{\frac{(1-h)\pi}{2}}, h)$}
\end{center}  
 \end{samepage} 

{\noindent}\hl{the embedding $h_t$ is defined as: $h_t:S^1\times D^{1} \xrightarrow{\zeta} S^1\times D^{1} \hookrightarrow M $. This operation corresponds to fixing the circle $S^1$ bounding the top of the cylinder $S^1 \times D^1$,  rotating the circle $S^1$ bounding the bottom of the cylinder by 180{\si{\degree}} and letting the rotation propagate from bottom to top. This twisting of the cylinder can be seen by comparing the second instance of Fig~{\ref{Fig5}~({b\textsubscript{2}}}) with  the second instance of Fig~{\ref{Fig6}~({b\textsubscript{2}}}),  but also by comparing the third instance of Fig~{\ref{Fig5}~({b\textsubscript{1}}}) with  the third instance of Fig~{\ref{Fig6}~({b\textsubscript{1}}}).}

It follows from Definition~\ref{surgery} that a dual 2-dimensional 0-surgery is a 2-dimensional 1-surgery and vice versa. Hence,  Fig~\ref{Fig5}~({b\textsubscript{1}}) shows that a 2-dimensional 0-surgery on a sphere is the reverse process of a 2-dimensional 1-surgery on a torus. Similarly,  as illustrated in Fig~\ref{Fig5}~({b\textsubscript{2}}),  a 2-dimensional 1-surgery on a sphere is the reverse process of a 2-dimensional 0-surgery on two spheres. In the figure the symbol $\longleftrightarrow $ depicts surgeries from left to right and their corresponding dual surgeries from right to left.

\section{1-dimensional topological surgery}\label{1D}

\begin{figure}[!h]
\begin{center}
\centering
\captionsetup{justification=centering}
\includegraphics[width=13cm]{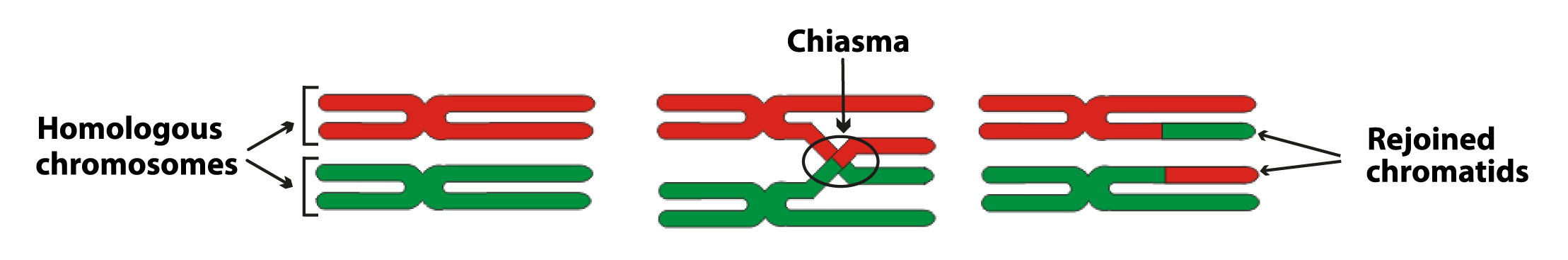}
\caption{{\bf Crossing over of chromosomes during meiosis.}}
\label{Fig7}
\end{center}
\end{figure}

\begin{figure}[!h]
\begin{center}
\centering
\captionsetup{justification=centering}
\includegraphics[width=11cm]{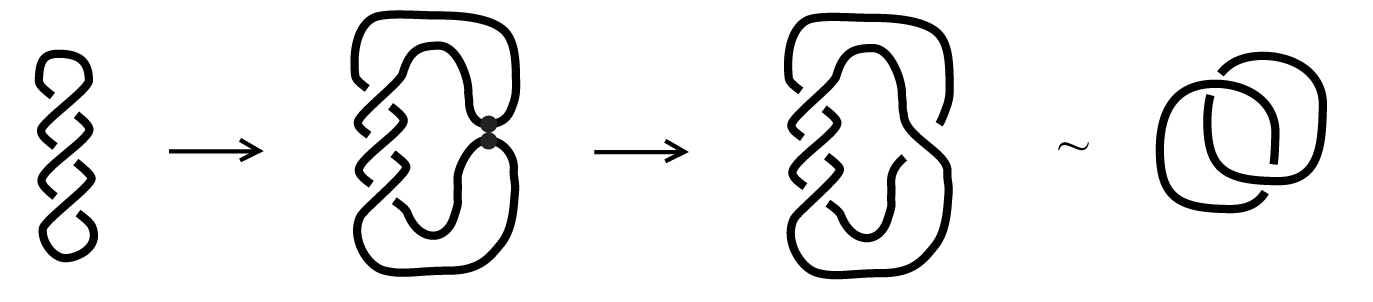}
\caption{{\bf DNA Recombination.}}
\label{Fig8}
\end{center}
\end{figure}

\begin{figure}[!h]
\begin{center}
\centering
\captionsetup{justification=centering}
\includegraphics[width=10cm]{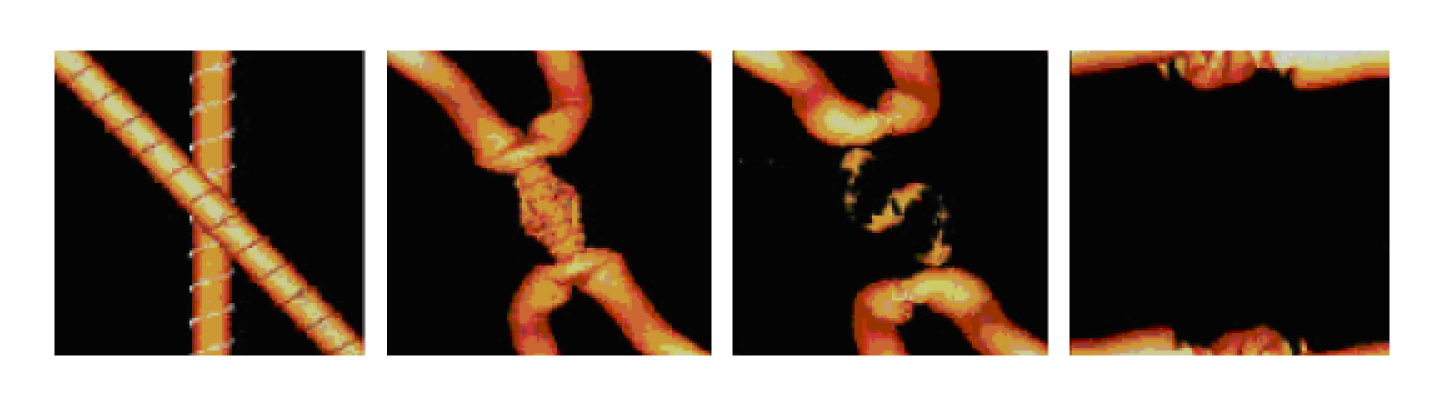}
\caption{{\bf The reconnection of cosmic magnetic lines.}}
\label{Fig9}
\end{center}
\end{figure}

1-dimensional 0-surgery happens in nature, in various scales, in phenomena where 1-dimensional splicing and reconnection occurs. For example, it happens on chromosomes during meiosis and produces new combinations of genes (see Fig~\ref{Fig7}), in site-specific DNA recombination (see Fig~\ref{Fig8}) whereby nature alters the genetic code of an organism, either by moving a block of DNA to another position on the molecule or by integrating a block of alien DNA into a host genome (see \cite{Su}), in magnetic reconnection, the phenomenon whereby cosmic magnetic field lines from different magnetic domains are spliced to one another, changing their patterns of connectivity with respect to the sources (see Fig~\ref{Fig9} from \cite{DaAn}) and in the reconnection of vortices in classical and quantum fluids (see \cite{LaRiSu}).

In this section we introduce dynamics which explains the process of 1-dimensional surgery, define the notion of solid 1-dimensional surgery and examine in more details the aforementioned natural phenomena.

\subsection{Introducing dynamics}\label{1DDynamics}

The formal definition of 1-dimensional 0-surgery gives a static description of the initial and the final stage whereas natural phenomena exhibiting 1-dimensional 0-surgery follow a continuous process. In order to address such phenomena or to understand how 1-dimensional 0-surgery happens, we need a non-static description.

Furthermore, in nature, 1-dimensional 0-surgery often happens locally, on arcs or segments. That is, the initial manifold is often bigger and we remove from its interior two segments $S^0 \times D^1$. Therefore, we also need dynamics that act locally.

\begin{figure}[!h]
\begin{center}
\centering
\captionsetup{justification=centering}
\includegraphics[width=13.4cm]{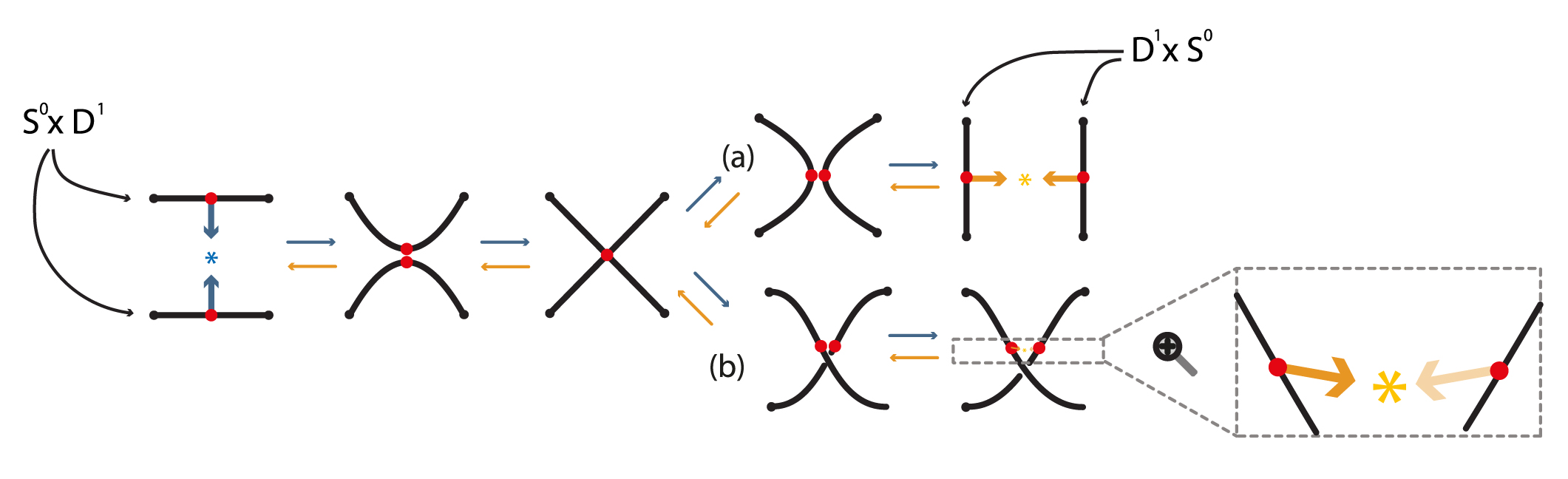}
\caption{{\bf Introducing dynamics to 1-dimensional surgery.}}
\label{Fig10}
\end{center}
\end{figure}

In Fig~\ref{Fig10}, we introduce dynamics which explain the intermediate steps of the formal definition and extend surgery to a continuous process caused by local forces. The process starts with the two points specified on the manifold (in red), on which attracting forces are applied (in blue). We assume that these forces are created by an attracting center (also in blue). Then, the two segments $S^0\times D^1$, which are neighborhoods of the two points, get close to one another. When the specified points (or centers) of two segments reach the attracting center, they touch and recoupling takes place giving rise to the two final segments $D^1 \times S^0$, which split apart. As mentioned in previous section, we have two cases (a) and (b), depending on the homemorphism $h$.

As mentioned in Section~\ref{1D_FormalE}, the dual case is also a 1-dimensional 0-surgery as it removes segments $D^1 \times S^0$ and replace them by segments $S^0\times D^1$. This is the reverse process which starts from the end and is illustrated in Fig~\ref{Fig10} as a result of the orange forces and attracting center which are applied on the `complementary' points.

\begin{rem}\label{Morse} \rm
It is worth mentioning that the intermediate steps of surgery presented in Fig~\ref{Fig10} can also be viewed in the context of Morse theory \cite{Mil}. By using the local form of a Morse function, we can visualize the process of surgery by varying parameter $t$ of equation $x^2-y^2=t$. For $t=-1$ it is the hyperbola shown in the second instance of Fig~\ref{Fig10} where the two segments get close to one another. For $t=0$ it is the two straight lines where the reconnection takes place as shown in the third instance of Fig~\ref{Fig10} while for $t=1$ it represents the hyperbola of the two final segments shown in case (a) of the fourth instance of Fig~\ref{Fig10}. This sequence can be generalized for higher dimensional surgeries as well, however, in this paper we will not use this approach as we are focusing on the introduction of forces and of the attracting center. 
\end{rem}

\begin{figure}[!h]
\begin{center}
\centering
\captionsetup{justification=centering}
\includegraphics[width=13.5cm]{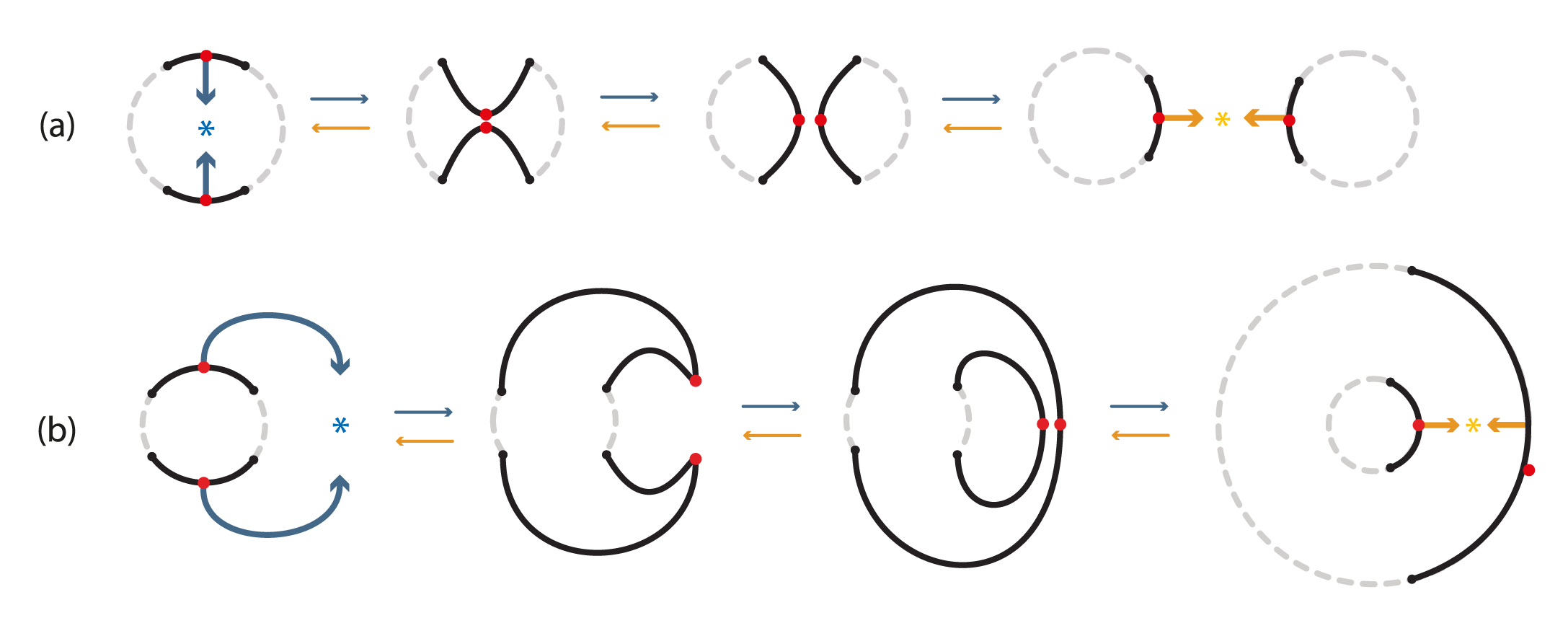}
\caption{{\bf 1-dimensional surgery on one and two circles.}}
\label{Fig11}
\end{center}
\end{figure}

These local dynamics produce different manifolds depending on where the initial neighborhoods are embedded. Taking the known case of the standard embedding $h_s$ and $M=S^{1}$, we obtain $S^1 \times S^0$ (for both regular and dual surgery), see  Fig~\ref{Fig11} (a). Furthermore, as shown in Fig~\ref{Fig11} (b), we also obtain $S^1 \times S^0$ even if the attracting center is outside $S^1$. Note that these outcomes are not different than the ones shown in formal surgery (recall Fig~\ref{Fig5} (a)) but we can now see the intermediate instances.

\subsection{Explaining 1-dimensional phenomena via dynamics}\label{1DPhenomena}

Looking closer at the aforementioned phenomena, the described dynamics and attracting forces are present in all cases. Namely, {\bf magnetic reconnection} (Fig~\ref{Fig9}) corresponds to a dual 1-dimensional 0-surgery (see Fig~\ref{Fig10} (b)) where $g:D^{1}\times S^{0}\hookrightarrow  M'$ is a dual embedding of the twisting homeomorphism $h_t$ defined in Section~\ref{1D_FormalE} of Section~\ref{definitions}. The tubes are viewed as segments and correspond to an initial manifold $M=S^0 \times D^1$ (or $M=S^1$ if they are connected) on which the local dynamics act on two smaller segments $S^0 \times D^1$. Namely, the two magnetic flux tubes have a nonzero parallel net current through them, which leads to attraction of the tubes (cf. \cite{KoFe}). Between them, a localized diffusion region develops where magnetic field lines may decouple. Reconnection is accompanied with a sudden release of energy and the magnetic field lines break and rejoin in a lower energy state.

In the case of {\bf chromosomal crossover} (Fig~\ref{Fig7}), we have the same dual 1-dimensional 0-surgery as magnetic reconnection (see Fig~\ref{Fig10} (b)). During this process, the homologous (maternal and paternal) chromosomes come together and pair, or synapse, during prophase. The pairing is remarkably precise and is caused by mutual attraction of the parts of the chromosomes that are similar or homologous. Further, each paired chromosomes divide into two chromatids. The point where two homologous non-sister chromatids touch and exchange genetic material is called chiasma. At each chiasma, two of the chromatids have become broken and then rejoined (cf. \cite{Pu}). In this process, we consider the initial manifold to be one chromatid from each chromosome, hence the initial manifold is $M=S^0 \times D^1$ on which the local dynamics act on two smaller segments $S^0 \times D^1$.

For {\bf site-specific DNA recombination} (see Fig~\ref{Fig8}), we have a 1-dimensional 0-surgery (see Fig~\ref{Fig10} (b)) with a twisted homeomorphism $h_t$ as defined in Section~\ref{1D_FormalE} of Section~\ref{definitions}. Here the initial manifold is a knot which is an embedding of $M=S^1$ in 3-space but this will be detailed in Section~\ref{es3sp}. As mentioned in \cite{JoLe}, enzymes break and rejoin the DNA strands, hence in this case the seeming attraction of the two specified points is realized by the enzyme. Note that, while both are genetic recombinations, there is a difference between chromosomal crossover and site-specific DNA recombination. Namely, chromosomal crossover involves the homologous recombination between two similar or identical molecules of DNA and we view the process at the chromosome level regardless of the knotting of DNA molecules. 

Finally, {\bf vortices reconnect} following the steps of 1-dimensional 0-surgery with a standard embedding shown in Fig~\ref{Fig10} (a). The initial manifold is again $M=S^0 \times D^1$. As mentioned in \cite{Ke}, the interaction of the anti-parallel vortices goes from attraction before reconnection, to repulsion after reconnection.

\subsection{Defining solid 1-dimensional surgery}

There are phenomena which undergo the process of 1-dimensional 0-surgery but happen on surfaces, such as {\bf tension on membranes or soap films} and the {\bf merging of oil slicks}. In order to model topologically such phenomena we introduce the notion of solid 1-dimensional 0-surgery. \textit{ Solid 1-dimensional 0-surgery on the $2$-disc $D^2$} is the topological procedure whereby a ribbon $D^1 \times D^1$ is being removed, such that the closure of the remaining manifold comprises two discs  $D^2 \times S^0$. The reader is referred to Fig~\ref{Fig5} (a) where the interior is now supposed to be filled in. This process is equivalent to performing 1-dimensional 0-surgeries on the whole continuum of concentric circles included in $D^2$. More precisely, and introducing at the same time dynamics, we define:

\begin{defn} \rm  We start with the $2$-disc of radius 1 with polar layering: 
$$
D^2 = \cup_{0<r\leq 1} S^1_r \cup \{P\},
$$ 
where $r$ the radius of a circle and $P$ the limit point of the circles, that is, the center of the disc. We specify colinear pairs of antipodal points \hl{, all on the same diameter,} with neighborhoods of analogous lengths, on which the same colinear attracting forces act, see Fig~\ref{Fig12} (1) where these forces and the corresponding attracting center are shown in blue. Then, in (2), antipodal segments get closer to one another or, equivalently, closer to the attracting center. Note that here, the attracting center coincides with the limit point of all concentric circles, which is shown in green from instance (2) and on. Then, as shown from (3) to (9), we perform 1-dimensional 0-surgery on the whole continuum of concentric circles. The natural order of surgeries is as follows: first, the center of the segments that are closer to the center of attraction touch, see (4). After all other points have also reached the center, see (5), decoupling starts from the central or limit point. We define 1-dimensional 0-surgery on the limit point $P$ to be the two limit points of the resulting surgeries. That is, the effect of \textit{solid 1-dimensional 0-surgery on a point is the creation of two new points} , see (6). Next, the other segments reconnect, from the inner, see (7), to the outer ones, see (8), until we have two copies of $D^2$, see (9) and (10). Note that the proposed order of reconnection, from inner to outer, is the same as the one followed by skin healing, namely, the regeneration of the epidermis starts with the deepest part and then migrates upwards.

\begin{figure}[!h]
\begin{center}
\centering
\captionsetup{justification=centering}
\includegraphics[width=13.5cm]{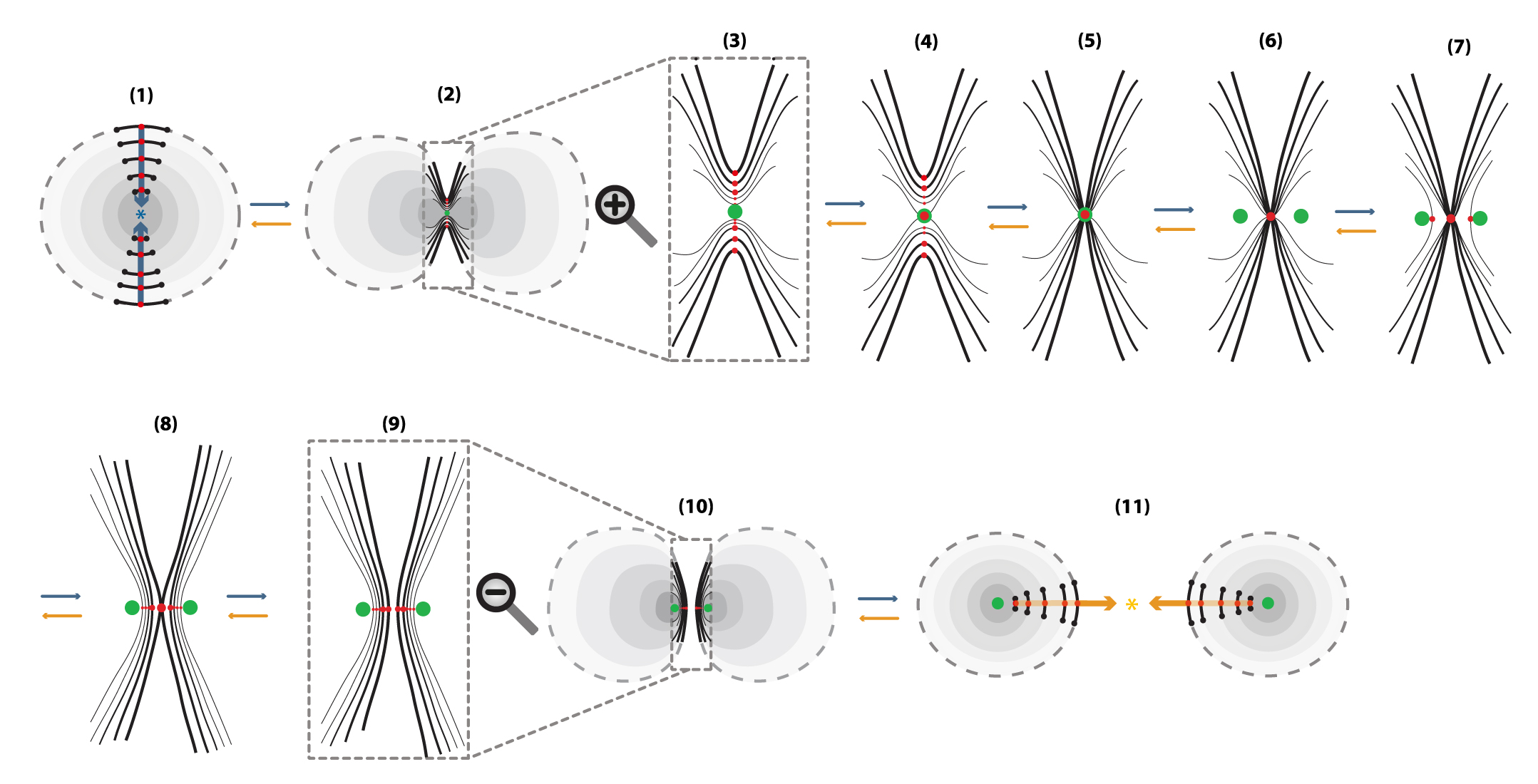}
\caption{{\bf Solid 1-dimensional surgery.}}
\label{Fig12}
\end{center}
\end{figure}

The above process is the same as first removing the center $P$ from $D^2$, doing the  1-dimensional 0-surgeries and then taking the closure of the resulting space. The resulting manifold is 
$$
\chi(D^2) := \cup_{0<r\leq 1}\chi(S^1_r) \cup \chi(P),
$$
which comprises two copies of $D^2$.

\smallbreak

We also have the reverse process of the above, namely,  \textit{ Solid 1-dimensional 0-surgery on two discs $D^2 \times S^0$}  is the topological procedure whereby a ribbon  $D^1 \times D^1$  joining the discs is added, such that the closure of the remaining manifold comprise one disc $D^2$, as illustrated in Fig~\ref{Fig12}. This process is the result of the orange forces and attracting center which are applied on the `complementary' points. This operation is equivalent to performing 1-dimensional 0-surgery on the  whole continuum of concentric circles in $D^2\amalg D^2$. We only need to define solid 1-dimensional 0-surgery on two limit points to be the limit point $P$ of the resulting surgeries. That is, the effect of \textit{ solid 1-dimensional 0-surgery on two points is their merging into one point}. The above process is the same as first removing the centers from the $D^2 \times S^0$, doing the 1-dimensional 0-surgeries and then taking the closure of the resulting space. The resulting manifold is 
$$
\chi^{-1}(D^2 \times S^0) := \cup_{0<r\leq 1}\chi^{-1}(S^1_r \times S^0) \cup \chi^{-1}(P \times S^0),
$$
which comprises one copy of $D^2$.

\end{defn}

\section{2-dimensional topological surgery}\label{2D}

\begin{figure}[!h]
\begin{center}
\centering
\captionsetup{justification=centering}
\includegraphics[width=8cm]{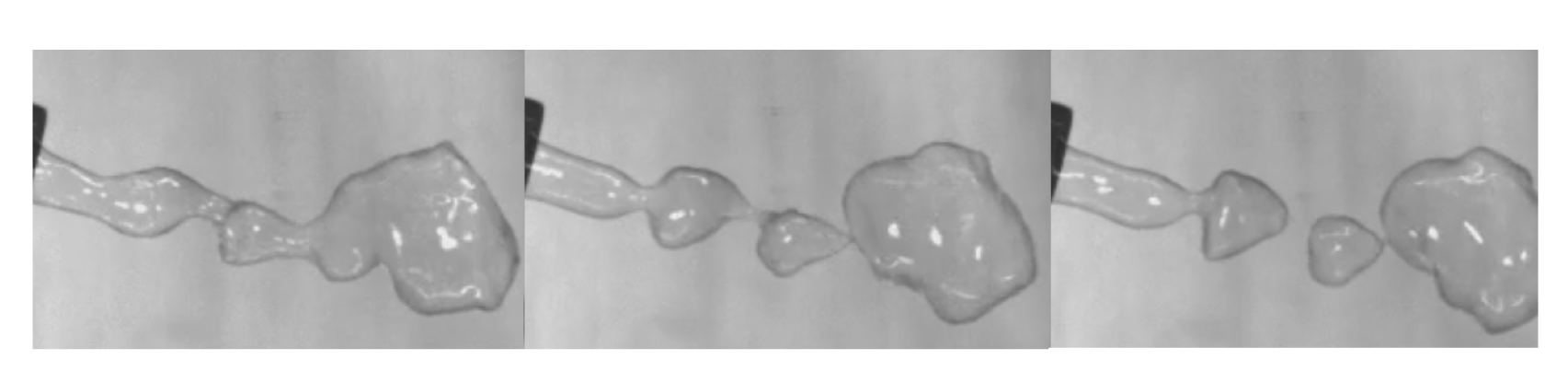}
\caption{{\bf Soap bubble splitting.} An example of 2-dimensional 1-surgery.}
{\footnotesize \textit{Source}: www.soapbubble.dk}
\label{Fig13}
\end{center}
\end{figure}

Both types of 2-dimensional surgeries are present in nature, in various scales, in phenomena where 2-dimensional merging and recoupling occurs. Natural processes undergoing \textit{ 2-dimensional 0-surgery} comprise, for example, drop coalescence, the formation of tornadoes and Falaco solitons, gene transfer in bacteria and the formation of black holes (for illustrations see Section~\ref{2D0}) . On the other hand, phenomena undergoing \textit{ 2-dimensional 1-surgery} comprise soap bubble splitting (see Fig~\ref{Fig13}), the biological process of mitosis and fracture as a result of tension on metal specimen (for illustrations see Section~\ref{2D1}). In this section we introduce dynamics which explains the process of 2-dimensional surgery, define the notions of solid 2-dimensional surgery and examine in more details the aforementioned natural phenomena. 

Note that except for soap bubble splitting which is a phenomena happening on surfaces, the other mentioned phenomena involve all three dimensions and are, therefore, analyzed after the introduction of solid 2-dimensional surgery, in Sections \ref{2D0} and \ref{2D1}.

\subsection{Introducing dynamics}\label{2DDynamics}
\begin{figure}[!h]
\begin{center}
\centering
\captionsetup{justification=centering}
\includegraphics[width=13.5cm]{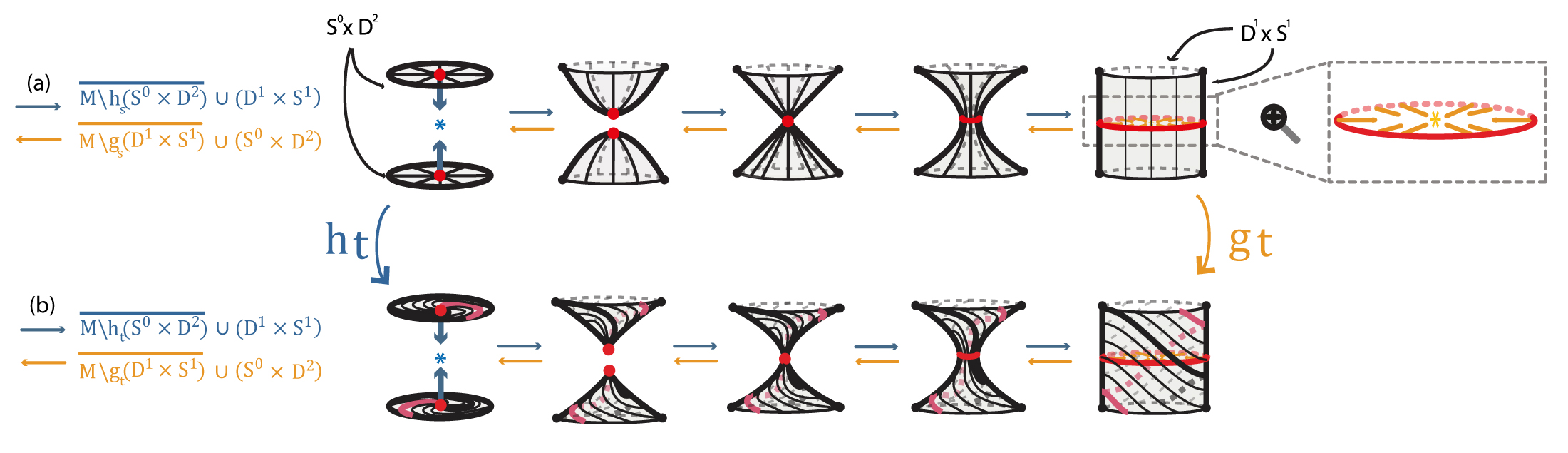}
\caption{{\bf Introducing dynamics to 2-dimensional surgery.} (a) 2-dimensional surgeries with standard embeddings (b) 2-dimensional surgeries with twisted embeddings.}
\label{Fig14}
\end{center}
\end{figure}

In order to model topologically phenomena exhibiting 2-dimensional surgery or to understand 2-dimensional surgery through continuity we need, also here, to introduce dynamics. In Fig~\ref{Fig14} (a), the 2-dimensional 0-surgery starts with two points, or poles, specified on the manifold (in red) on which attracting forces created by an attracting center are applied (in blue).  Then, the two discs $S^0\times D^2$,  neighborhoods of the two poles, approach each other. When the centers of the two discs touch, recoupling takes place and the discs get transformed into the final cylinder $D^1\times S^1$. 

As mentioned in Example~\ref{2D_FormalE1}, the dual case of 2-dimensional 0-surgery is the 2-dimensional 1-surgery and vice versa. This is also shown in Fig~\ref{Fig14} (a) where the reverse process is the \textit{2-dimensional 1-surgery} which starts with the cylinder and a specified cyclical region (in red) on which attracting forces created by an attracting center are applied (in orange). A `necking' occurs in the middle which degenerates into a point and finally tears apart creating two discs $S^0\times D^2$. As also seen in  Fig~\ref{Fig14} (a), in the case of 2-dimensional 0-surgery, forces (in blue) are applied on two points, or $S^0$, while in the case of the 2-dimensional 1-surgery, forces (in orange) are applied on a circle $S^1$.

In Fig~\ref{Fig14} (b), we have an example of \textit{twisted 2-dimensional 0-surgery} where the two discs $S^0\times D^2$  are embedded via a twisted homemorphism $h_t$ while, in the dual case, the cylinder $D^1\times S^1$ is embedded via a twisted homemorphism $g_t$. Here $h_t$ rotates the two discs while $g_t$ rotates the top and bottom of the cylinder by $3\pi/4$ and $-3\pi/4$ respectively. More specifically, if we define the homeomorphism  $\omega_1,\omega_2:D^{2}\to D^{2}$ to be rotations by $3\pi/4$ and $-3\pi/4$ respectively, then $h_t$ is defined as the composition $h_t:S^{0}\times D^2 \xrightarrow{\omega_1 \amalg  \omega_2}  S^0\times D^{2}  \xrightarrow{h} M $. The homeomorphism $g_t:D^1\times S^1  \rightarrow M $ is defined analogously.  

\begin{figure}[!h]
\begin{center}
\centering
\captionsetup{justification=centering}
\includegraphics[width=13cm]{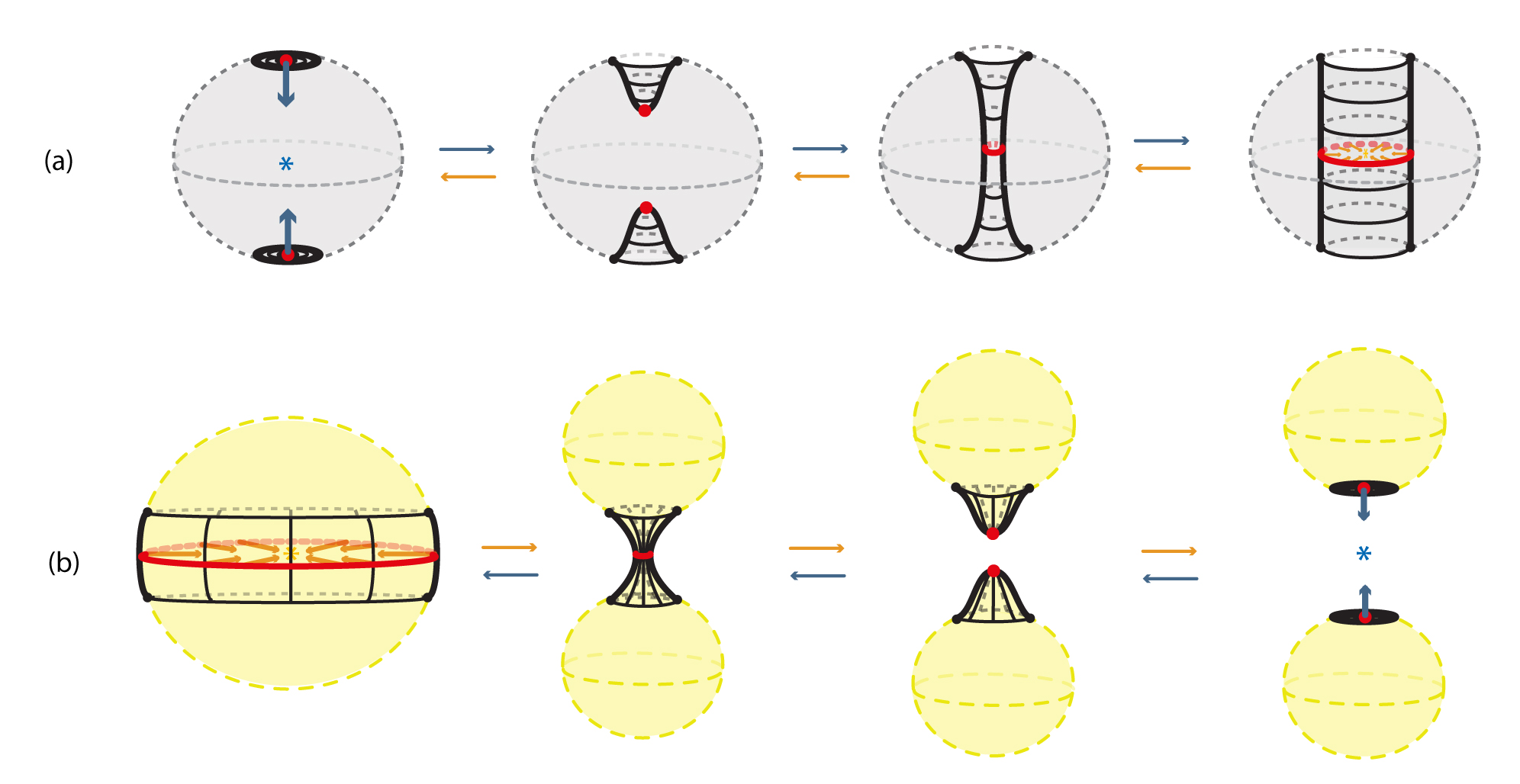}
\caption{{\bf } (a) 2-dimensional 0-surgery on $M=S^{2}$ and 2-dimensional 1-surgery on $M'=S^{0}\times S^{2}$ (b) 2-dimensional 1-surgery on $M=S^{2}$ and 2-dimensional 0-surgery on $M'=S^{0}\times S^{2}$.}
\label{Fig15}
\end{center}
\end{figure}

These local dynamics produce different manifolds depending on the initial manifold where the neighborhoods are embedded. Taking $M=S^{2}$, the local dynamics of Fig~\ref{Fig14} (a) are shown in Fig~\ref{Fig15} (a) and (b) producing the same manifolds seen in formal 2-dimensional surgery (recall Fig~\ref{Fig5}~(b\textsubscript{1})(b\textsubscript{2})). Note that, as also seen in 1-dimensional surgery (Fig~\ref{Fig11} (b)), if the blue attracting center in Fig~\ref{Fig15} (a) was outside the sphere and the cylinder was attached on $S^2$ externally, the result would still be a torus.

Looking back at the natural phenomema happening on surfaces, an example is {\bf soap bubble splitting} during which a soap bubble splits into two smaller bubbles. This process is the 2-dimensional 1-surgery on $M=S^{2}$ shown in Fig~\ref{Fig15} (b). The orange attracting force in this case is the surface tension of each bubble that pulls molecules into the tightest possible groupings.

\subsection{Defining solid 2-dimensional surgery} \label{solid2D}

Most natural phenomena undergoing 2-dimensional surgery do not happen on surfaces but are three-dimensional. Therefore we introduce, also here, the notion of \textit{ solid 2-dimensional surgery}. There are two types of solid 2-dimensional surgery on the $3$-ball, $D^3$, analogous to the two types of 2-dimensional surgery.

The first one is the \textit{ solid 2-dimensional 0-surgery} which is the topological procedure of removing a solid cylinder homeomorphic to the product set $D^1 \times D^2$, $h(D^1 \times D^2)$ (such that the part $S^0 \times D^2$ of its boundary lies in the boundary of $D^3$) and taking the closure of the remaining manifold $D^3 \setminus h(D^1 \times D^2)$, which is a regular (or twisted) solid torus. See  Fig~\ref{Fig5} (b\textsubscript{1}) where the interior is supposed to be filled in.
The second type is the \textit{ solid 2-dimensional 1-surgery} which is the topological procedure of removing a solid cylinder homeomorphic to the product set $D^2 \times D^1$, $h(D^2 \times D^1)$, (such that the part $S^1 \times D^1$ of its boundary lies in the boundary of $D^3$) and taking the closure of the remaining manifold $D^3 \setminus h(D^2 \times D^1)$, which is two copies of $D^3$.  See  Fig~\ref{Fig5} (b\textsubscript{2}) where the interior is supposed to be filled in. Those processes are equivalent to performing 2-dimensional surgeries on the whole continuum of concentric spheres included in $D^3$. More precisely, and introducing at the same time dynamics, we define:

\begin{defn}\label{continuum2D} \rm  Start with the $3$-ball of radius 1 with polar layering: 
$$
D^3 = \cup_{0<r\leq 1} S^2_r \cup \{P\},
$$ 
where $r$ the radius of a 2-sphere and $P$ the limit point of the spheres, that is, the center of the ball. 
\textit{Solid 2-dimensional 0-surgery on $D^3$} is the  topological procedure shown in Fig~\ref{Fig16} (a): on all spheres $S^2_r$  colinear pairs of antipodal points are specified, \hl{all on the same diameter,} on which the same colinear attracting forces act. The poles have disc neighborhoods of analogous areas. Then, 2-dimensional 0-surgeries are performed on the whole continuum of the concentric spheres using the same homeomorphism $h$. Moreover, 2-dimensional 0-surgery on the limit point $P$ is defined to be the limit circle of the nested tori resulting from the continuum of 2-dimensional surgeries. That is,  the effect of \textit{2-dimensional 0-surgery on a point is the creation of a circle}. The process is characterized on one hand by the 1-dimensional core $L$ of the removed solid cylinder joining the antipodal points on the outer shell and intersecting each spherical layer in the two antipodal points and, on the other hand, by the homeomorphism $h$, resulting in the whole continuum of layered tori. The process can be viewed as drilling out a tunnel along $L$ according to $h$.  For a twisted embedding $h$, this agrees with our intuition that, for opening a hole, \textit{drilling with twisting} seems to be the easiest way.

\begin{figure}[!h]
\begin{center}
\centering
\captionsetup{justification=centering}
\includegraphics[width=13.5cm]{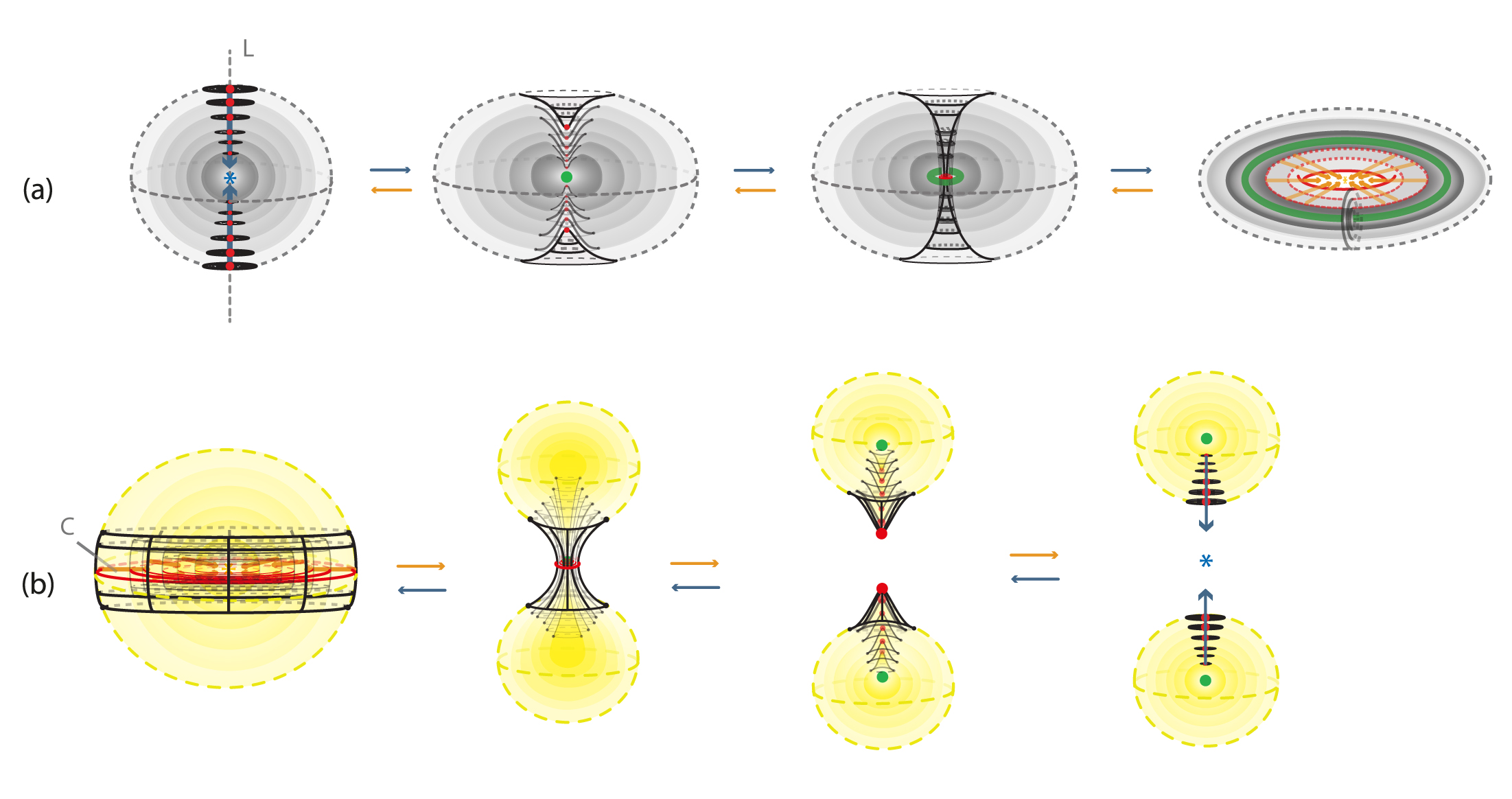}
\caption{{\bf Solid 2-dimensional surgery on the 3-ball.} (a) 2-dimensional 0-surgery with the standard embedding (b) 2-dimensional 1-surgery with the standard embedding.}
\label{Fig16}
\end{center}
\end{figure}

\noindent On the other hand, \textit{solid 2-dimensional 1-surgery on $D^3$} is the  topological procedure where: on all spheres $S^2_r$ nested annular peels of the solid annulus of analogous areas are specified and the same coplanar attracting forces act on all spheres, see  Fig~\ref{Fig16} (b). Then, 2-dimensional 1-surgeries are performed on the whole continuum of the concentric spheres using the same homeomorphism $h$. Moreover, 2-dimensional 1-surgery on the limit point $P$  is defined to be the two limit points of the nested pairs of 2-spheres resulting from  the continuum of 2-dimensional surgeries. That is, the effect of \textit{2-dimensional 1-surgery on a point is the creation of two new points}.
The process is characterized by the 2-dimensional central disc of the solid annulus and the homeomorphism $h$, and it can be viewed as squeezing the central disc $D$ or, equivalently, as pulling apart the upper and lower hemispheres with possible twists if $h$ is a twisted embedding. This agrees with our intuition that for cutting a solid object apart, \textit{pulling with twisting} seems to be the easiest way.

For both types, the above process is the same as: first removing the center $P$ from $D^3$, performing the  2-dimensional surgeries and then taking the closure of the resulting space. Namely we obtain:  
$$
\chi(D^3) := \cup_{0<r\leq 1}\chi(S^2_r) \cup \chi(P),
$$
which is a solid torus in the case of solid 2-dimensional 0-surgery and  two copies of $D^3$ in the case of solid 2-dimensional 1-surgery.
\end{defn}

\smallbreak
As seen in Fig~\ref{Fig16}, we also have the two dual solid 2-dimensional surgeries, which represent the reverse processes. As already mentioned in in Section~\ref{2D_FormalE1} the dual case of 2-dimensional 0-surgery is the 2-dimensional 1-surgery and vice versa. More precisely:

\begin{defn}  \rm The dual case of solid 2-dimensional 0-surgery on $D^3$ is the \textit{solid 2-dimensional 1-surgery on a solid torus $D^2 \times S^1$} whereby a solid cylinder $D^1 \times D^2$ filling the hole is added, such that the closure of the resulting manifold comprises one 3-ball $D^3$. This is the reverse process shown in  Fig~\ref{Fig16} (a) which results from the orange forces and attracting center. It only remain to define the solid 2-dimensional 1-surgery on the limit circle to be the limit point $P$ of the resulting surgeries. That is, the effect of \textit{ solid 2-dimensional 1-surgery on the core circle is that it collapses into one point}. The above process is the same as first removing the core circle from $D^2 \times S^1$, doing the 2-dimensional 1-surgeries on the nested tori, with the same coplanar acting forces, and then taking the closure of the resulting space. Given that the solid torus can be written as a union of nested tori together with the core circle: $D^2 \times S^1=(\cup_{0<r\leq 1}S^1_r \cup \{0\}) \times S^1$, the resulting manifold is 

$$
\chi^{-1}(D^2 \times S^1) := \cup_{0<r\leq 1}\chi^{-1}(S^1_r \times S^1) \cup \chi^{-1}(\{0\} \times S^1),
$$
which comprises one copy of $D^3$.

Further, the dual case of solid 2-dimensional 1-surgery on $D^3$ is the \textit{solid 2-dimensional 0-surgery on two $3$-balls $D^3$} whereby a solid cylinder $D^2 \times D^1$ joining the balls is added, such that the closure of the resulting manifold comprise of one 3-ball $D^3$. This is the reverse process shown in Fig~\ref{Fig16} (b) which results from the blue forces and attracting center. We only need to define the solid 2-dimensional 0-surgery on two limit points to be the limit point $P$ of the resulting surgeries. That is, as in solid 1-dimensional surgery, the effect of \textit{ solid 2-dimensional 0-surgery on two points is their merging into one point}. The above process is the same as first removing the centers from the $D^3 \times S^0$, doing the 2-dimensional 0-surgeries on the nested spheres, with the same colinear forces, and then taking the closure of the resulting space. The resulting manifold is 
$$
\chi^{-1}(D^3 \times S^0) := \cup_{0<r\leq 1}\chi^{-1}(S^2_r \times S^0) \cup \chi^{-1}(P \times S^0),
$$
which comprises one copy of $D^3$.
\end{defn}

\begin{rem} \rm The notions of 2-dimensional (resp. solid 2-dimensional) surgery, can be generalized from $S^2$ (resp. $D^3$) to a surface (resp. a handlebody) of genus $g$ creating a surface (resp. a handlebody) of genus $g+1$.
\end{rem} 

\subsection{Natural phenomena exhibiting solid 2-dimensional 0-surgery}\label{2D0}

Solid 2-dimensional 0-surgery is often present in natural phenomena where attracting forces between two poles are present, such as the formation of tornadoes, the formation of Falaco solitons, the formation of black holes, gene transfer in bacteria and drop coalescence. We shall  discuss these phenomena in some detail pinning down their exhibiting of topological surgery.

\begin{figure}[!h]
\begin{center}
\centering
\captionsetup{justification=centering}
\includegraphics[width=8cm]{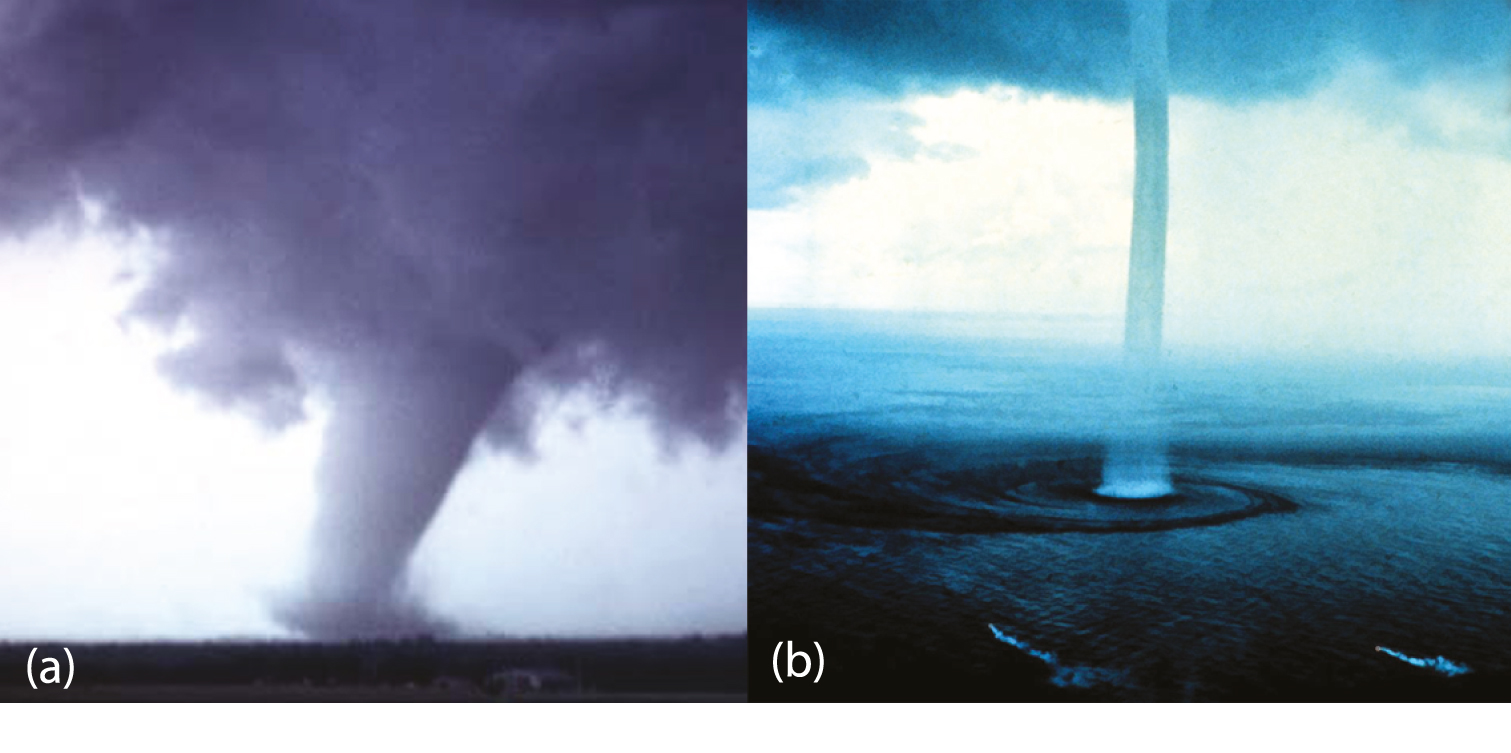}
\caption{{\bf (a) Funnel clouds drilling and tornado formation (b) waterspout.} An example of solid 2-dimensional 0-surgery.}
\label{Fig17}
\end{center}
\end{figure}

\smallbreak
Regarding {\bf tornadoes}: except from their shape (see Fig~\ref{Fig17}) which fits the cylinder $D^1\times S^1$ that gets attached in the definition of 2-dimensional 0-surgery, the process by which they are formed also follows the dynamics introduced in Section~\ref{solid2D}. Namely, if certain meteorological conditions are met, an attracting force between the cloud and the earth beneath is created and funnel-shaped clouds start descending toward the ground. Once they reach it, they become tornadoes. In analogy to solid 2-dimensional 0-surgery, first the poles are chosen, one on the tip of the cloud and the other on the ground, and they seem to be joined through an invisible line. Then, starting from the first point, the wind revolves in a  helicoidal motion toward the second point, resembling `hole drilling' along the line until the hole is drilled. Therefore, tornado formation undergoes the process of solid 2-dimensional 0-surgery with a twisted embedding, as in Fig~\ref{Fig14} (b). The initial manifold can be considered as $M=D^3 \times S^0$, that is, one 3-ball on the cloud and one on the ground. Note that in this realization of solid 2-dimensional 0-surgery, the attracting center coincides with the ground and we only see helicoidal motion in one direction.

\begin{figure}[!h]
\begin{center}
\centering
\captionsetup{justification=centering}
\includegraphics[width=6cm]{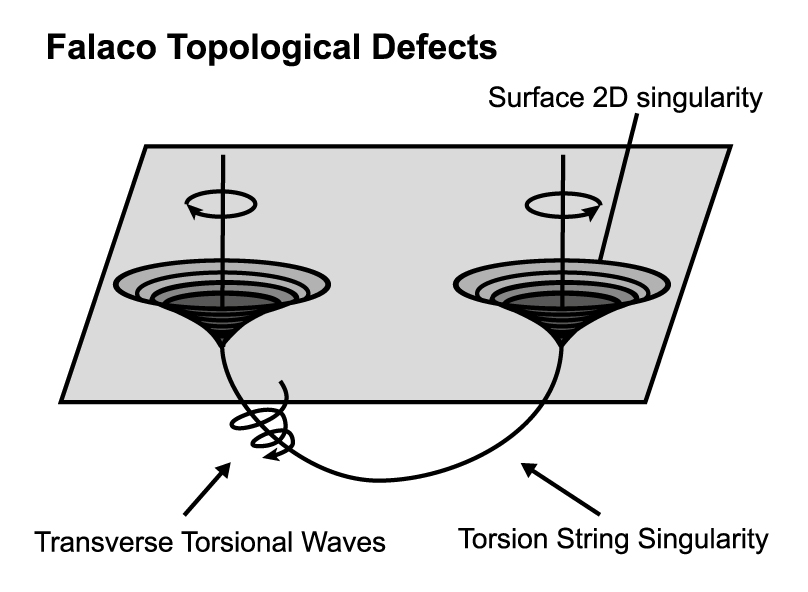}
\caption{{\bf Pairs of Falaco solitons.} An example of solid 2-dimensional 0-surgery.}
\label{Fig18}
\end{center}
\end{figure}

Another natural phenomenon exhibiting solid 2-dimensional 0-surgery is the formation of {\bf Falaco solitons}, see Fig~\ref{Fig18} (for photos of pairs of Falaco solitons in a swimming pool, see \cite{Ki}). Note that the term `Falaco Soliton' appears in 2001 in \cite{Ki2}. Each Falaco Soliton consists of a pair of locally unstable but globally stabilized contra-rotating identations in the water-air discontinuity surface of a swimming pool. These pairs of singular surfaces (poles) are connected by means of a stabilizing thread. This thread corresponds to the `invisible line' mentioned in the process of tornado formation which is visible in this case. The two poles get connected and their rotation propagates below the water surface along the joining thread and the tubular neighborhood around it. This process is a solid 2-dimensional 0-surgery with a twisted embedding (see Fig~\ref{Fig14} (b)) where the initial manifold is the water contained in the volume of the pool where the process happens, which is homeomorphic to a 3-ball, that is $M=D^3$. Two differences compared to tornadoes are: here the helicoidal motion is present in both poles and the attracting center is not located on the ground but between the poles, on the topological thread joining them. 

It is also worth mentioning that the creation of Falaco solitons is immediate and does not allow us to see whether the transitions of the 2-dimensional 0-surgery shown in Fig~\ref{Fig14} (b) are followed or not. However, these dynamics are certainly visible during the annihilation of Falaco solitons. Namely, when the topological thread joining the poles is cut, the tube tears apart and slowly degenerates to the poles until they both stops spinning and vanish. Therefore, the continuity of our dynamic model is clearly present during the reverse process which corresponds to a solid 2-dimensional 1-surgery on a pair of Falaco solitons, that is, a solid torus $D^2 \times S^1$ degenerating into a still swimming pool $D^3$.
 
Note that it is conjectured in \cite{Ki} that the coherent topological features of the Falaco solitons and, by extension, the process of solid 2-dimensional 0-surgery appear in both macroscopic level (for example in the Wheeler's wormholes) and microscopic level (for example in the spin pairing mechanism in the microscopic Fermi surface). For more details see \cite{Ki}.

\begin{figure}[!h]
\begin{center}
\centering
\captionsetup{justification=centering}
\includegraphics[width=9cm]{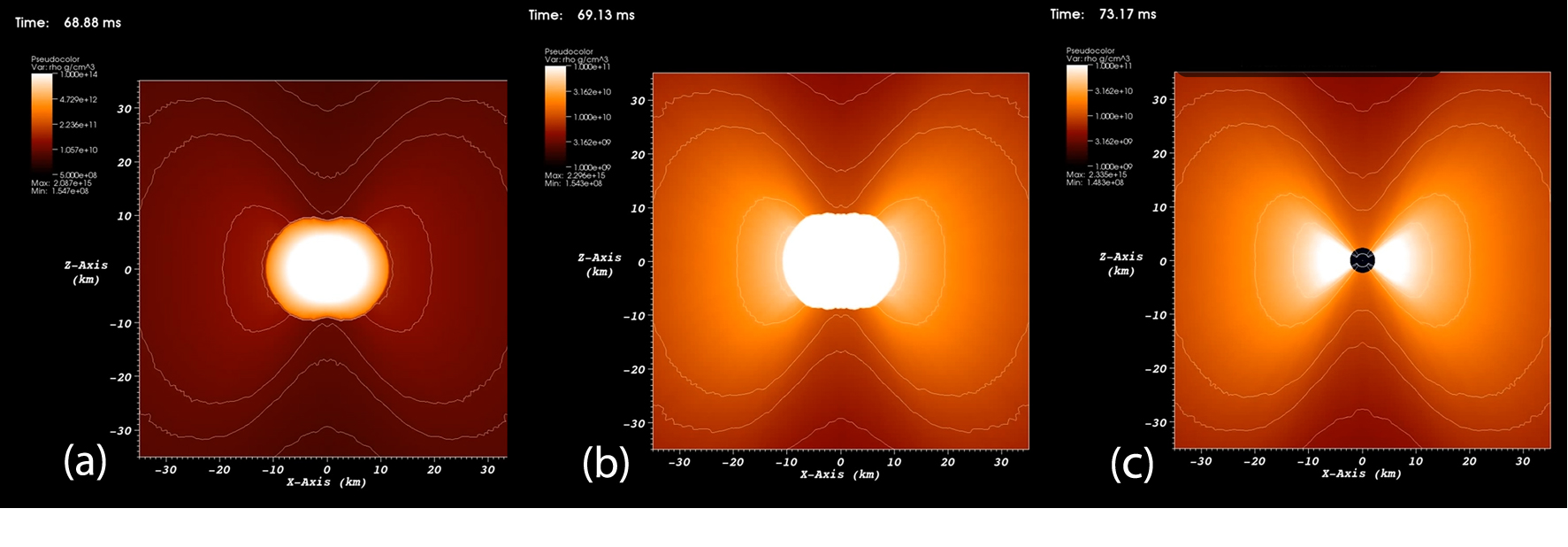}
\caption{{\bf The formation of a black hole.} An example of solid 2-dimensional 0-surgery.}
\label{Fig19}
\end{center}
\end{figure}

Another phenomenon undergoing solid 2-dimensional 0-surgery is the formation of a {\bf black hole}. Most black holes form from the remnants of a large star that dies in a supernova explosion and have a gravitational field so strong that not even light can escape. In the simulation of a black hole formation (see \cite{Ott}), the density distribution at the core of a collapsing massive star is shown. In Fig~\ref{Fig19} matter performs solid 2-dimensional 0-surgery as it collapses into a black hole. Matter collapses at the center of attraction of the initial manifold $M=D^3$ creating the singularity, that is, the center of the black hole, which is surrounded by the toroidal accretion disc (shown in white in Fig~\ref{Fig19} (c)).

\begin{figure}[!h]
\begin{center}
\centering
\captionsetup{justification=centering}
\includegraphics[width=5cm]{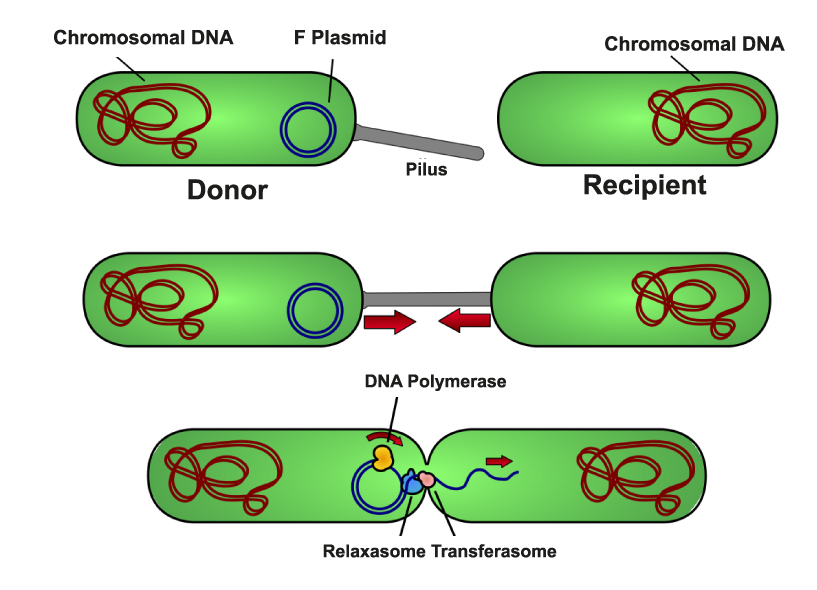}
\caption{{\bf Gene transfer in bacteria.} An example of solid 2-dimensional 0-surgery.}
\label{Fig20}
\end{center}
\end{figure}

Solid 2-dimensional 0-surgery is also found in the mechanism of {\bf gene transfer in bacteria}. See Fig~\ref{Fig20} (also, for description and instructive illustrations see \cite{HHGRSV}). The donor cell produces a connecting tube called a `pilus' which attaches to the recipient cell, brings the two cells together and transfers the donor's DNA. This process is similar to the one shown earlier in Fig~\ref{Fig16} (b) as two copies of $D^3$ merge into one, but here the attracting center is located on the recipient cell. This process is a solid 2-dimensional 0-surgery on two 3-balls $M=D^3 \times S^0$. 


Finally, {\bf drop coalescence} is the merging of two dispersed drops into one. As gene transfer in bacteria, this process is also a solid 2-dimensional 0-surgery on two 3-balls $M=D^3 \times S^0$, see Fig~\ref{Fig16} (b). The process of drop coalescence also exhibits the forces of our model. Namely, the surfaces of two drops must be in contact for coalescence to occur. This surface contact is dependent on both the van der Waals attraction and the surface repulsion forces between two drops. When the van der Waals forces cause rupture of the film, the two surface films are able to fuse  together, an event more likely to occur in areas where the surface film is weak. The liquid inside each drop is now in direct contact, and the two drops are able to merge into one.

\begin{rem}\label{DSGenerality} \rm
Although in this Section some natural processes were viewed as a solid 2-dimensional topological surgery on $M=D^3 \times S^0$, we could also consider the initial manifold as being a 3-ball surrounding the phenomena and view it as a surgery on $M=D^3$. Concerning the process of tornado formation, this approach also has a physical meaning. Namely, as the process is triggered by the difference in the conditions of the lower and upper atmosphere, the initial manifold can be considered as the 3-ball containing this air cycle. 
\end{rem} 

\subsection{Natural phenomena exhibiting solid 2-dimensional 1-surgery}\label{2D1}

As already mentioned, the collapsing of the central disc of the sphere caused by the orange attracting forces in Fig~\ref{Fig16} (b) can also be caused by pulling apart the upper and lower hemispheres of the 3-ball $D^3$, that is, the causal forces can also be repelling. For example, during fracture of metal specimens under tensile forces, solid 2-dimensional 1-surgery is caused by forces that pull apart each end of the specimen. On the other hand, in the biological process of mitosis, both attracting and repelling forces forces are present.

\begin{figure}[!h]
\begin{center}
\centering
\captionsetup{justification=centering}
\includegraphics[width=3cm]{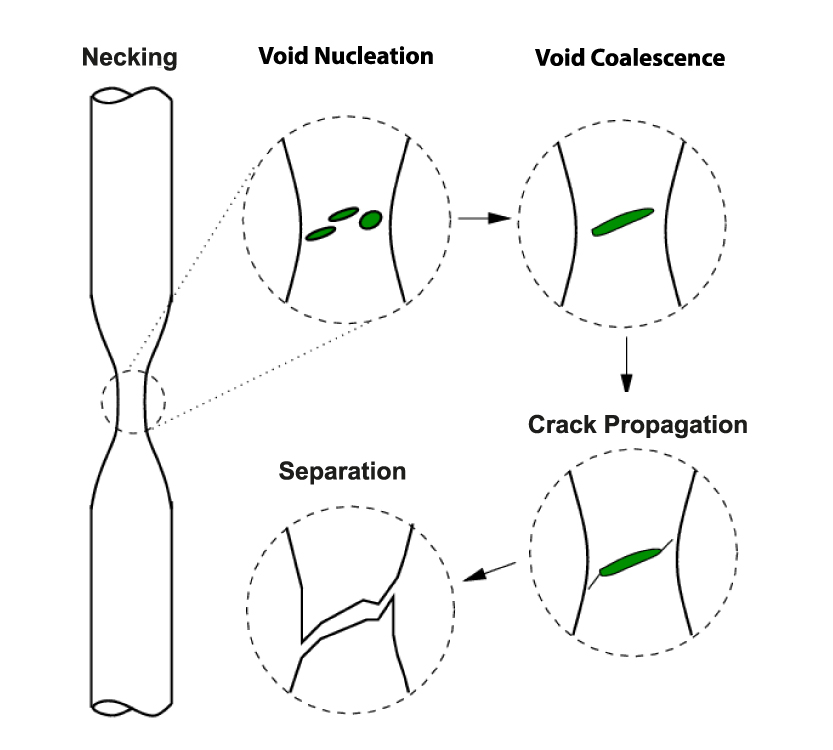}
\caption{{\bf Tension and the necking phenomenon.} An example of solid 2-dimensional 1-surgery.}
\label{Fig21}
\end{center}
\end{figure}

When the tension applied on metal specimens by tensile forces results in {\bf necking} and then {\bf fracture}, the process exhibits solid 2-dimensional 1-surgery. More precisely, in experiments in mechanics, tensile forces (or loading) are applied on a cylindrical specimen made of dactyle material (steel, aluminium, etc.). Up to some critical value of the force the deformation is homogeneous (the cross-sections have the same area). At the critical value the deformation is localized within a very small area where the cross-section is reduced drastically, while the sections of the remaining portions increase slightly. This is the `necking phenomenon'. Shortly after, the specimen is fractured (view \cite{UV} for details). In Fig~\ref{Fig21} are the the basic steps of the process: void formation, void coalescence (also known as crack formation), crack propagation, and failure. Here, the process is not as smooth as our theoretical model and the tensile forces applied on the specimen are equivalent to repelling forces. The specimen is homeomorphic to the sphere shown in Fig~\ref{Fig16} (b) hence the initial manifold is $M=D^3$.


\begin{figure}[!h]
\begin{center}
\centering
\captionsetup{justification=centering}
\includegraphics[width=7cm]{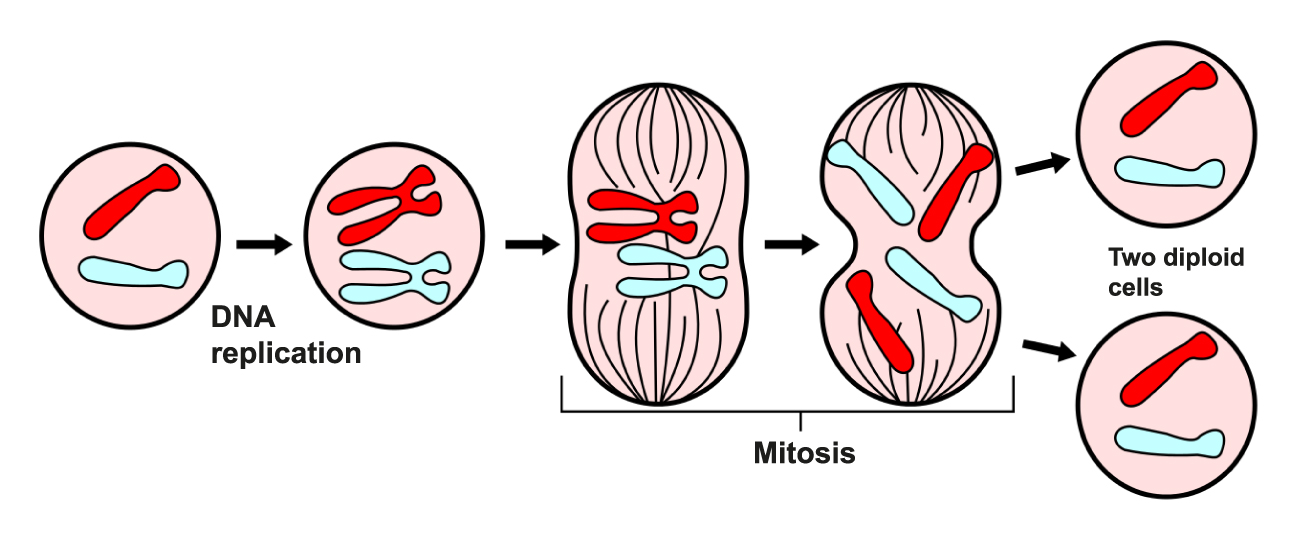}
\caption{{\bf The process of mitosis.} An example of solid 2-dimensional 1-surgery.}
\label{Fig22}
\end{center}
\end{figure}

Solid 2-dimensional 1-surgery on $M=D^3$ also happens in the biological process of {\bf mitosis}, where a cell splits into two new cells. See Fig~\ref{Fig22} (for description and instructive illustrations see for example \cite{KeFa}). We will see that both aforementioned forces are present here. During mitosis, the chromosomes, which have already duplicated, condense and attach to fibers that pull one copy of each chromosome to opposite sides of the cell (this pulling is equivalent to repelling forces). The cell pinches in the middle and then divides by cytokinesis. The structure that accomplishes cytokinesis is the contractile ring, a dynamic assembly of filaments and proteins which assembles just beneath the plasma membrane and contracts to constrict the cell into two (this contraction is equivalent to attracting forces). In the end, two genetically-identical daughter cells are produced.

 
\begin{rem} \rm 
It is worth noting that the splitting of the cell into two coincide with the fact that 2-dimensional 1-surgery on a point is the creation of two new points (see Definition ~\ref{continuum2D}).
\end{rem}

\section{Connecting 1- and 2-dimensional surgeries}\label{Connecting}

\begin{figure}[!h]
\begin{center}
\centering
\captionsetup{justification=centering}
\includegraphics[width=13.5cm]{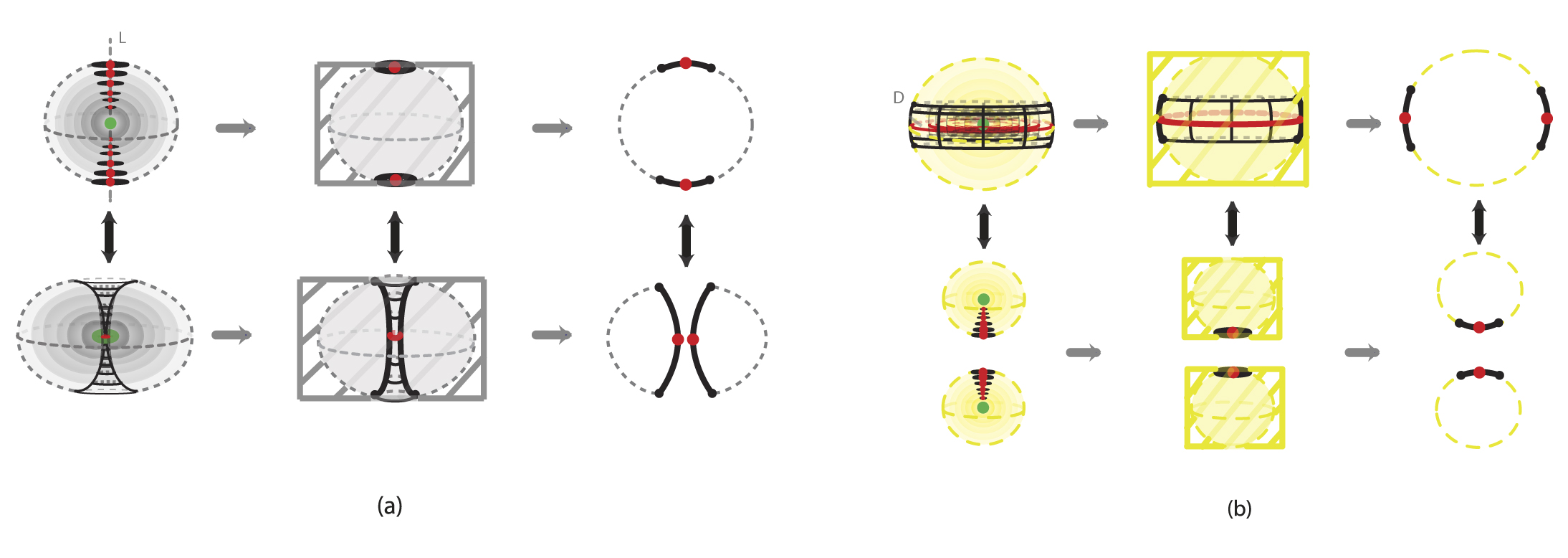}
\caption{{\bf Connecting low-dimensional surgeries.} From left to right we pass from solid 2-dimensional to 2-dimensional to 1-dimensional surgery for (a) 0-surgeries and (b) 1-surgeries}
\label{Fig23}
\end{center}
\end{figure}

As illustrated in Fig~\ref{Fig23}, a 1-dimensional surgery is a planar cross-section of the corresponding 2-dimensional surgery which, in turn, is a spherical/toroidal crossection of the corresponding type of solid 2-dimensional surgery. This is true for both 1 or 0-surgeries (see Fig~\ref{Fig23} (a) and (b) respectively).

On the left-hand top and bottom pictures of Fig~\ref{Fig23} (a) and (b) we see the initial and final stage of solid 2-dimensional surgery. Taking the intersection with the boundary of the 3-ball $D^3$ we pass to the middle pictures where we see the the initial and final pictures of 2-dimensional surgery. Taking finally the intersection with a meridional plane gives rise to the initial and final stages of 1-dimensional surgery (rightmost illustrations). The above concerns 0-surgeries in Fig~\ref{Fig23} (a) and 1-surgeries in Fig~\ref{Fig23} (b).

\begin{figure}[!h]
\begin{center}
\centering
\captionsetup{justification=centering}
\includegraphics[width=13.5cm]{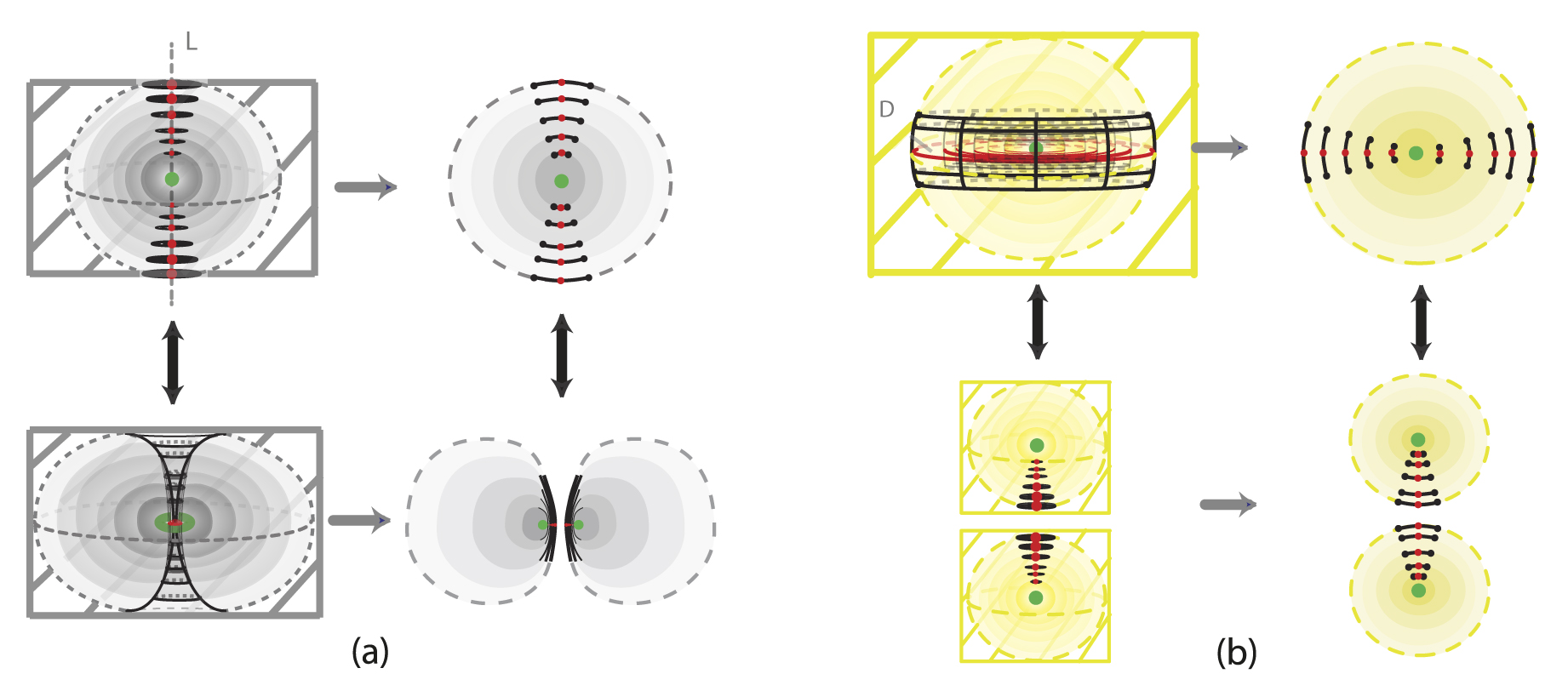}
\caption{{\bf Connecting solid surgeries.} From left to right we pass from solid 2-dimensional to solid 1-dimensional surgery via a cutting meridional plane.}
\label{Fig24}
\end{center}
\end{figure}

Furthermore, in Fig~\ref{Fig24} we see the relation between solid surgeries in dimensions 2 and 1. \hl{Namely, we see that solid $1$-dimensional surgery is a cross-section of solid $2$-dimensional surgery via a cutting meridional plane. In particular, we see that} solid $2$-dimensional 0-surgery on the central point of the spherical nesting results in the central circle of the toroidal nesting. This circle has two intersecting points with the plane which are the result of solid $1$-dimensional $0$-surgery on the central point, see Fig~\ref{Fig24} (a). On the other hand, both solid $2$-dimensional $1$-surgery and solid $1$-dimensional $0$-surgery on the central point creates two points, see Fig~\ref{Fig24} (b).

\hl{It is worth adding that all types of $1$- and $2$-dimensional surgeries can be also connected via appropriate horizontal and vertical rotations in the $3$-space as demonstrated in} \cite{SS2}.

\section{The ambient space \texorpdfstring{$S^3$} {S3}}\label{S3}

All natural phenomena exhibiting surgery (1- or 2-dimensional, solid or usual) take place in the ambient 3-space. As we will see in the next section, the ambient space can play an important role in the process of surgery. By {\it 3-space} we mean here the compactification of ${\mathbb R}^3$ which is the 3-sphere  $S^3$. This choice, as opposed to ${\mathbb R}^3$, takes advantage of the duality of the descriptions of $S^3$. In this section we present the three most common descriptions of $S^3$ in which this duality is apparent and which will set the ground for defining the notion of embedded surgery in $S^3$. Beyond that, we also demonstrate how the descriptions are interrelated via solid 2-dimensional 0-surgery which, due to the duality of the dimensions, takes place in both the initial 3-ball and its complement.

\subsection{Descriptions of \texorpdfstring{$S^3$} {S3}}\label{decsr}

In dimension 3, the simplest c.c.o. 3-manifolds are: the 3-sphere $S^3$ and the lens spaces $L(p,q)$. In this paper however, we will focus on  $S^3$. We start by recalling its three most common descriptions: 
\subsubsection{Via \texorpdfstring{${\mathbb R}^3$} {R3}}\label{ViaR3}
\begin{figure}[!h]
\begin{center}
\centering
\captionsetup{justification=centering}
\includegraphics[width=8cm]{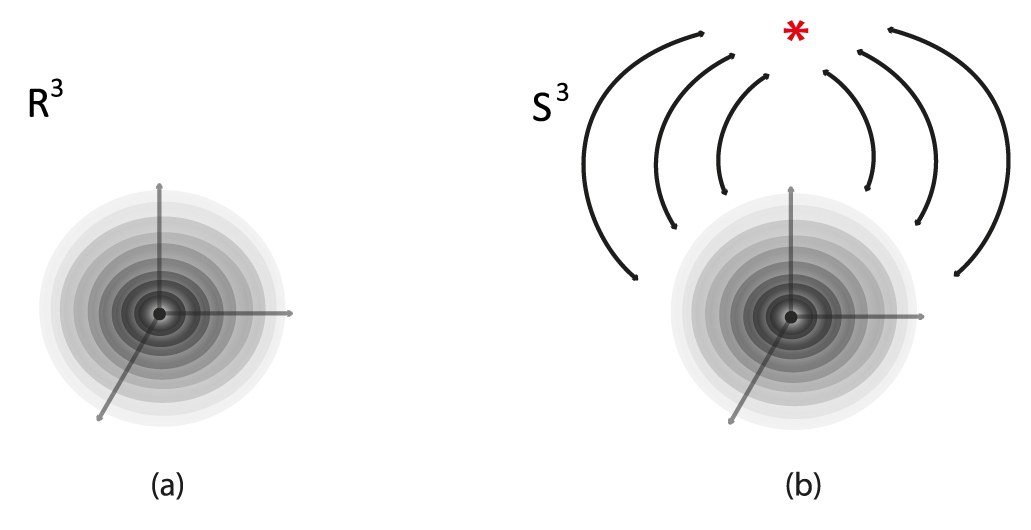}
\caption{{\bf $S^3$ is the compactification of ${\mathbb R}^3$.}}
\label{Fig25}
\end{center}
\end{figure}

$S^3$ can be viewed as ${\mathbb R}^3$ with all points at infinity compactified to one single point: $S^3 = {\mathbb R}^3 \cup \{\infty\}$. See Fig~\ref{Fig25} (b). 
 ${\mathbb R}^3$ can be viewed as an unbounded continuum of nested 2-spheres centered at the origin, together with the point at the origin, see Fig~\ref{Fig25} (a),  and also as the de-compactification of $S^3$.  So, $S^3$ minus the point at the origin and the point at infinity can be viewed as a continuous nesting of 2-spheres.

\subsubsection{Via two 3-balls}\label{ViaTwoBalls}
\begin{figure}[!h]
\begin{center}
\centering
\captionsetup{justification=centering}
\includegraphics[width=9cm]{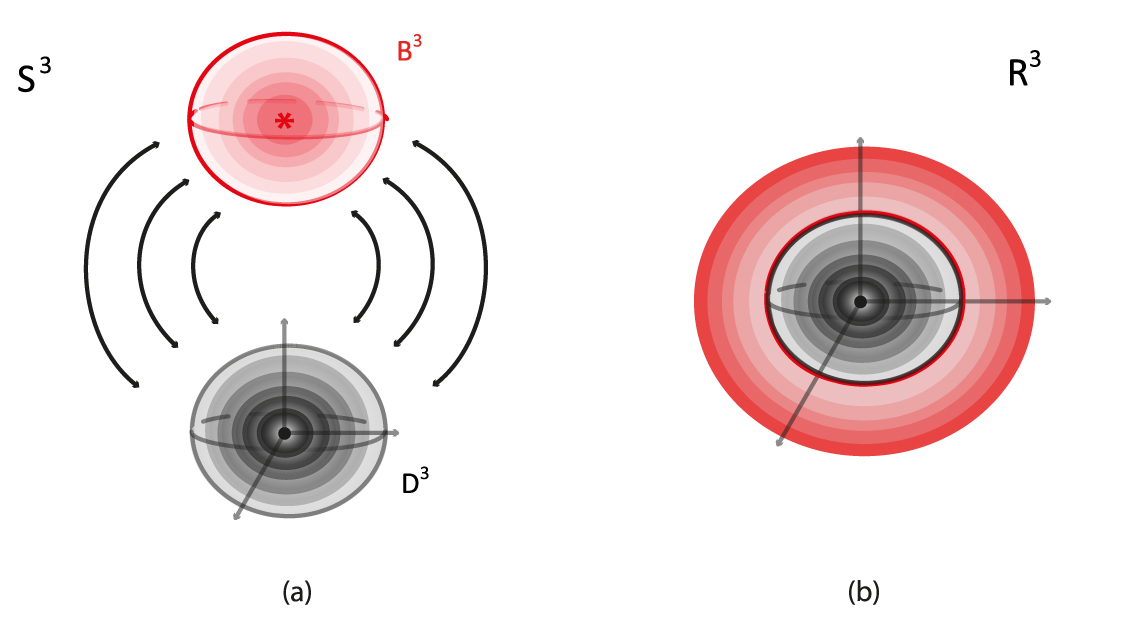}
\caption{{\bf $S^3$ is the result of gluing two 3-balls.}}
\label{Fig26}
\end{center}
\end{figure}

$S^3$ can be viewed as  the union of two $3$-balls: $S^3 = B^3 \cup D^3$, see Fig~\ref{Fig26} (a). This second description of $S^3$ is clearly related to the first one, since a  (closed) neighborhood of the point at infinity can stand for one of the two $3$-balls.  Note that, when removing the point at infinity in Fig~\ref{Fig26}  (a) we can see the concentric spheres of the 3-ball $B^3$ (in red) wrapping around the concentric spheres of the 3-ball $D^3$, see Fig~\ref{Fig26}     
 (b). This is another way of viewing ${\mathbb R}^3$ as the de-compactification of $S^3$.  This picture is the analogue of the stereographic projection of $S^2$ on the plane ${\mathbb R}^2$, whereby the projections of the concentric circles of the south hemisphere together with the projections of the concentric circles of the north hemisphere form the well-known polar description of ${\mathbb R}^2$ with the unbounded continuum of concentric circles.

\subsubsection{Via two solid tori}\label{ViaTwoTori}
\begin{figure}[!h]
\begin{center}
\centering
\captionsetup{justification=centering}
\includegraphics[width=13.5cm]{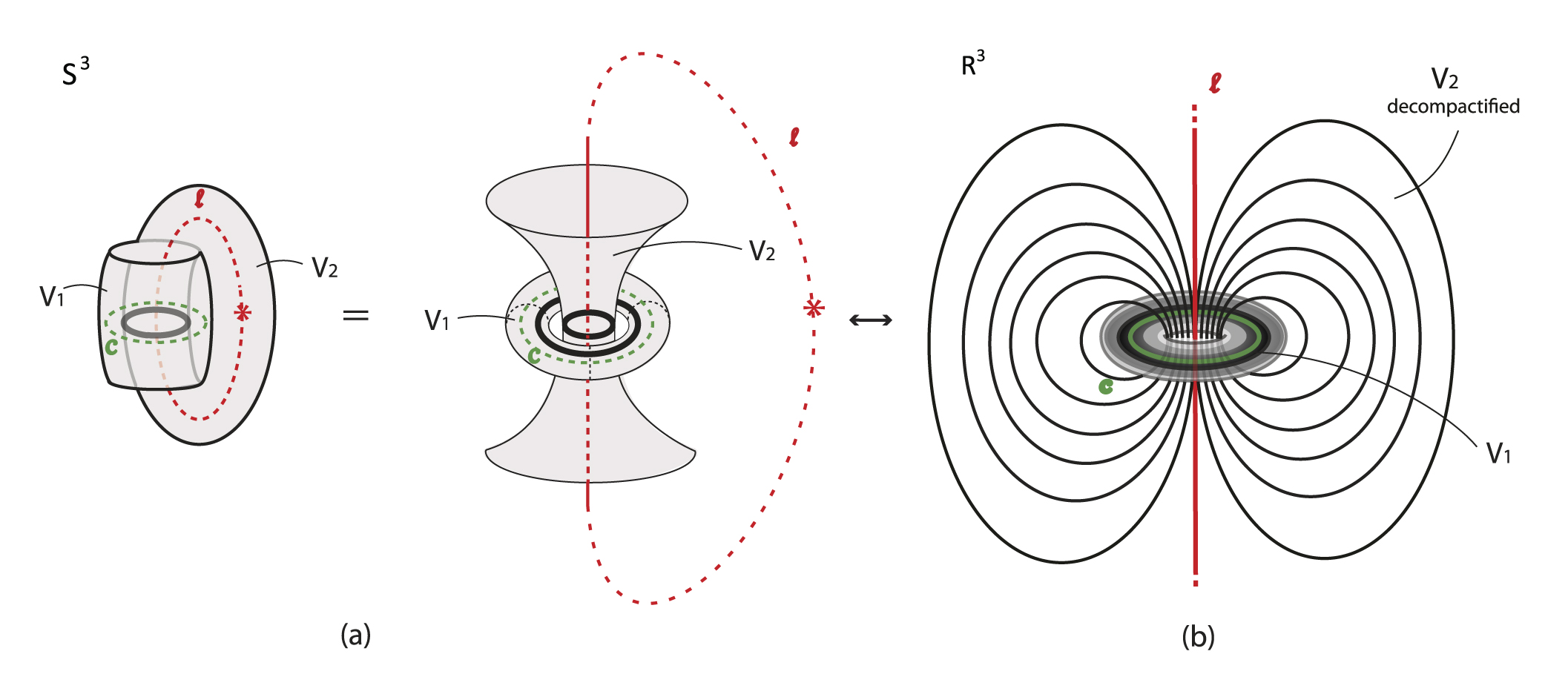}
\caption{{\bf (a) $S^3$ as a union of two solid tori (b) De-compactificated view. }}
\label{Fig27}
\end{center}
\end{figure}

The third well-known representation of $S^3$ is as the union of two solid tori, $S^3 = V_1\cup_\vartheta V_2$, via the torus homeomorphism $\vartheta$ along the common boundary.  $\vartheta$ maps a meridian of $V_2$ to a longitude of $V_1$ which has linking number zero with the core curve $c$ of $V_1$.
The illustration in Fig~\ref{Fig27} (a) gives an idea of this splitting of $S^3$. In the figure, the core curve of $V_1$ is in dashed green. So, the complement of a solid torus $V_1$ in $S^3$ is another solid torus $V_2$ whose core curve $l$ (in dashed red) may be assumed to pass by the point at infinity. Note that, $S^3$ minus the core curves $c$ and $l$ of  $V_1$ and $V_2$ (the green and red curves in Fig~\ref{Fig27}) can be  viewed as a continuum of nested tori.   
When removing the point at infinity in the representation of $S^3$ as a union of two solid tori, the core of the solid torus $V_2$ becomes an infinite line $l$ and the nested tori of  $V_2$ can now be seen wrapping around the nested tori of  $V_1$. See Fig~\ref{Fig27} (b).  Therefore, ${\mathbb R}^3$ can be viewed as an unbounded continuum of nested tori, together with the core curve $c$ of  $V_1$  and the infinite line  $l$. This line $l$ joins pairs of antipodal points of all concentric spheres of the first description. Note that in the nested spheres description (Fig~\ref{Fig25}) the  line $l$ pierces all spheres while in the nested tori description the  line $l$ is the `untouched' limit circle of all tori.

\begin{rem}\label{Glatz} \rm
It is worth observing the resemblance of Fig~\ref{Fig27} (b) with the well-known representation of the {\bf Earth magnetic field}. A numerical simulation of the Earth magnetic field via the Glatzmaier-Roberts geodynamo model was made in \cite{Glatzmaier}, see Fig~\ref{Fig28}. The magnetic field lines are lying on nested tori and comprise a visualization of the decompactified view of $S^3$ as two solid tori.

\begin{figure}[!h]
\begin{center}
\centering
\captionsetup{justification=centering}
\includegraphics[width=4cm]{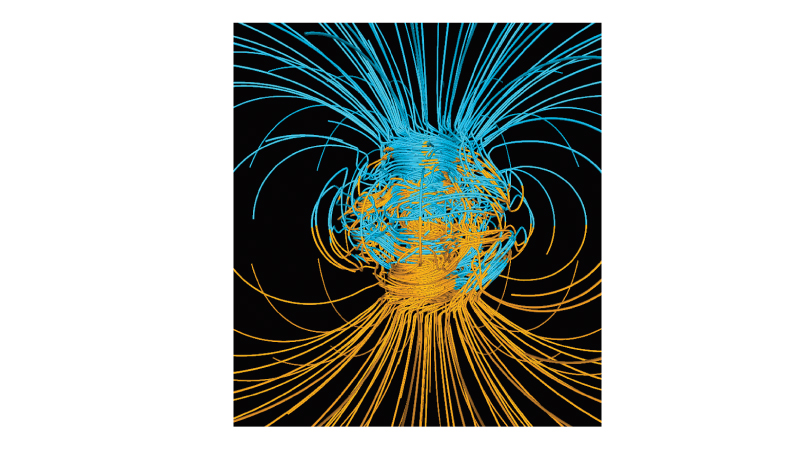}
\caption{{\bf The Earth magnetic field as a decompactified view of $S^3$ as two solid tori}}
\label{Fig28}
\end{center}
\end{figure}

\end{rem} 

\begin{rem}\label{Hopf} \rm
It is also worth mentioning that another way to visualize $S^3$ as two solid tori is the {\bf Hopf fibration}, which is a map of $S^3$ into $S^2$. The parallels of $S^2$ correspond to the nested tori of $S^3$, the north pole of $S^2$ correspond to the core curve $l$ of $V_2$ while the south pole of $S^2$ corresponds to the core curve $c$ of $V_1$. An insightful animation of the Hopf fibration can be found in \cite{NilesJ}.
\end{rem}

\subsection{Connecting the descriptions of \texorpdfstring{$S^3$} {S3} } \label{condecsr}

\subsubsection{Via corking} \label{corking}
\begin{figure}[!h]
\begin{center}
\centering
\captionsetup{justification=centering}
\includegraphics[width=9cm]{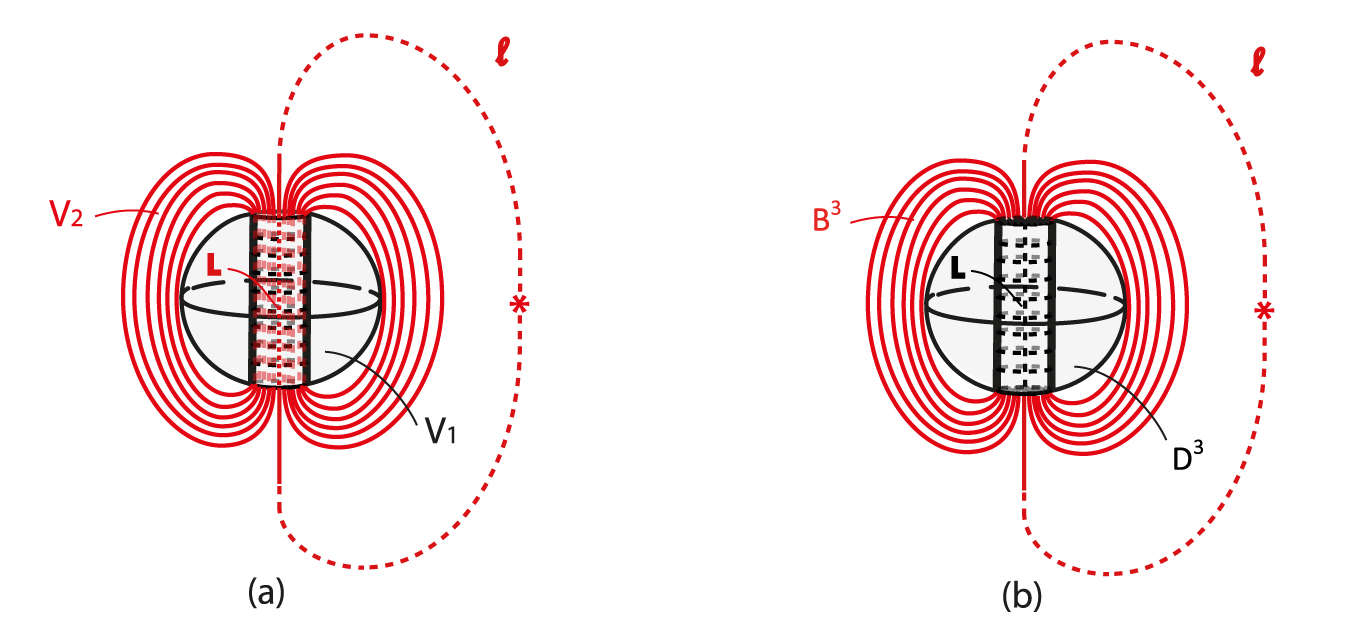}
\caption{{\bf Passing from (a) $S^3$ as two \hl{solid} tori  to (b) $S^3$ as two balls.}}
\label{Fig29}
\end{center}
\end{figure}

The connection between the first two descriptions of $S^3$ was already discussed in previous Section. The third  description is a bit harder to connect with the first two. We shall do this here. A way to see this connection is the following. 
Consider the description of  $S^3$ as the union of two 3-balls, $B^3$ and $D^3$ (Fig~\ref{Fig26}). 
Combining with the third description of  $S^3$ (Fig~\ref{Fig27}) we notice that both 3-balls are pierced by  the core curve $l$ of the solid torus $V_2$. 
 Therefore, $D^3$ can be viewed as the solid torus $V_1$ to which  a solid cylinder $D^1\times D^2$  is attached  via the homeomorphism $\vartheta$: 
$$
D^3 = V_1\cup_\vartheta (D^1\times D^2).
$$
This solid cylinder is part of the solid torus $V_2$,  a `cork'  filling the hole of $V_1$. Its core curve is an arc $L$, part of the core curve $l$ of $V_2$. View  Fig~\ref{Fig29}. The second ball $B^3$ (Fig~\ref{Fig26}) can be viewed as the remaining of $V_2$ after removing the `cork' $D^1\times D^2$: 
$$
B^3 = \overline{ V_2 \setminus_\vartheta (D^1\times D^2)}.
$$
In other words the solid torus $V_2$ is cut into two solid cylinders, one comprising the `cork' of $V_1$ and the other comprising the 3-ball $B^3$.

\begin{rem}\label{truncate} \rm
If we remove a whole neighborhood $B^3$ of the point at infinity and focus on the remaining  3-ball $D^3$, the line $l$ of the previous picture is truncated to the arc $L$ and the solid cylinder $V_2$ is truncated to  the `cork' of $D^3$.   
\end{rem}

\subsubsection{Via surgery}\label{ConSurg}

\begin{figure}[!h]
\begin{center}
\centering
\captionsetup{justification=centering}
\includegraphics[width=13.5cm]{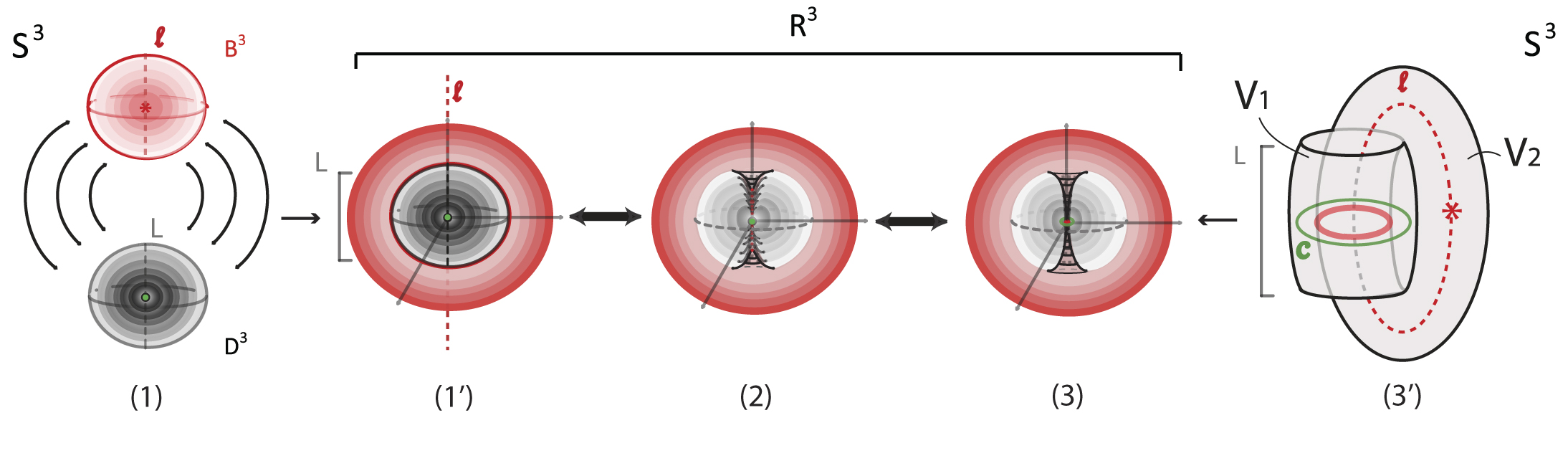}
\caption{{\bf Passing from the two balls description to the two \hl{solid} tori description of $S^3$ via solid 2-dimensional 0-surgery }}
\label{Fig30}
\end{center}
\end{figure}

We will now examine how we can pass from the two-ball description to the two-tori description of $S^3$ via solid 2-dimensional 0-surgery. We start with two points that have a distance $L$ between them. Let $M=D^3$ be the solid ball having arc $L$ as a diameter. We define this 3-ball as the `truncated' space on which we will focus. When the center of $D^3$ becomes attracting, forces are induced on the two points of $D^3$ and solid 2-dimensional 0-surgery is initiated. The complement space is the other solid ball $B^3$ containing the point at infinity, recall Fig~\ref{Fig26}. This joining arc $L$ is seen as part of a simple closed curve $l$ passing by the point at infinity. In Fig~\ref{Fig30} (1) this is shown in $S^3$ while Fig~\ref{Fig30} (1$'$) shows the corresponding decompactified view in ${\mathbb R}^3$.  

In Fig~\ref{Fig30} (2), we see the `drilling' along $L$ as a result of the attracting forces. This is exactly the same process as in Fig~\ref{Fig16} if we restrict it to $D^3$. But since we are in $S^3$, the complement space $B^3$ participates in the process and, in fact, it is also undergoing solid 2-dimensional 0-surgery. In Fig~\ref{Fig30} (3), we can see that, as surgery transforms the solid ball $D^3$ into the solid torus $V_1$, $B^3$ is transformed into $V_2$. That is, the nesting of concentric spheres of $D^3$ (respectively $B^3$) is transformed into the nesting of concentric tori in the interior of $V_1$ (respectively $V_2$). This is a double surgery with one attracting center which is inside the first 3-ball $D^3$ (in grey) and outside the second 3-ball $B^3$ (in red). By Definition~\ref{continuum2D}, the point at the origin (in green) turns into the core curve $c$ of $V_1$ (in green). Fig~\ref{Fig30} (3) is exactly the decompactified view of $S^3$ as two \hl{solid} tori as shown in Fig~\ref{Fig27} (b) while Fig~\ref{Fig30} (3$'$) is the corresponding view in $S^3$ as shown in Fig~\ref{Fig27} (1). 

Fig~\ref{Fig30} shows that one can pass from the second description of $S^3$ to the third by performing solid 2-dimensional 0-surgery (with the standard embedding homeomorphism) along the arc $L$ of $D^3$. It is worth mentioning that this connection between the descriptions of $S^3$ and solid 2-dimensional 0-surgery is a dynamic way to visualize the connection established in Section \ref{corking}.

\section{Embedding surgery in \texorpdfstring{$S^3$} {S3}}\label{es3sp}

In this section we define the notion of \textit{embedded surgery in 3-space}. As we will see, when embedded surgery occurs, depending on the dimension of the manifold, the ambient space either leaves `room' for the initial manifold to assume a more complicated configuration or it participates more actively in the process.

\subsection{Defining embedded $m$-dimensional $n$-surgery}\label{es3sp2}

We will now concretely define the notion of embedded $m$-dimensional $n$-surgery in some sphere $S^d$ and we will then focus on the case $d=3$.

\begin{defn} \rm   \textit{An embedded $m$-dimensional $n$-surgery} is a $m$-dimensional $n$-surgery where the initial manifold is an $m$-embedding $e:M \hookrightarrow S^d, d\geq m$ of some $m$-manifold $M$. Namely, according to Definition~\ref{surgery}: 
\[M' = \chi(e(M)) = \overline{e(M)\setminus h(S^n\times D^{m-n})} \cup_{h|_{S^n\times S^{m-n-1}}} D^{n+1}\times S^{m-n-1}. \]
\end{defn}

From now on we fix $d=3$. Embedding surgery allows to view it as a process happening in 3-space instead of abstractly. In the case of embedded 1-dimensional 0-surgery on a circle $M=S^1$, the ambient space gives enough `room' for the initial 1-manifold to become any type of knot. Hence, embedding {\it allows the initial manifold to assume a more complicated homeomorphic configuration}. This will be analyzed further in Section \ref{Embedded1D}.

Passing now to 2-dimensional surgeries, let us first note that embedded 2-dimensional surgery is often used a theoretical tool in various proofs in low dimensional topology. Further, an embedding of a sphere $M=S^2$ in $S^3$ presents no knotting because knots require embeddings of codimension 2. However, in this case the ambient space plays a different role. Namely, embedding 2-dimension surgeries {\it allows the complementary space of the initial manifold to participate actively in the process}. Indeed, while some natural phenomena undergoing surgery can be viewed as `local', in the sense that they can be considered independently from the surrounding space, some others are intrinsically related to the surrounding space. This relation can be both \textit{causal}, in the sense that the ambient space is involved in the triggering of the forces causing surgery, and \textit{consequential}, in the sense that the forces causing surgery, can have an impact on the ambient space in which they take place. This will be analyzed in Sections \ref{Embedded2D0} and \ref{Embedded2D1}.

\subsection{Embedded 1-dimensional 0-surgery and related phenomena} \label{Embedded1D}

We will now get back to site-specific {\bf DNA recombination} (see Section \ref{1DPhenomena}), in order to better define this type of surgery. As seen in this process (recall Fig~\ref{Fig8}) the initial manifold of 1-dimensional 0-surgery can be a knot, in other words, an embedding of the circle $M=S^1$ in 3-space. We therefore introduce the notion of \textit{ embedded 1-dimensional 0-surgery} whereby the initial manifold $M$ is embedded in the 3-space. This notion allows the topological modeling of phenomena with more complicated initial 1-manifolds. As mentioned, for our purposes, we will consider $S^3$ as our standard 3-space. For details on the descriptions of $S^3$, see Section~\ref{decsr}. Since a knot is by definition an embedding of $M=S^1$ in $S^3$ or ${\mathbb R}^3$, in this case embedded 1-dimensional surgery is the so-called \textit {knot surgery}. It is worth mentioning that there are infinitely many knot types and that 1-dimensional surgery on a knot may change the knot type or even result in a two-component link. A good introductory book on knot theory is \cite{Ad} among many other.

Looking back to the process of DNA recombination which exhibits embedded 1-dimensional 0-surgery, a DNA knot is the self-entanglement of a single circular DNA molecule. With the help of certain enzymes, site-specific recombination can transform supercoiled circular DNA into a knot or link.  The first electron microscope picture of knotted DNA was presented in \cite{WaDuCo}. In this experimental study, we see how genetically engineered circular DNA molecules can form DNA knots and links through the action of a certain recombination enzyme. A similar picture is presented in Fig~\ref{Fig8}, where site-specific recombination of a DNA molecule produces the Hopf link.
 
\begin{figure}[!h]
\begin{center}
\centering
\captionsetup{justification=centering}
\includegraphics[width=9cm]{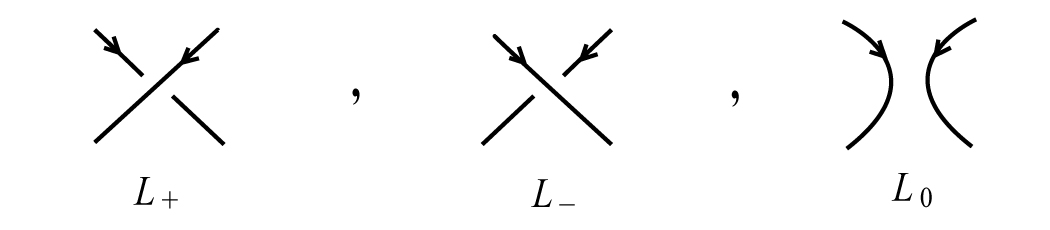}
\caption{One can pass from one of these three links to another via knot surgery.} 
\label{Fig31}
\end{center}
\end{figure}

Another theoretical example of knot surgery comprises the knot or link diagrams involved in the \textit{skein relations} satisfied by {\bf knot polynomials}, such as the Jones polynomial \cite{Jones2} and the Kauffman bracket polynomial \cite{K1}. For example, the illustration in Fig~\ref{Fig31} represents  a so-called `Conway triple', that is, three knot or link diagrams  $L_+$,$L_-$ and $L_0$ which are identical everywhere except in the region of a crossing and the polynomials of these three links satisfy a given linear relation.

\begin{rem}\label{Seifert} \rm
In analogy to embedded 1-dimensional 0-surgery, we also have the notion of embedded solid 1-dimensional 0-surgery. As $S^1$ is the boundary of $D^2$, any knot is the boundary of a, so-called, Seifert surface, so embedded solid 1-dimensional 0-surgery could be extended to a Seifert surface.   
\end{rem}

\subsection{Embedded solid 2-dimensional 0-surgery and related phenomena} \label{Embedded2D0}

In Section~\ref{ConSurg} we showed how we can pass from the two-ball description to the two-tori description of $S^3$. Although we had not yet defined it at that point, the process we described is, of course, an embedded solid 2-dimensional 0-surgery in $S^3$ on an initial manifold $M=D^3$. It is worth mentioning that all natural processes undergoing embedded solid 2-dimensional 0-surgery on an initial manifold $M=D^3$ can be also viewed in this context. For example, if one looks at the formation of black holes and examines it as an independent event in space, this process shows a decompactified view of the passage from a two 3-ball description of $S^3$, that is, the core of the star and the surrounding space, to a two torus description, that is, the accretion disc surrounding the black hole (shown in white in the third instance of Fig~\ref{Fig19}) and the surrounding space. In this Section, we will see how some natural phenomena undergoing solid 2-dimensional 0-surgery exhibit the causal or consequential relation to the ambient space mentioned in Section \ref{es3sp2} and are therefore better described by considering them as embedded in $S^3$.

For example, during the formation of {\bf tornados}, recall Fig~\ref{Fig17} (a), the process of solid 2-dimensional 0-surgery is triggered by the difference in the conditions of the lower and upper atmosphere. Although the air cycle lies in the complement space of the initial manifold $M=D^3 \times S^0$, it is involved in the creation of funnel-shaped clouds that will join the two spherical neighborhood (one in the cloud and one in the ground). Therefore \textit{the cause of the phenomenon extends beyond its initial manifold and surgery is the outcome of global changes}.

We will now discuss phenomena where \textit{the outcome of the surgery process propagates beyond the final manifold}. A first example are {\bf waterspouts}. After their formation, the tornado's cylindrical `cork', that is, the solid cylinder homeomorphic to the product set $D^1 \times D^2$, has altered the whole surface of the sea (recall Fig~\ref{Fig17} (b)). In other words, the spiral pattern on the water surface extends beyond the initial spherical neighborhood of the sea, which is represented by one of the two 3-balls of the initial manifold. 

As another example, during the formation of {\bf black holes}, the strong gravitational forces have altered the space surrounding the initial star and the singularity is created outside the final solid torus. In all these phenomena, the process of surgery alters matter outside the manifold in which it occurs. In other words, the effect of the forces causing surgery propagates to the complement space, thus causing a more global change in 3-space.

\begin{rem}\label{duality2d0} \rm
Looking back at Fig~\ref{Fig30}, it is worth pinning down the following duality of embedded solid 2-dimensional 0-surgery for $M=D^3$: the attraction of two points lying on the boundary of segment $L$ by the center of $D^3$ can be equivalently viewed in the complement space as the repulsion of these points by the  center of $B^3$ (that is, the point at infinity) on the boundary of curve (or line, if viewed in ${\mathbb R}^3$) $l-L$.
\end{rem}

\subsection{Embedded solid 2-dimensional 1-surgery and related phenomena} \label{Embedded2D1}
 
\begin{figure}[!h]
\begin{center}
\centering
\captionsetup{justification=centering}
\includegraphics[width=13.5cm]{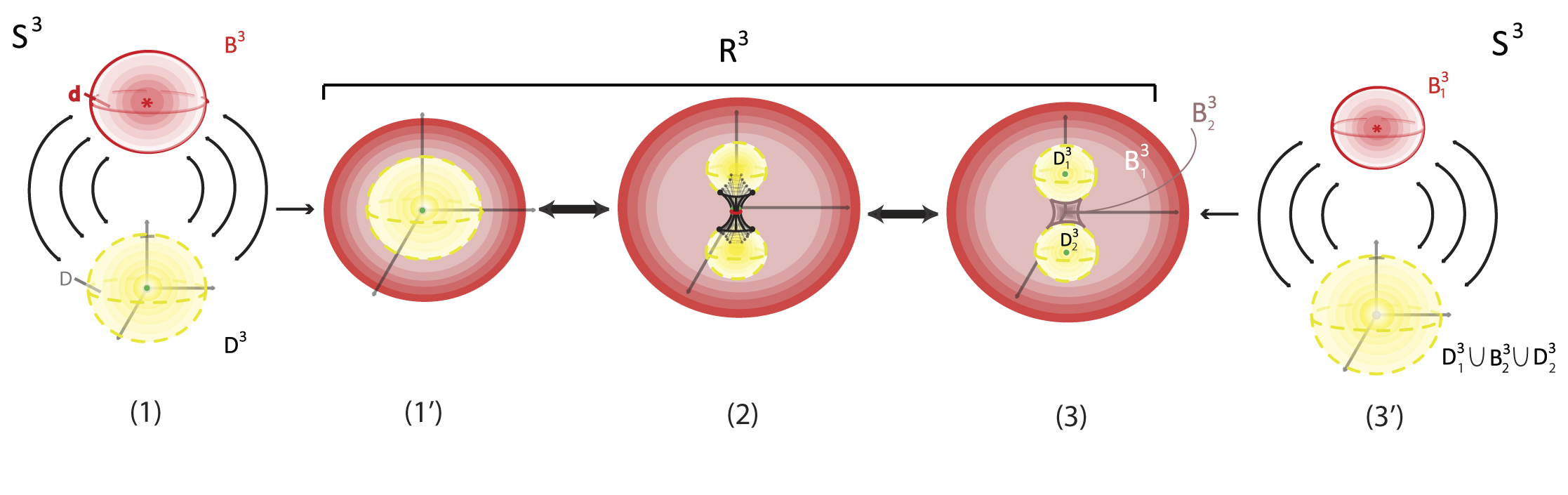}
\caption{{\bf Embedded solid 2-dimensional 1-surgery}}
\label{Fig32}
\end{center}
\end{figure}

We will now discuss the process of embedded solid 2-dimensional 1-surgery in $S^3$ in the same way we did for the embedded solid 2-dimensional 0-surgery in $S^3$, recall Fig~\ref{Fig30}. Taking again $M=D^3$ as the initial manifold, embedded solid 2-dimensional 1-surgery is illustrated in Fig~\ref{Fig32}. The process begins with disc $D$ in the 3-ball $D^3$ on which colinear attracting forces act, see instances (1) and (1$'$) for the decompactified view. In (3), the initial 3-ball $D^3$ is split in two new 3-balls  $D_1^3$ and $D_2^3$. By Definition~\ref{continuum2D}, the point at the origin (in green) evolves into the two centers of $D_1^3$ and $D_2^3$ (in green). This is exactly the same process as in Fig~\ref{Fig16} if we restrict it to $D^3$, but since we are in $S^3$, the complement space $B^3$ is also undergoing, by symmetry, solid 2-dimensional 1-surgery. Again, this is a double surgery with one attracting center which is inside the first 3-ball (in yellow) and outside the second 3-ball (in red). This process squeezes the central disc $D$ of $D^3$ while the central disc $d$ of $B^3$ engulfs disc $D$ and becomes the separating plane $d \cup D$.

As seen in instance (3) of Fig~\ref{Fig32}, the process alters the existing complement space $B^3$ to $B_1^3$ and creates a new space $B_2^3$ which can be considered as the 'void' between $D_1^3$ and $D_2^3$. By viewing the process in this way, we pass from a two 3-balls description of $S^3$ to another one, that is, from $S^3 = B^3 \cup D^3$ to $S^3 = (D_1^3 \cup B_2^3 \cup D_2^3) \cup B_1^3$.

\begin{rem}\label{duality2d1} \rm
The duality described in 2-dimensional 0-surgery is also present in 2-dimensional 1-surgery. Namely, the attracting forces from the circular boundary of the central disc $D$ to the center of $D^3$ can be equivalently viewed in the complement space as repelling forces from the center of $B^3$ (that is, the point at infinity) to the boundary of the central disc $d$, which coincides with the boundary of $D$. 
\end{rem}

All natural phenomena undergoing embedded solid 2-dimensional 1-surgery take place in the ambient 3-space. However, we do not have many examples of such phenomena which demonstrate the causal or consequential effects discussed in Section~\ref{es3sp2}. Yet one could, for example,  imagine taking a solid material specimen that has started necking and immerse it in some liquid  until its pressure causes fracture to the specimen. In this case the complement space is the liquid and it triggers the process of surgery. 

Finally, the annihilation of Falaco solitons is also a case of embedded solid 2-dimensional 1-surgery. The topological thread can be cut by many factors but in all cases these are related to the complement space.

\section{A dynamical system modeling embedded solid 2-dimensional 0-surgery}\label{ds}

So far, inspired by natural processes undergoing surgery, we have extended the formal definition of topological surgery by introducing new notions such as forces, solid surgery and embedded surgery. However, in our schematic models, time and dynamics were not introduced by equations. In this section we connect topological surgery, enhanced with these notions, with a dynamical system. We will see that, with a small change in parameters, the trajectories of its solutions are performing embedded solid 2-dimensional 0-surgery. Therefore, this dynamical system constitutes a specific set of equations modeling natural phenomena undergoing embedded solid 2-dimensional 0-surgery. More specifically, we will see that the change of parameters of the system affects the eigenvectors and induces a flow along a segment joining two steady state points. This segment corresponds to the segment $L$ introduced in Section~\ref{S3}  and the induced flow represents the attracting forces shown in Fig~\ref{Fig16} (a). Finally, we will see how our topological definition of solid 2-dimensional 0-surgery presented in Section~\ref{solid2D} is verified by our numerical simulations and, in particular, that surgery on a steady point becomes a limit cycle.

\hl{In} \cite{SS3}, \hl{the reader can find a more detailed analysis of the connection between the dynamical system and embedded solid 2-dimensional 0-surgery, with more numerical simulations and figures, and with a special focus on natural phenomena that can be modeled through this system, with emphasis on tornado formation.}

\subsection{The dynamical system and its steady state points} \label{DS}
In \cite{SaGr1}, N.Samardzija and L.Greller study the behavior of the following dynamical system ($\Sigma$) that generalizes the classical Lotka--Volterra problem \cite{Lo,Vo} into three dimensions: 

\begin{equation} \tag{\texorpdfstring{$\Sigma$}}
\left\{
\begin{array}{l}
\frac{dX}{dt}=X-XY+CX^2-AZX^2 \\
\\
\frac{dY}{dt}=-Y+XY \\
\\
\frac{dZ}{dt}=-BZ+AZX^2 \\
\end{array}
\right\} \ A, B, C > 0
\end{equation} 
\bigbreak

\noindent In subsequent work \cite{SaGr2}, the authors present a slightly different model, provide additional numerical simulations and deepen the qualitative analysis done in \cite{SaGr1}. Since both models coincide in the parametric region we are interested in, we will use the original model and notation and will briefly present some key features of the analyses done in \cite{SaGr1} and \cite{SaGr2}. 

The system ($\Sigma$) is a two-predator and one-prey model, where the predators $Y,Z$ do not interact directly with one another but compete for prey~$X$. As $X,Y,Z$ are populations, only the positive solutions are considered in this analysis. It is worth mentioning that, apart from a population model, ($\Sigma$) may also serve as a biological model and a chemical model, for more details see \cite{SaGr1}. 

The parameters $A,B,C$ are analyzed in order to determine the bifurcation properties of the system, that is, to study the changes in the qualitative or topological structure of the family of differential equations ($\Sigma$). As parameters $A,B,C$ affect the dynamics of constituents  $X,Y,Z$, the authors were able to determine conditions for which the ecosystem of the three species results in steady, periodic or chaotic behavior. More precisely, the authors derive five steady state solutions for the system but only the three positive ones are taken into consideration. These points are:

\begin{equation*}
S_{1}=\left(
\begin{array}{c}
0 \\
0 \\
0
\end{array}
\right) , \ \
S_{2}=\left(
\begin{array}{c}
1 \\
1+C \\
0
\end{array}
\right) , \ \
S_{3}=\left(
\begin{array}{c}
\sqrt{B/A} \\
0 \\
\frac{1+C\sqrt{B/A}}{\sqrt{AB}}
\end{array}
\right) 
\end{equation*}

It is worth reminding here that a steady state (or singular) point of a dynamical system is a solution that does not change with time. 

\subsection{Local behavior and numerical simulations} \label{LBAndSimu}
Let, now, $J(S_{i})$ be the Jacobian of ($\Sigma$) evaluated at  $S_{i}$ for $i=1,2,3$ and let the sets $\Gamma\{J(S_{i})\}$ and $W\{J(S_{i})\}$ to be, respectively, the eigenvalues and the corresponding associated eigenvectors of $J(S_{i})$. These are as follows: 

\begin{equation*}
\Gamma\{J(S_{1})\}=\{1,-1,-B \}
 ; \ \
W\{J(S_{1})\}=\left\{
\left[\begin{array}{c}
1 \\
0 \\
0
\end{array}\right],
\left[\begin{array}{c}
0 \\
1 \\
0
\end{array}\right],
\left[\begin{array}{c}
0 \\
0 \\
1
\end{array}\right]
\right\}
\end{equation*}

\begin{equation*}
\Gamma\{J(S_{2})\}=\{A-B,(C+\sqrt{(C-2)^2-8})/2,(C-\sqrt{(C-2)^2-8})/2 \}
\end{equation*}

\begin{equation*}
W\{J(S_{2})\}=\left\{
\left[\begin{array}{c}
1 \\
(C+1)/(A-B) \\
\frac{B+C-A+(C+1)/(B-A)}{A}
\end{array}\right],
\left[\begin{array}{c}
1 \\
\frac{C-\sqrt{(C-2)^2-8}}{2} \\
0
\end{array}\right],
\left[\begin{array}{c}
1 \\
\frac{C+\sqrt{(C-2)^2-8}}{2} \\
0
\end{array}\right]
\right\}
\end{equation*}
 
\begin{adjustwidth}{-0.14in}{0in}
\begin{equation*}
\Gamma\{J(S_{3})\}= \left\{\sqrt{\frac{B}{A}}-1,\frac{-1+\sqrt{1-8B(1+C\sqrt{B/A})}}{2},\frac{-1-\sqrt{1-8B(1+C\sqrt{B/A})}}{2} \right\}
\end{equation*}
\end{adjustwidth}

\begin{adjustwidth}{-1in}{0in}
\begin{equation*}
W\{J(S_{3})\}=\left\{
\left[\begin{array}{c}
1 \\
 -1-\frac{2\sqrt{AB}(1+C\sqrt{B/A})}{\sqrt{B/A}-1}\\
\frac{2(1+C\sqrt{B/A})}{\sqrt{B/A}-1}
\end{array}\right],
\left[\begin{array}{c}
1 \\
0 \\
\frac{-1-\sqrt{1-8B(1+C\sqrt{B/A})}}{2B}
\end{array}\right],
\left[\begin{array}{c}
1 \\
0 \\
\frac{-1+\sqrt{1-8B(1+C\sqrt{B/A})}}{2B}
\end{array}\right]
\right\}
\end{equation*}
\end{adjustwidth}

\bigbreak
\bigbreak

Using the sets of eigenvalues and eigenvectors presented above, the authors characterize in \cite{SaGr1}, \cite{SaGr2} the local behavior of the dynamical system around these three points using the Hartman-Grobman (or linearization) Theorem. Since $1>0$ and $-1,-B<0$, $S_{1}$ is a saddle point for all values of parameters $A,B,C$. However, the behavior around $S_{2}$ and $S_{3}$ changes as parameters $A,B,C$ are varied. The authors show that the various stability conditions can be determined by only two parameters: $C$ and $B/A$. It is also shown in \cite{SaGr1} that stable solutions are generated left of and including the line $B/A=1$ while chaotic/periodic regions appear on the right of the line $B/A=1$. We are interested in the behavior of ($\Sigma$) as it passes from stable to chaotic/periodic regions. Therefore we will focus and analyze the local behavior around $S_{2}$ and $S_{3}$ and present numerical simulations for: stable region (a) where $B/A=1$ and  $(1/8B-1)\sqrt{A/B}<C\le 2(1+\sqrt2)$ and chaotic/periodic region (b)  where $B/A>1$ and $(1/8B-1)\sqrt{A/B}<C\le 2(1+\sqrt2)$.

\bigbreak

{\noindent} $\bullet$ {\bf Region (a)} 
\smallbreak

\begin{figure}[!h]
\begin{center}
\centering
\captionsetup{justification=centering}
\includegraphics[width=13.5cm]{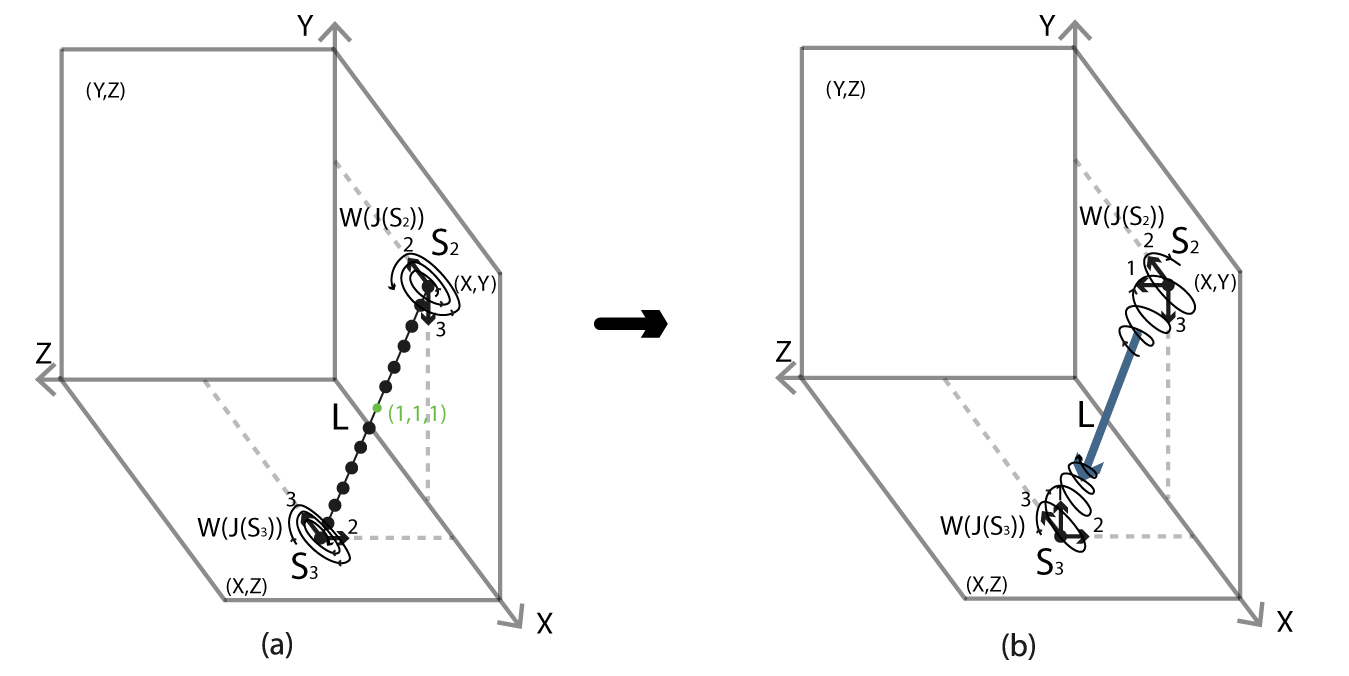}
\caption{{\bf Local behavior.} Flow induced along $L$ by changing parameter space from (a) $B/A=1$ to (b) $B/A>1$. Indices 1,2 and 3 indicate the first, second and third component in $W(J(S_{2}))$ and $W(J(S_{3}))$.}
\label{Fig33}
\end{center}
\end{figure}

\begin{figure}[!h]
\begin{center}
\centering
\captionsetup{justification=centering}
\includegraphics[width=13.5cm]{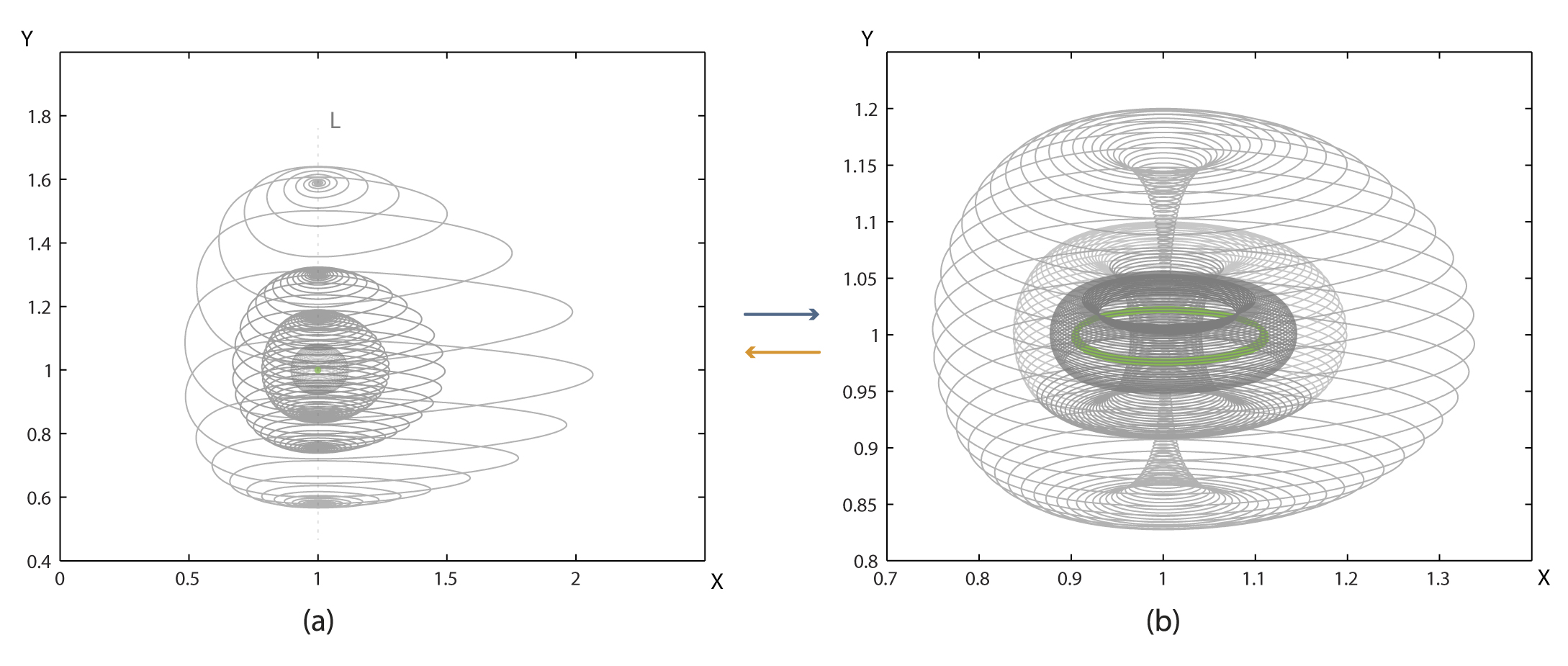}
\caption{{\bf Embedded solid 2-dimensional 0-surgery by changing parameter space from (a) $B/A=1$ to (b) $B/A>1$.}}
\label{Fig34}
\end{center}
\end{figure}

{\noindent} Setting $B/A=1$ and equating the right side of ($\Sigma$) to zero, one finds as solution the one-dimensional singular manifold: 
\[L=\{(X,Y,Z); \ X=1, Z=(1+C-Y)/A \}\]
that passes through the points $S_{2}$ and $S_{3}$. Since all points on $L$ are steady state points, there is no motion along it. For $(1/8B-1)\sqrt{A/B}<C\le 2(1+\sqrt2)$, $S_{2}$ {\it is an unstable center} while $S_{3}$ {\it is a stable center} (for a complete analysis of all parametric regions see \cite{SaGr1}). This means that if $\lambda_1, \lambda_2, \lambda_3$ denote the eigenvalues of either $S_{2}$ or $S_{3}$ with $\lambda_1\in \Bbb{R}$ and $\lambda_2, \lambda_3\in \Bbb{C}$, then $\lambda_1=0$ and  $Re(\lambda_2)=Re(\lambda_3)>0$ for $S_{2}$ while $\lambda_1=0$ and  $Re(\lambda_2)=Re(\lambda_3)<0$ for $S_{3}$. Moreover, the point $(X,Y,Z)=(1,1,C/A)$ is the center of $L$. The line segment $X=1$, $0<Y<1$ and $(1+C)/A<Z<C/A$ supports attracting type singularities (and includes $S_{3}$) while the line segment defined by  $X=1$, $1<Y<1+C$ and $0<Z<C/A$ supports unstable singularities (and includes $S_{2}$), for details see \cite{SaGr2}. More precisely, each attracting point corresponds to an antipodal repelling point, the only exception being the center of $L$ which can be viewed as the spheroid of 0-diameter. The local behavior of ($\Sigma$) around $S_{2}$ and $S_{3}$ in this region together with line $L$ are shown in Fig~\ref{Fig33} (a). A trajectory (or solution)  initiated near $L$ in the repelling segment expands until it gets trapped by the attracting segment, forming the upper and lower hemisphere of a distinct sphere. Hence, a nest of spherical shells surrounding line $L$ is formed, see Fig \ref{Fig34} (a). Moreover, the nest fills the entire positive space with stable solutions.

\bigbreak
{\noindent} $\bullet$ {\bf Region (b)}
\smallbreak

{\noindent} For $B/A>1$ and $(1/8B-1)\sqrt{A/B}<C\le 2(1+\sqrt2)$, $S_{2}$ {\it is an inward unstable vortex} and $S_{3}$ {\it is an outward stable vortex}. This means that in both cases they must satisfy the conditions $\lambda_1\in \Bbb{R}$ and $\lambda_2, \lambda_3\in \Bbb{C}$ with $\lambda_3= \lambda_2^*$, the conjugate of $\lambda_2$. The eigenvalues of $S_{2}$ must further satisfy $\lambda_1<0$ and  $Re(\lambda_2)=Re(\lambda_3)>0$, while the eigenvalues of $S_{3}$ must further satisfy $\lambda_1>0$ and $Re(\lambda_2)=Re(\lambda_3)<0$. The local behaviors around $S_{2}$ and $S_{3}$ for this parametric region are shown in Fig~\ref{Fig33} (b). It is worth mentioning that Fig~\ref{Fig33} (b) reproduces Fig 1 of \cite{SaGr1} with a change of the axes so that the local behaviors of $S_{2}$ and $S_{3}$ visually correspond to the local behaviors of the trajectories in Fig~\ref{Fig34} (b) around the north and the south pole. 

Note now that the point $S_{2}$ as well as the eigenvectors corresponding to its two complex eigenvalues, all lie in the $xy$--plane. On the other hand, the point $S_{3}$ and also the  eigenvectors corresponding to its two complex eigenvalues all lie in the $xz$--plane. The flow along line $L$ produced by the actions of these eigenvectors forces trajectories initiated near 
$S_{2}$ to wrap around $L$ and move toward $S_{3}$ in a motion reminiscent of hole drilling. The connecting manifold $L$ is also called the `slow manifold' in \cite{SaGr1} due to the fact that trajectories move slower when passing near it. As trajectories reach $S_{3}$, the eigenvector corresponding to the real eigenvalue of $S_{3}$ breaks out of the $xz$--plane and redirects the flow toward $S_{2}$. As shown in Fig~\ref{Fig34} (a) and (b), as $B/A=1$ moves to $B/A>1$, this process transforms each spherical shell to a toroidal shell. The solutions scroll down the toroidal surfaces until a limit cycle (shown in green in Fig~\ref{Fig34} (b)) is reached. It is worth pointing out that this limit cycle is a torus of 0-diameter and corresponds to the sphere of 0-diameter, namely, the central steady point of $L$ also shown in green in Fig~\ref{Fig34} (a). 

However, as the authors elaborate in \cite{SaGr2}, while for $B/A=1$ the entire positive space is filled with nested spheres, when $B/A>1$, only spheres up to a certain volume become tori. More specifically, quoting the authors: ``to preserve uniqueness of solutions, the connections through the slow manifold $L$ are made in a way that higher volume shells require slower, or higher resolution, trajectories within the bundle". As they further explain, to connect all shells through $L$, ($\Sigma$) would need to possess an infinite resolution. As this is never the case, the solutions evolving on shells of higher volume are `choked' by the slow manifold. This generates solution indetermination, which forces higher volume shells to rapidly collapse or dissipate. The behavior stabilizes when trajectories enter the region where the choking becomes weak and weak chaos appears. As shown in both \cite{SaGr1} and \cite{SaGr2}, the outermost shell of the toroidal nesting is a fractal torus. Note that in Fig~\ref{Fig34} (b) we do not show the fractal torus because we are interested in the interior of the fractal torus which supports a topology stratified with toroidal surfaces. Hence, all trajectories are deliberately initiated in its interior where no chaos is present. 

It is worth pointing out that Fig~\ref{Fig34} reproduces the numerical simulations done in \cite{SaGr2}. More precisely, Fig~\ref{Fig34} (a) represents solutions of ($\Sigma$) for $A=B=C=3$ and trajectories initiated at points $[1,1.59,0.81]$, $[1,1.3,0.89]$, $[1,1.18;0.95]$, $[1,1.08,0.98]$ and $[1,1,1]$. Fig~\ref{Fig34} (b) represents solutions of ($\Sigma$) for $A=2.9851,B=C=3$ and trajectories initiated at points $[1.1075,1,1]$, $[1,1,0.95]$, $[1,1,0.9]$ and $[1,1,1]$. 

As already mentioned, as $B/A=1$ changes to $B/A>1$, $S_{2}$ changes from an unstable center to an inward unstable vortex and $S_{3}$ changes from a stable center to an outward stable vortex. It is worth reminding that this change in local behavior is true not only for the specific parametrical region simulated in Fig~\ref{Fig34}, but applies to all cases satisfying $(1/8B-1)\sqrt{A/B}<C\le 2(1+\sqrt2)$. For details we refer the reader to Tables II and III in \cite{SaGr1} that recapitulate the extensive diagrammatic analysis done therein.

Finally, it is worth observing the changing of the local behavior around $S_{2}$ and $S_{3}$ in our numerical simulations. In Fig~\ref{Fig34} (a), for $B/A=1$ we have:
\begin{adjustwidth}{-1.6in}{0in}
\begin{equation*}
\Gamma\{J(S_{2})\}=\{{\bf0.0000},1.500-1.3229i,1.500+1.3229i \} , 
\Gamma\{J(S_{3})\}= \{{\bf0.0000},-1.000+4.8780i,-1.000-4.878i \}
\end{equation*}  
\end{adjustwidth}
  while in Fig~\ref{Fig34} (b), for $B/A>1$, both centers change to vortices (inward unstable and outward stable) through the birth of the first eigenvalue shown in bold (negative and positive respectively): 

\begin{adjustwidth}{-1.6in}{0in}
\begin{equation*}
\Gamma\{J(S_{2})\}=\{{\bf-0.0149},1.500-1.3229i,1.500+1.3229i \} , 
\Gamma\{J(S_{3})\}= \{{\bf0.0025},-1.000+4.8780i,-1.000-4.878i \}
\end{equation*}
\end{adjustwidth}

\begin{rem}\label{NumMethods} \rm
The use of different numerical methods may affect the shape of the attractor. For example, as mentioned in \cite{SaGr2}, higher resolution produces a larger fractal torus and a finer connecting manifold. However, the ‘hole drilling’ process and the creation of a toroidal nesting is always a common feature.
\end{rem}

\subsection{Connecting the dynamical system with embedded solid 2-dimensional 0-surgery} \label{DSConn}
 
In this section, we will focus on the process of embedded solid 2-dimensional 0-surgery on a 3-ball $D^3$ viewed as a continuum of concentric spheres together with their common center: $D^3=\cup_{0<r\leq 1} S^2_r \cup \{P\}$. Recall from Section~\ref{solid2D} that the process is defined as the union of 2-dimensional 0-surgeries on the whole continuum of concentric spheres $S^2_r$ and on the limit point $P$. For each spherical layer, the process starts with attracting forces acting between $S^0 \times D^2$, i.e two points, or poles, centers of two  discs. In natural phenomena undergoing solid 2-dimensional 0-surgery, such as tornadoes (recall Fig~\ref{Fig17}) or Falaco solitons (recall Fig~\ref{Fig18}), these forces often induce a helicoidal motion from one pole to the other along the line $L$ joining them.

Having presented the dynamical system ($\Sigma$) in Section~\ref{DS} and its local behavior in  Section~\ref{LBAndSimu}, its connection with embedded solid 2-dimensional 0-surgery on a 3-ball is now straightforward. To be precise, surgery is performed on the manifold formed by the trajectories of ($\Sigma$). Indeed, as seen in Fig~\ref{Fig34} (a) and (b), with a slight perturbation of parameters, trajectories pass from spherical to toroidal shape through a `hole drilling' process along a slow manifold $L$ which pierces all concentric spheres. The spherical and toroidal nestings in Figs~\ref{Fig16} (a) and ~\ref{Fig34} are analogous. The attracting forces acting between the two poles shown in blue in the first instance of Fig~\ref{Fig16} (a) are realized by the flow along $L$ (also shown in blue in Fig~\ref{Fig33} (b)). When $B/A>1$, the action of the eigenvectors is an attracting force between $S_{2}$ and $S_{3}$ acting along $L$, which drills each spherical shell and transforms it to a toroidal shell. 

Furthermore, in order to introduce solid 2-dimensional 0-surgery on $D^3$ as a new topological notion, we had to define that 2-dimensional 0-surgery on a point is the creation of a circle. The same behavior is seen in ($\Sigma$). Namely, surgery on the limit point $P$, which is a steady state point, creates the limit cycle which is the limit of the tori nesting. As mentioned in \cite{SaGr2}, this type of bifurcation is a `Hopf bifurcation', so we can say that we see surgery creating a Hopf bifurcation.

Hence, instead of viewing surgery as an abstract topological process, we may now view it as a property of a dynamical system. Moreover, natural phenomena exhibiting 2-dimensional topological surgery through a `hole-drilling' process, such as the creation of Falaco solitons, the formation of tornadoes, of whirls, of wormholes, etc, may be modeled mathematically by the dynamical system ($\Sigma$). For more details, see also \cite{SS2}. \hl{This system enhances the schematic topological model presented in Fig~{\ref{Fig16}~(a)} with analytical formulation of the underlying dynamics. Indeed, if we link the three time-dependent quantities $X,Y,Z$ to physical parameters of related phenomena undergoing 2-dimensional 0-surgery,  system ($\Sigma$) can provide time forecasts for these phenomena.}

\begin{rem}\label{DSANDS3} \rm
It is worth pointing out that ($\Sigma$) is also connected with the 3-sphere $S^3$. We can view the spherical nesting of Fig~\ref{Fig34} (a) as the 3-ball $D^3$ shown in Fig~\ref{Fig30} (1) and (1$^\prime$). Surgery on its central point creates the limit cycle which is the core curve $c$ of $V_1$ shown in Fig~\ref{Fig30} (3) and (3$^\prime$). If we extend the spherical shells of Fig~\ref{Fig34} to all of ${\mathbb R}^3$ and assume that the entire nest resolves to a toroidal nest, then the slow manifold $L$ becomes the infinite line $l$. In the two-ball description of $S^3$, $l$ pierces all spheres, recall Fig~\ref{Fig30} (1$^\prime$), while in the two-tori description, it is the core curve of $V_2$ or the `untouched' limit circle of all tori, recall Fig~\ref{Fig30} (3) and (3$^\prime$).
\end{rem} 

\begin{rem}\label{Kiehn} \rm
In \cite{Ki} R.M. Kiehn studies how the Navier-Stokes equations admit bifurcations to Falaco solitons. In other words, the author looks at another dynamical system modeling this natural phenomenon which, as we showed in Section \ref{2D0}, exhibits solid 2-dimensional 0-surgery. To quote the author: ``It is a system that is globally stabilized by the presence of the connecting 1-dimensional string" and ``The result is extraordinary for it demonstrates a global stabilization is possible for a system with one contracting direction and two expanding directions coupled with rotation". It is also worth quoting Langford \cite{Lang} which states that computer simulations indicate that ``the trajectories can be confined internally to a sphere-like surface, and that Falaco Soliton minimal surfaces are visually formed at the North and South pole". One possible future research direction would be to investigate the similarities between this system and ($\Sigma$) in relation to surgery. 
\end{rem}

\section{Conclusions}
Topological surgery occurs in numerous natural phenomena of various scales where a sphere of dimension 0 or 1 is selected and attracting forces are applied. Examples of such phenomena comprise: chromosomal crossover, magnetic reconnection, mitosis, gene transfer, the creation of Falaco solitons, the formation of whirls and tornadoes, magnetic fields and the formation of black holes. 

In this paper we explained these natural processes via topological surgery. To do this we first enhanced the usual static description of topological surgery  of dimensions~1 and 2 by introducing dynamics, by means of attracting forces. We then filled in the interior spaces in 1- and 2-dimensional surgery, introducing the notions of solid 1- and 2-dimensional surgery. This way more natural phenomena can fit with these topologies. Further, we introduced the notion of embedded surgery, which leaves room for the initial manifold to assume a more complicated configuration and \hl{describes how} the complementary space of the initial manifold participates in the process. \textit{Thus, instead of considering surgery as a formal and static topological process, it can now be viewed as an intrinsic and dynamic property of many natural phenomena.} 

\hl{Apart from the examples studied in this paper, there are several other phenomena exhibiting surgery, and our topological models indicate where to look for the forces causing surgery and what deformations should be observed in the local submanifolds involved. Also, our modeling of the changes occurring in the complement space during embedded surgery provides a} `\hl{global}’ \hl{explanation of the phenomenon, which can be of great physical significance. Similarly, our descriptions of the duality of forces in embedded surgery could potentially lead to new physical explanations. For instance, it would interesting to investigate from the physical point of view, whether the forces collapsing a star to a black hole could be equally viewed as repelling forces from the} `\hl{point at infinity}’.

Equally important, all these new notions resulted in pinning down the connection of solid 2-dimensional 0-surgery with a dynamical system. This connection gives us on the one hand {\it a mathematical model for 2-dimensional surgery} and, on the other hand, {\it a dynamical system modeling natural  phenomena exhibiting 2-dimensional topological surgery through a `hole-drilling' process}. \hl{The provided dynamical system presents significant common features with our schematic topological model of 2-dimensional 0-surgery, in the sense that eigenvectors act as the attracting forces, trajectories lie on the boundaries of the manifolds undergoing surgery and surgery on the steady state point (which is the central point of the spherical nesting) creates a limit cycle (which is central circle of the toroidal nesting). Furthermore, this system enhances our modeling as it can provide time forecasts for these phenomena.} 

\hl{This subject together with the search of other dynamical systems realizing other types of surgery will be the subject of future work.} Another possible future research direction includes using the proposed dynamical system as a base for establishing a more general theoretical connection between topological surgery and bifurcation theory.
Currently we are working with Louis H.Kauffman on generalizing the notions presented in this paper to 3-dimensional surgery and higher dimensional natural processes. \hl{A first step toward this generalization can be found in} \cite{SS2}.

We hope that through this study, topology and dynamics of natural phenomena, as well as  topological surgery itself, will be better understood and that our connections will serve as ground for many more insightful observations.





\section*{Acknowledgments}
We are grateful to Louis H.Kauffman and Cameron McA.Gordon for many fruitful conversations on Morse theory and 3-dimensional surgery. We would also like to thank Colin C.Adams for fruitful communications and to acknowledge early discussions with Nick Samardzija on topological aspects of dynamical systems as well as a comment by Tim Cochran pointing out the connection of our new notions with Morse theory. 
\hl{Finally, we would like to thank the Reviewers for their positive comments and for helping us through their suggestions to open the manuscript to broader scientific audiences.}


\newpage
\appendix
\section{Appendix: Mathematical notions} \label{Appendix}
{\noindent \bf Topological spaces}  
\begin{enumerate}
  \item A {\it topological space} is a set $X$ with a distinguished family $\tau$ of subsets possessing the following properties:
\begin{itemize}
  \item the empty set and the whole set $X$ belong to $\tau$
  \item the intersection of a finite number of elements of $\tau$ belongs to $\tau$
  \item the union of any subfamily of elements of $\tau$ belongs to $\tau$  
\end{itemize}
The family $\tau$ is said to be the {\it topology} on $X$. Any set belonging to $\tau$ is called {\it open}. A {\it neighborhood}  of a point $x \in X$ is any open set containing $x$. Any set whose complement is open is called {\it closed}. The minimal closed set (with respect to  inclusion) containing a given set $A \subset X$ is called the 
$closure$ of $A$ and is denoted by $\bar{A}$. The maximal open set contained in a given set $A \subset X$ is called the {\it interior} of $A$ and is denoted by $Int(A)$. 
\item If $(X,\tau)$ is a topological space, a {\it base} of the space $X$ is a subfamily $\tau' \subset \tau$ such that any element of $\tau$ can be represented as the union of elements of $\tau'$. In other words, $\tau'$ is a family of open sets such that any open set of $X$ can be represented as the union of sets from this family. In the case when at least one base of $X$ is countable, we say that $X$ is a space with {\it countable base}.
\item If $X \times Y$ is the Cartesian product of the topological spaces $X$ and $Y$ (regarded as sets), then $X \times Y$ becomes a topological space (called the {\it product} of the spaces $X$ and $Y$) if we declare open all the products of open sets in $X$ and in $Y$ and all possible unions of these products.
\end{enumerate}

{\noindent \bf Manifolds}
\begin{enumerate}
\setcounter{enumi}{3}
\item A topological space is said to be a {\it Hausdorff space} if any two distinct points of the space have nonintersecting neighborhoods.
\item A Hausdorff space $M^{n}$ with countable base is said to be an {\it n-dimensional topological manifold} if any point $x \in  M^{n}$ has a neighborhood homeomorphic to $\mathbb{R} ^{n}$ or to $\mathbb{R} ^{n}_+$, where $\mathbb{R} ^{n}_+ = \{(x_1,...,x_n)\mid x_i \in  \mathbb{R}, x_1 \geq 0 \}$. For example, a surface is a 2-dimensional manifold.
\item The set of all points $x \in  {M} ^{n}$ that have no neigbourhoods homeomorphic to $\mathbb{R} ^{n}$ is called the {\it boundary} of the manifold ${M} ^{n}$ and is denoted by $\partial {M} ^{n}$. When  $\partial {M} ^{n} = \emptyset$, we say that ${M} ^{n}$ is a {\it manifold without boundary}. It is easy to verify that if the boundary of a manifold ${M} ^{n}$ is nonempty, then it is an $(n-1)$-dimensional manifold.
\end{enumerate}

{\noindent \bf \hl{Properties of manifolds}}  
\begin{enumerate}
  \setcounter{enumi}{6}
\item  \hl{A topological space is called \textit{connected} if it cannot be presented as a union of two nonintersecting nonempty sets each of which is simultaneously open and closed.}
\item  \hl{A topological space $X$ is called \textit{compact} if any open covering of $X$ (i.e any collection of open sets of $X$ whose union is $X$) has a finite subcovering.}
\item \hl{An atlas $\{(U_{\alpha},{\varphi}_\alpha)\}$ of a smooth manifold is called orienting if the Jacobians of all the maps  ${\varphi}_\beta \circ {{\varphi}_\alpha}^{-1}$, where $U_{\beta} \cap U_{\alpha}\neq \emptyset$, are positive. A manifold possessing an orienting atlas is called \textit{orientable}. If an orienting atlas of the manifold $M$ is chosen, we say that an orientation is given on $M$ (for details see for example} \cite{PS}) .
\end{enumerate}

{\noindent \bf The standard topology of $\mathbb{R} ^{n}$ and its one-point compactification}  
 \begin{enumerate}
 \setcounter{enumi}{9}
\item In general, to define the topology $\tau$, it suffices to indicate a base of the space. For the space $\mathbb{R} ^{n} = \{(x_1,...,x_n)\mid x_i \in  \mathbb{R} \}$, the standard topology is given by the base $U_{a,\epsilon}=\{x \in  \mathbb{R} ^{n} \mid  \lvert x-a \lvert < \epsilon \}$, where $a \in \mathbb{R} ^{n}$ and $\epsilon > 0$. We can additionally require that all the coordinates of the point $a$, as well as the number $\epsilon$, be rational; in this case we obtain a countable base.
\item  To the set $\mathbb{R} ^{n}$ let us add the element $\infty$ and introduce in $\mathbb{R} ^{n}  \cup \{ \infty\} $ the topology whose base is the base of $\mathbb{R} ^{n}$ to which we have added the family of sets $U_{\infty,R}=\{x \in  \mathbb{R} ^{n} \mid  \lvert x \lvert > R \} \cup \{\infty\}$. The topological space thus obtained is called the {\it one-point compactification} of $\mathbb{R} ^{n}$; it can be shown that this space is homeomorphic to the $n-$dimensional sphere $S^n=\{x \in  \mathbb{R} ^{n+1} \mid  \lvert x \lvert =1\}$.
\end{enumerate}

{\noindent \bf Homeomorphisms and gluings}  
\begin{enumerate}
  \setcounter{enumi}{11}
 \item The map of one topological space into another is called {\it continuous} if the preimage of any open set is open. A map $f: X \rightarrow Y$ is said to be a {\it homeomorphism} if it is bijective and both $f$ and $f^{-1}$ are continuous; the spaces $X$ and $Y$ are then called {\it homeomorphic} or {\it topologically equivalent}.
 \item An injective continuous map between topological spaces $f: X \hookrightarrow Y$ is called an {\it embedding} if $f$ is an homeomorphism between $X$ and $f(X)$.
 \item Suppose $X$ and $Y$ are topological spaces without common elements, $A$ is a subset of $X$, and $f: X \rightarrow Y$ is a continuous map. In the set $X \cup Y$, let us introduce the relation $a \sim  f(a)$. The resulting quotient space $(X \cup Y)/\sim$ is denoted by $X \cup_f Y$; the procedure of constructing this space is called {\it gluing} or {\it attaching} $Y$ to $X$ along the map $f$.
 \end{enumerate}

The above definitions were taken from \cite{PS}. For more details, the reader is referred to \cite{PS, Ro, Pr}.

\nolinenumbers

%
%
 



\newpage
\singlespacing


\begin{thebibliography}{00}

\bibitem{An} Antoniou S. The chaotic attractor of a 3-dimensional Lotka--Volterra dynamical system and its relation to the process of topological surgery. Diplom Thesis, National Technical University of Athens. 2005. 

\bibitem{SS} Lambropoulou S, Antoniou S. Dynamical systems and topological surgery; 2008. Preprint. Available from: \url{http://arxiv.org/pdf/0812.2367v1.pdf}.


\bibitem{SSNI} Lambropoulou S. Topological Surgery and Dynamics (joint work with Samardzija N, Diamantis I, Antoniou S).  Mathematisches Forschungsinstitut Oberwolfach Report, Workshop: Algebraic Structures in Low-Dimensional Topology. 2014; No. 26/2014. Available from: \url{https://www.mfo.de/document/1422a/OWR_2014_26.pdf}.


\bibitem{SaGr1} Samardzija N, Greller L. Explosive route to chaos through a fractal torus in a generalized Lotka-Volterra model. Bull Math Biol 50. 1988; No. 5: 465--491. doi: 10.1007/BF02458847.

\bibitem{PS} Prasolov VV, Sossinsky AB. Knots, links, braids and 3-manifolds. AMS Translations of Mathematical Monographs 154; 1997.

\bibitem{Ro} Rolfsen D. Knots and links. Publish or Perish Inc. AMS Chelsea Publishing; 2003.

\bibitem{Ra} Ranicki A. Algebraic and Geometric Surgery. Oxford Mathematical Monographs, Clarendon Press; 2002.

\bibitem{Su} Sumners D. Lifting the Curtain: Using Topology to Probe the Hidden Action of Enzymes. Not Am Math Soc 42. 1995; No. 5: 528--537. Available from: \url{http://www.ams.org/notices/199505/sumners.pdf}.

\bibitem{DaAn} Dahlburg RB, Antiochos SK. Reconnection of antiparallel magnetic flux tubes. J Geophys Res 100. 1995; No. A9: 16991--16998. doi: 10.1029/95JA01613.

\bibitem{LaRiSu} Laing CE, Ricca RL, Sumners D. Conservation of writhe helicity under anti-parallel reconnection. Scientific Reports 5. 2014; No. 9224. doi: 10.1038/srep09224. Available from: \url{http://www.nature.com/srep/2015/150328/srep09224/full/srep09224.html}. 

\bibitem{Mil} Milnor J. Morse Theory. Princeton University Press; 1963. 



\bibitem{KoFe} Kondrashov D, Feynman J, Liewer PC, Ruzmaikin A. Three-dimensional Magnetohydrodynamic Simulations of the Interaction of Magnetic Flux Tubes. The Astrophysical Journal 519. 1999; 884--898. doi: 10.1086/307383.

\bibitem{Pu} Pujari S. Useful Notes on the Mechanism of Crossing Over. Your Article Library; 2015. Available from: \url{http://www.yourarticlelibrary.com/biology/useful-notes-on-the-mechanism-of-crossing-over-biology-810-words/6634/}.

\bibitem{JoLe} Johnson AB, Lewis J et al. Molecular Biology of the Cell. Garland Science; 2002. Available from: \url{http://www.ncbi.nlm.nih.gov/books/NBK26845/}.

\bibitem{Ke} Kerr RM. Fully developed hydrodynamic turbulence from a chain reaction of reconnection events. Procedia IUTAM. 2013; No. 9: 57-–68. doi: 10.1016/j.piutam.2013.09.006. Available from: \url{http://www.sciencedirect.com/science/article/pii/S2210983813001284}. 

\bibitem{Ki} Kiehn RM. Non-equilibrium systems and irreversible processes - Adventures in applied topology vol. 1 - Non equilibrium thermodynamics. University of Houston Copyright CSDC Inc; 2004. pp. 147,150.


\bibitem{Ki2} Kiehn RM. Falaco Solitons, Cosmic Strings in a Swimming Pool; 2001. Preprint. Available from: \url{http://arxiv.org/pdf/gr-qc/0101098v1}.

\bibitem{Ott} Ott CD, et al. Dynamics and Gravitational Wave Signature of Collapsar Formation. Phys Rev Lett 106. 2011. doi: 10.1103/PhysRevLett.106.161103.

\bibitem{HHGRSV} Hartwell LH, Hood L, Goldberg ML, Reynolds AE, Silver LM, Veres RC. Genetics, from genes to genomes. McGraw Hill; 2000. pp. 486.


\bibitem{UV} Zhigilei LV, Lecture notes, Chapter 8: Failure - University of Virginia, Department of Materials Science; 2010. Available from: \url{http://people.virginia.edu/~lz2n/mse209/Chapter8.pdf}.

\bibitem{KeFa} Keeton WT, McFadden CH. Elements of biological science. W.W. Norton \& Company Inc; 1983. pp. 395.

\bibitem{SS2} Antoniou S, Lambropoulou S. Topological Surgery in Nature. Algebraic Modeling of Topological and Computational Structures and Applications, book at Springer Proceedings in Mathematics and Statistics (PROMS), S. Lambropoulou, D. Theodorou, P. Stefaneas and L. H. Kauffman Eds. Forthcoming.

\bibitem{Glatzmaier} Glatzmaier GA, Roberts PH. A three-dimensional self-consistent computer simulation of a geomagnetic field reversal, Nature 377 (6546). 1995; 203–-209. doi:10.1038/377203a0.

\bibitem{NilesJ} Johnson N, A visualization of the Hopf fibration. Available from: \url{http://nilesjohnson.net/hopf.html}.

\bibitem{Ad} Adams C. Why Knot?: An Introduction to the Mathematical Theory of Knots with Tangle. Key College Publishing; 2004.


\bibitem{WaDuCo} Wasserman SA, Dungan JM, Cozzarelli NR. Discovery of a predicted DNA knot substantiates a model for site-specific recombination. Science 229. 1985; 171--174.

\bibitem{Jones2} Jones VFR. A new knot polynomial and von Neumann algebras,  Not Am Math Soc 33. 1986; No. 2: 219--225.  

\bibitem{K1} Kauffman LH. State models and the Jones polynomial, Topology 26. 1987; No. 3:  395--407. doi: 10.1016/0040-9383(87)90009-7. Available from: \url{www.maths.ed.ac.uk/~aar/papers/kauffmanjones.pdf}.

\bibitem{SS3} Lambropoulou S, Antoniou S. Topological Surgery,  Dynamics and Applications to Natural Processes.  J. Knot Theory Ramifications 26. 2017. doi: 10.1142/S0218216517430027.

\bibitem{Lo} Lotka AJ. Undamped Oscillations Derived from the Law of Mass Action, J. American Chem. Soc. 1920;  No. 42: 1595--1599.

\bibitem{Vo} Volterra V. Le\c{c}ons sur la Th\'{e}orie Math\'{e}matique de la lutte pour la vie, Paris, Gauthier-Villars. Gabay, J.; ed. 1931, reissued 1990. doi: 10.1090/S0002-9904-1936-06292-0.

\bibitem{SaGr2} Samardzija N, Greller L. Nested tori in a 3-variable mass action model. Proc R Soc Lond A Math Phys Sci 439. 1992; No. 1907: 637--647. doi: 10.1098/rspa.1992.0173.

\bibitem{Lang} Langford LD. A review of interactions of Hopf and steady-state bifurcations, Non-linear dynamics and turbulence, G.I. Barrenblatt, G.Iooss, D.D Joseph, Pitman, Boston. 1983; 215-237.

\bibitem{Pr} Prasolov VV. Intuitive Topology. American Mathematical Society, Mathematical World vol. 4; 1995.

\end{thebibliography}
\end{document}